# Five lectures on regularity structures and SPDEs

I. BAILLEUL

**Abstract**

The following set of five lectures provides an introduction to regularity structures and their use for the study of singular stochastic partial differential equations. Two appendices provide some additional informations that enter in the main text either as some technical results or as some results that deepen the context within which we set these lectures.

*Please do not hesitate to send me your comments.*[1] This document was typed with TeXmacs.

# Lecture 1 – Singular stochastic PDEs

Throughout we denote by $(\Omega, \mathscr{F}, \mathbb{P})$ the probability space on which all the stochastic objects are defined, and write $\omega$ for a generic element of $\Omega$. We denote by $\mathbb{T}^d$ the $d$-dimensional torus.

## 1 What are *singular* stochastic PDEs?

**§1. Product operation** – For $r < 0$ we define for the moment the distribution space[2]

$$C^r([0,1]) \equiv \{g = \partial^{[|r|]+1} h; h \in C^{1-\{|r|\}}([0,1])\}$$

with $\|g\|_r$ defined as the $(1-\{|r|\})$-Hölder norm of $h$. Here is a fact from analysis.

**Theorem (Bony)** – *For any $r_1, r_2 \in \mathbb{R}$, the map*

$$(u_1, u_2) \in C^{r_1}([0,1]) \times C^{r_2}([0,1]) \mapsto u_1 u_2 \in C^{\min(r_1, r_2)}([0,1])$$

*is well-defined and continuous iff $r_1 + r_2 > 0$.*

This statement has a higher dimensional analogue that involves the same condition $r_1 + r_2 > 0$, independently of the dimension of the ambiant space and independently of its scaling property. (You can find a proof of this statement in Appendix 1.) As you will see, this fundamental fact is the starting point of our story. Let us illustrate that with an elementary example.

**§2. Controlled differential equations** – Let $(e_j)_{j \geq 1}$ be an orthonormal basis of $L^2([0,1])$ and $(\gamma^j)_{j \geq 1}$ be a sequence of independent normal random variables. One defines the real-valued *Brownian motion* over $[0,1]$ by the formula $B_t = \sum_{j \geq 1} \gamma^j \int_0^t e_j$. One can prove that $B$ is almost surely continuous.[3] Better, one has $B \in C^{1/2-\epsilon}([0,1])$ for all $\epsilon > 0$ and $B \notin C^{1/2}([0,1])$ almost surely. An $\ell$-dimensional Brownian motion is a family $(B^1, \ldots, B^\ell)$ of $\ell$ independent Brownian motions.

If you have had a course on stochastic calculus, you know that it is possible to make sense of the stochastic differential equations in $\mathbb{R}^d$ driven by some smooth bounded vector fields $f = (f_1, \ldots, f_\ell)$ and an $\ell$-dimensional Brownian motion

$$dx_t = f(x_t) dB_t = \sum_{i=1}^\ell f_i(x_t) dB_t^i \tag{1}$$

---

1. Univ Brest, CNRS, LMBA - UMR 6205, F-29238 Brest, France. *E-mail*: ismael.bailleul@univ-brest.fr
2. We denote by $\{a\}$ the fractional part of a positive real number $a$ and by $[a]$ its integer part.
3. While the path $B \in C([0,1], \mathbb{R})$ depends on the choice of orthonormal basis $(e_j)_{j \geq 1}$ its law turns out to be independent of that choice.

and solve it uniquely. How is it possible while $dB \in C^{-1/2-\epsilon}$ and $x \in C^{1/2-\epsilon}$, so $(-1/2-\epsilon) + (1/2-\epsilon) < 0$ and $f(x)dB$ is not well-defined pathwise, from Bony's theorem?! The point is that this equation is actually formulated in a purely probabilistic way via an integral formulation that involves a stochastic integral: a purely probabilistic object whose definition involves martingales and requires no regularity conditions on the integrand and the integrator. Rather it involves some integrability condition with respect to the measure $\mathbb{P} \times$ Leb. It is only after the stochastic integration map has been defined that we can formulate the equation as a fixed point problem on a space of random variables that involves this stochastic integration map.

Can we still say something if $B$ in the above equation is not Brownian motion but another sort of random path? A fractional Brownian motion with Hurst index $H < 1/2$ for instance. Can we say anything if $B \in C^{1/2-\epsilon}$ is not random but deterministic?

**§3. Classical stochastic partial differential equations** – Let $(e_i)_{i \geq 1}$ stand for an orthonormal basis of $L^2(\mathbb{T}^d)$ and $(B^1, B^2, \ldots)$ stand for a sequence of independent Brownian motions. Classical stochastic partial differential equations over the $d$-dimensional torus $\mathbb{T}^d$ involve the notion of cylindrical Brownian motion $B_t = \sum_{i \geq 1} B_t^i e_i$. Almost surely this is not an element of $L^2(\mathbb{T}^d)$, but $dB$ is almost surely a distribution on $[0,1] \times \mathbb{T}^d$, called *spacetime white noise*.

A typical example of classical parabolic SPDE reads informally

$$(\partial_t - D^2)u = f(u)dB + g(u, \partial u) \qquad (2)$$

with $D^2$ the Laplace-Beltrami operator[4], and some initial condition is given. It has a well-defined integral formulation involving the stochastic integral operator

$$(\phi \text{ adapted random variable}) \mapsto \int_0^t \phi_s(dB_s) = \sum_{i \geq 1} \int_0^t \phi_s(e_i) dB_s^i$$

with respect to the cylindrical Brownian motion $B$, in the particular case where $\phi$ takes the form $\phi_s(e) = e^{(t-s)D^2}(f(u_s)e)$ for some adapted process $u$. Here again, the purely probabilistic notion of martingale plays a fundamental role in the definition of this integration notion, and in the definition of a solution to Equation (2). This type of definition can be extended to the larger class of semimartingales.

Here are a two classes of examples of classical SPDEs.

– Stochastic perturbations of reaction-diffusion equations where $f(u) = 1$. The equation of stochastic quantization in Euclidean quantum field theory corresponds to $f(u) = 1$ and $g$ the functional gradient of an energy functional. In neurosciences, Fitz-Nagumo's equations describes the activity of nerve cells and correspond to reaction-diffusion equations with polynomial nonlinearities $g$.

– In population genetics the solution of Equation (2), with $g(u, \partial u) = au - b$ and $f(u) = \sqrt{u_+}$ or $f(u) = \sqrt{u_+(1-u)_+}$, models some large scale population structure over space and time.

**§4. Singular stochastic PDEs** – Not all equations driven by a spacetime white noise make sense in the classical SPDE/Itô setting. This is for instance the case of the KPZ equation on $[0,1] \times \mathbb{T}$

$$(\partial_t - D^2)u = -(u')^2 + dB.$$

Since $u$ is expected to be of the form $(\partial_t - D^2)^{-1}(dB) + v$, with $v$ of better regularity than $(\partial_t - D^2)^{-1}(dB)$, while this term has parabolic regularity $1/2 - \epsilon$ for all $\epsilon > 0$, the term $u'$ has regularity $-1/2 - \epsilon$ and the term $(u')^2$ is ill-defined, from Bony's theorem.

---

[4]. We do not use the classical notation $\Delta$ for the Laplace-Beltrami operator as this symbol will stand later for a different object.

Moreover spacetime white noise is not an appropriate model of noise in a large number of situations. What can we do for instance if $dB$ in (2) is not a spacetime white noise but only a distribution-valued Gaussian process with some covariance structure? Is there any hope to do anything if $dB$ is not even a Gaussian random distribution? We asked the same question in §2 about the controlled differential equation (1). One can answer these questions with the tools of regularity structures[5]. This setting was introduced by M. Hairer and some co-authors in four groundbreaking works [11,6,7,5].

An *equation* is said to be *singular* if it involves a priori an ill-defined product $u_1 u_2$ with $u_1 \in C^{r_1}$ and $u_2 \in C^{r_2}$ and $r_1 + r_2 \leq 0$. Singular PDEs are meaningless from a classical point of view. Spacetime white noise has almost surely parabolic regularity $-3/2 - \epsilon$ for all $\epsilon > 0$. All classical SPDEs are then singular from a pathwise point of view as soon as the noise is multiplied by some function of $u$.

**Examples of singular stochastic PDEs** – Here are three typical examples of singular stochastic PDEs.

– The 2-dimensional parabolic Anderson model equation

$$(\partial_t - D^2)u = u \xi_x^{(2)} \qquad (3)$$

is set on the 2-dimensional torus $\mathbb{T}^2$, with unknown a real-valued function $u$, with $\xi_x^{(2)}$ a space white noise of spatial Besov-Hölder regularity $-1 - \eta$, for all $\eta > 0$ almost surely. Its solution models the propagation of heat in an impure medium. Note that the noise here is time-independent, so it is definitely outside of the realm of stochastic calculus, which only deals with semimartingales.

– The $\Phi_3^4$ equation from stochastic quantization

$$(\partial_t - D^2)u = -u^3 + \xi_{t,x}^{(3)}$$

is set on the 3-dimensional torus, with unknown a real-valued function $u$, with $\xi_{t,x}^{(3)}$ a *spacetime* white noise of parabolic Besov-Hölder regularity $-5/2 - \eta$, for all $\eta > 0$ almost surely. It models, amongst others, the local average magnetization of a field of magnet in a particular scaling limit.

– The (KPZ) equation from random interface growth models

$$(\partial_t - D^2)u = -(u')^2 + \xi_{t,x}^{(1)}$$

is set on the 1-dimensional torus, with unknown a real-valued function $u$, with $\xi_{t,x}^{(1)}$ a *spacetime* white noise of parabolic Besov-Hölder regularity $-3/2 - \eta$, for all $\eta > 0$ almost surely.

**§5. What we can learn from controlled differential equations.** The model situation of controlled differential equations (1) already contains the fundamental product problem that characterizes singular equations when $B$ is a *deterministic* control of Hölder regularity strictly smaller than $1/2$. This model situation was investigated in depth in some groundbreaking work of T. Lyons in the mid-nineties. You will find in Appendix 2 a light introduction to the setting of rough paths that Lyons introduced.[6] Hairer's regularity structures have their roots in Lyons's theory of rough paths.

We can learn a number of things from this simpler setting.

A solution path to a controlled ordinary differential equation

$$z_t = z_0 + \sum_{i=1}^{\ell} \int_0^t f_i(z_r) \, \dot{h}_r^i \, dr$$

---

5. It is also possible to use some other tools than regularity structures.
6. It is not necessary to know anything about rough paths to follow the present lectures.

driven by a *deterministic control* $h \in C^1([0,1], \mathbb{R}^\ell)$ has a local expansion property

$$z_t = z_s + f(z_s)(h_t - h_s) + (f'_j f_i)(z_s)\left\{\int_s^t (h_r - h_s)^i (h'_r)^j\, dr\right\} + O(|t-s|^3).$$

around a arbitrary time $s$ obtained by Picard iteration. The control $h$ enters in this formula via the coefficients $(h_t - h_s, (\int_s^t (h_r - h_s)^i \otimes (h'_r)^j\, dr_{i,j})_{1 \leq i,j \leq \ell}) \in \mathbb{R}^\ell \times (\mathbb{R}^\ell \times \mathbb{R}^\ell)$ of the expansion. This type of local expansion actually characterizes a solution path.

There is a simple and versatile idea behind this type of characterization, encoded in a statement called the *sewing lemma*. *Given a family $(\mu_{ts})_{0 \leq s \leq t \leq 1}$ of vectors in $\mathbb{R}^d$, there is a simple coherence relation on the elements of this family that ensures the existence of a unique path whose increments between the times $s$ and $t$ are (quantatively) close to $\mu_{ts}$ for $t$ near $s$, for all $s$.*

The iterated integrals $\int_s^t (h_r - h_s)^i (h'_r)^j\, dr$ do not make sense if $h$ is only a *deterministic* $\alpha$-Hölder, with $\alpha < 1/2$. Let us say here that $1/3 < \alpha < 1/2$. A rough path $\boldsymbol{X}$ is a placeholder $(X_t - X_s, \mathbb{X}_{ts})_{0 \leq s \leq t \leq 1}$ for the non-existing family of coefficients $(h_t - h_s, (\int_s^t (h_r - h_s)^i\, dh_r^j)_{1 \leq i,j \leq \ell})_{0 \leq s \leq t \leq 1}$, that satisfies some size and algebraic relations that mimick those satisfied by the coefficients when $h \in C^1$.

Here is the road to make sense of the controlled ordinary differential equation

$$z_t = z_0 + \int_0^t f(z_r) d\boldsymbol{X}_r$$

driven by a deterministic rough path $\boldsymbol{X}$. The main problem consists in making sense of the integral; then we can formulate the equation as a fixed point problem.

(1) We consider only some paths $z$ that have some local expansion property of the form

$$z_t \simeq z_s + z'_s(X_t - X_s)$$

for some $z'_s$ (not the classical derivative).

(2) Given such a path, we use the sewing lemma to see the existence of a unique path whose increments are close to

$$\mu_{ts} = f(z_s)(X_t - X_s) + (f'f)(z_s)\mathbb{X}_{ts}.$$

This path is our definition of the integration map $\int_0^\cdot f(z_r) d\boldsymbol{X}_r$.

This sketch is expanded in Appendix 2. We keep from that short discussion three fundamental ideas.

– One can define uniquely some objects from some families of local and approximate descriptions.

– This leads us to describe some objects by their local expansion properties.

– If ever some coefficients in some expansion do not make sense when some side object (here a control) is not regular enough, take a placeholder for it!

## 2 A blackbox for a class of singular stochastic PDEs

We consider here for simplicity some problems set in the spatial domain $\mathbb{T}^d$ with real-valued unknown $u: [0,T] \times \mathbb{T}^d \to \mathbb{R}$, on a time interval to be fixed with the solution, and assume that the random noise $\xi$ on $\mathbb{R} \times \mathbb{T}^d$ is scalar valued. We consider some parabolic equation of the form

$$(\partial_t - D^2)u = f(u)\xi + g(u, \partial u) \tag{4}$$

with initial condition $u_0$ in some function space. A naive application of Schauder estimate tells us that if the right hand side of Equation (4) has the same parabolic regularity $-\beta_0 < 0$ as the noise $\xi$ then $u$ needs to have parabolic regularity $2 - \beta_0$. The product $f(u)\xi$ is then meaningful only if $2 - 2\beta_0 > 0$, i.e. $0 < \beta_0 < 1$, by a generalization of Theorem 1 to a parabolic setting. Spacetime white noise has almost surely parabolic regularity $-3/2 - \epsilon$ for all $\epsilon > 0$. All classical SPDEs are then singular from a pathwise point of view.[7] What saves us in the classical setting is that we can define a solution as a fixed point of a integral formulation of the equation that involves a stochastic integral. The latter does not involve any regularity condition for its definition but rather some integrability condition with respect to $\mathbb{P} \times \text{Leb}$.

We assume from now on that $-\beta_0 \leq -1$, so Equation (4) is singular and a priori meaningless pathwise. A mollification procedure by convolution $\xi \mapsto \xi^\epsilon$ turns the noise $\xi$ into a smooth function $\xi^\epsilon$ on $\mathbb{R} \times \mathbb{T}^d$. A cheap way out for the study of (4) consists in replacing the noise $\xi$ by $\xi^\epsilon$, solve for the corresponding equation, and then send the regularization parameter $\epsilon > 0$ to 0. The following fundamental result tells us that is not exactly the right thing to do.

**Theorem 1.** *Under some mild conditions on the law of the random noise $\xi$ and the initial condition $u_0$, for any mollification procedure $\xi \mapsto \xi^\epsilon$ there exists a deterministic explicit function $c^\epsilon$ and a random positive time $T = T(\omega)$ such that the solutions to the classical and locally well-posed equation*

$$(\partial_t - D^2)u^\epsilon = f(u^\epsilon)\xi^\epsilon + g(u^\epsilon, \partial u^\epsilon) + c^\epsilon(u^\epsilon, \partial u^\epsilon) \tag{5}$$

*with initial condition $u_0$ converge in probability in $C([0,T] \times \mathbb{T}^d, \mathbb{R})$ to some limit function $u$ that does not depend on the mollification procedure.*

Equation (5) is called the *renormalized equation*; there is no counterterm ($c^\epsilon = 0$) if $\xi$ is almost surely continuous for instance. The semi-informal statement of Theorem 1 calls for a number of comments.

*1. The model case of stochastic differential equations* – The statement of Theorem 1 is reminiscent of a famous result by Wong & Zakai (1965) about stochastic differential equations. If one mollifies Brownian motion ($B \to B^\epsilon$) in an Itô equation $z_t = z_0 + \int_0^t f(z_s) dB_s$, with $f \in C_b^2$, and one solves instead $z_t^\epsilon = z_0 + \int_0^t (f(z_s^\epsilon) dB_s^\epsilon + c(z_s^\epsilon) ds)$, for some appropriate choice of non-diverging (!) counter-term $c(\cdot)$, then the paths $z^\epsilon$ are converging in probability to the solution path to the Itô equation.

*2. About the mild conditions* – Several sets of conditions can be imposed on the law of the noise. The main condition in Chandra & Hairer [7] involves the cumulants of $\xi$. Hairer & Steele [14] and Bailleul & Hoshino [3] require that the law of $\xi$ satisfies a spectral gap inequality. None of these conditions is stronger than the other. Spacetime white noise satisfies both conditions, as well as a whole class of Gaussian and non-Gaussian random noises.

One of the unspecified conditions in Theorem 1 involves both the noise and the equation via a notion of *subcriticality*. It essentially means here that no term in the right hand side of (2) has a formal[8] regularity smaller than the regularity of the noise.[9]

The fact that Theorem 1 applies beyond the realm of classical stochastic calculus is a major point.

We typically ask the initial condition $u_0$ to be $r$-Hölder regular for some $r > 0$.

---

7. The word 'pathwise' means here 'for $\mathbb{P}$-almost all $\omega \in \Omega$'.

8. The regularity is computed with the rule of thumb on the regularity of a product and Schauder estimates, as if all the operations were well defined: $\text{reg}(a) \times \text{reg}(b)$ has regularity $b$ if $a > 0$ and $b < 0$, and regularity $a + b$ if both $a$ and $b$ are negative.

9. It will be encoded in Lecture 3 as the requirement that the index set $A \subset \mathbb{R}$ in the $T$ space of a regularity structure is bounded below and has no accumulation points.

**3. Non-uniqueness of the counterterm $c^\epsilon$** – The function $c^\epsilon$ is not unique. There is actually a representation $\mathsf{C}(\cdot)$ of a finite dimensional Lie group in a function space such that $c^\epsilon = \mathsf{C}(k^\epsilon)$ for some (typically diverging) element $k^\epsilon$ of the group. It turns out that for any other element $k$ in the group the solution $u^\epsilon_{(k)}$ to Equation (5) with $\mathsf{C}(k^\epsilon k)$ in place of $c^\epsilon = \mathsf{C}(k^\epsilon)$, with initial condition $u_0$, also converges in probability in $C([0, T_{(k)}] \times \mathbb{T}^d, \mathbb{R})$ to some limit function $u_{(k)}$, for some random positive time $T_{(k)} = T_{(k)}(\omega)$. The functions $u$ and $u_{(k)}$ have no reason to coincide.

The group element $k^\epsilon$ and the counterterm $c^\epsilon(\cdot)$ are not always diverging.

To single out an appropriate choice of $k^\epsilon$ was one of the main results of Bruned, Hairer & Zambotti in [6], Theorem 6.17 therein. The first description of $\mathsf{C}(\cdot)$ was given by Bruned, Chandra, Chevyrev & Hairer in Theorem 2.21 of [5]. A much simplified proof of this result was given by Bailleul & Bruned in [1].

**4. Link with classical SPDEs** – When $f$ and $g$ are sufficiently regular and bounded, when $\xi$ is spacetime white noise and $g(u, \partial u) = g(u)$, then Equation (4) is a classical SPDE with a classical probabilistic solution. If $\xi^\epsilon = \rho^\epsilon \star \xi$ is obtained by spacetime convolution of $\xi$ with $\rho^\epsilon(t, x) = \epsilon^{-3}\rho(\epsilon^{-2}t, \epsilon^{-1}x)$ for some smooth, even and compactly supported function $\rho$ with integral 1, then there is an explicit function $h_{f,g,\rho} = h + (\cdots)$ that depends on $f, g$ and $\rho$ such that the limit function $u$ in Theorem 1 coincides almost surely with the solution of the classical SPDE (3) with $h_{f,g,\rho}$ in place of $g$. This fact proved by Hairer & Pardoux [13] is the equivalent of the Wong-Zakai theorem in an SPDE setting

**5. New sets of equations** – The parabolic Anderson model equation (3), with its purely space noise, can be dealt within the settign of regularity structures, as well as the other examples of equations from §2 in Section 1. A whole class of discrete one dimensional evolution equations can be proved to have the KPZ equation as their scaling limit. A particular case of the generalized KPZ equation

$$(\partial_t - D^2)u = g(u)|u'|^2 + h(u)\nabla u + f(u)\xi$$

models the motion of an elastic string in a viscous random environment.

**6. About the tools for proving Theorem 1** – Consider a fixed singular SPDE. We can try a naive approach based on our understanding of the model situation of controlled ordinary differential equations. Assume for that purpose that we have an equivalent of the sewing lemma that tells us that some sufficiently accurate approximations of any object in the equation is sufficient to characterize it. (See point **2** *(ii)* below.) A formal Picard iteration of the equation produces some local expansions for a solution of the equation. Some of their terms do not make sense if the noise in the equation has low regularity. As in the model setting of control differential equations: do not panic! Take a placeholder for it and keep going on. Then try and formulate your equation as a fixed point problem in a space of objects that have the kind of local expansion that you just derived.

A complete proof of Theorem 1 requires quite some work. One looks for a reformulation of the classical and well-posed renormalized equation (5) that is robust enough to allow for $\epsilon > 0$ to go to 0. One can obtain such a reformulation within the language of regularity structures. One of the marvels of their application to the study of singular stochastic PDEs is that it provides a decoupling between the probabilistic and dynamical sides of the study of Equation (4). This decoupling takes the following form.

1. There is a metric product space $\mathsf{N} \equiv \mathsf{N}_1 \times \cdots \times \mathsf{N}_a$ called the space of *enhanced noises* such that $\xi$ takes almost surely its values in $\mathsf{N}_1$. The integer $a$ depends on the equation and it is typical that all the spaces $\mathsf{N}_i$ are distribution spaces. The space $\mathsf{N}$ is the space of placeholders in the above picture.

    – Any deterministic continuous function $\zeta$ on spacetime has a canonical lift[10] $\zeta^{\rm can+} \equiv (\zeta_1^{\rm can}, \ldots, \zeta_a^{\rm can})$ into the space $\mathsf{N}$ with $\zeta_1^{\rm can} = \zeta$. This naive lifting map is a polynomial smooth function of $\zeta \in C([0,1] \times \mathbb{T}^d)$.

---

10. Because their is no reason to introduce any placeholder for some well defined quantities.

– For a non-continuous noise $\xi$ one can construct *in a clever way* an **N**-valued random variable[11] $\xi^+ = (\xi_1, \ldots, \xi_\ell)$ such that $\xi_1 = \xi$ and $\xi^+$ is the limit in **N** *in a probabilistic sense* of some random variables of the form $\xi^{\epsilon+} \equiv (\xi_1^\epsilon, \ldots, \xi_a^\epsilon)$ with $\xi_1^\epsilon = \xi^\epsilon$. This construction is what *renormalization* is about; it was first achieved in [7]. *The construction of this enhanced noise is the only place in the analysis of Equation* (2) *where probability is involved.*

2. The other elements that are involved in the proof of Theorem 1 are *deterministic*. Given any enhanced noise $n \in \mathbf{N}$, one can associate to the formal equation (4) and any positive finite time $T$ a non-classical *n*-dependent Banach function space $\mathcal{D}(n)$ over $[0,T] \times \mathbb{T}^d$ that represent the space of objects that have some local expansions whose coefficients are given by the different components of *n*. This space has the following properties.

   *(i)* One can associate to any $\boldsymbol{v} \in \mathcal{D}(n)$ and any $z = (t, x) \in (0, T) \times \mathbb{T}^d$ in the state space a distribution on spacetime denoted by $\Lambda^n(\boldsymbol{v})(z)$. It is meant to be the local a priori expansion of a classical object near an arbitrary point $z$.

   *(ii)* There is indeed a unique distribution that is well approximated near $z$ by the local expansion $\Lambda^n(\boldsymbol{v})(z)$, for all $z$. It is called the *reconstruction of $\boldsymbol{v}$*, and denoted $\mathscr{R}(\Lambda^n(\boldsymbol{v}))$.

   *(iii)* One can formulate an *n*-dependent fixed point equation $\boldsymbol{v} = \Phi^n(\boldsymbol{v})$ in $\mathcal{D}(n)$ with a unique solution $\boldsymbol{u}^n$, provided we work on a sufficiently small time interval. The function $\boldsymbol{u}^n$ depends continuously on the parameter *n* and the reconstruction[12]

   $$u^n \equiv \mathscr{R}(\Lambda^n(\boldsymbol{u}^n)) \in C([0,T] \times \mathbb{T}^d)$$

   of $\boldsymbol{u}^n$ also depends continuously on the enhanced noise $n \in \mathbf{N}$.

   *(iv)* If we choose for *n* the canonical enhancement $\zeta^{\mathrm{can}+}$ of an arbitrary continuous function $\zeta$ then $u^\zeta \equiv u^{\zeta^{\mathrm{can}+}}$ turns out to be the solution of the equation

   $$(\partial_t - D^2)u^\zeta = f(u^\zeta)\zeta + g(u^\zeta, \partial u^\zeta)$$

   with its corresponding initial condition. The *clever way* of constructing $\xi^+$ in point **1** refers to the fact that the (non-canonically) enhanced regularized noise $\xi^{\epsilon+}$ has the property that $u^\epsilon \equiv u^{\xi^{\epsilon+}}$ is the solution of the renormalized equation

   $$(\partial_t - D^2)u^\epsilon = f(u^\epsilon)\xi^\epsilon + g(u^\epsilon, \partial u^\epsilon) + c^\epsilon(u^\epsilon, \partial u^\epsilon)$$

   with its corresponding initial condition. The continuity of the map $n \to u^n$ and the probabilistic convergence of $\xi^{\epsilon+}$ to $\xi^+$ in **N** then prove the probabilistic convergence of $u^\epsilon \equiv u^{\xi^{\epsilon+}}$ to some limit function $u^{\xi^+} = \mathscr{R}(\Lambda^{\xi^+}(\boldsymbol{u}^{\xi^+}))$: This is the content of Theorem 1.

## 3 Plan of the lectures

This series of lectures aims at giving you a clear picture of a number of objects that appear in the setting of regularity structures, and to describe the mechanics of the proof of Theorem 1. Giving a full proof is beyond the scope of these lectures and we refer for instance to the *Tourist's guide to regularity structures and singular stochastic PDES* [2], and the work [3] on renormalisation, for this.

   Lecture 2 – Germs and PDEs

   Lecture 3 – The mechanics of local expansions

   Lecture 4 – Renormalised continuous models

---

11. A placeholder for some a priori ill-defined quantities.
12. One has indeed $u^n \in C([0,T] \times \mathbb{T}^d)$ under the assumptions of Theorem 1.

Lecture 5 – Back to the future

**Where else can you learn about regularity structures and singular stochastic PDEs?**
In chronological order: Hairer's lecture notes [12], Chandra & Weber's introductory article [8], the book [10] by Friz & Hairer, and Berglund's book [4] provide some accessible accounts of part of the theory. The work [9] of Corwin & Shen provides a nice non-technical overview of a particular context in which singular stochastic PDEs arise.

## References

[1] I. Bailleul and Y. Bruned. *Locality for singular stochastic PDEs*. arXiv:2109.00399. To appear in Ann. Probab. (2025+). *A short proof of the renormalised equation for a large class of renormalisation procedures.*

[2] I. Bailleul and M. Hoshino. *A tourist's guide to regularity structures and singular stochastic PDEs*. EMS Surv. Math. Sci. (2025). *A detailed presentation of the analytic and algebraic sides of regularity structures, applied to the generalized KPZ equation.*

[3] I. Bailleul and M. Hoshino. *Random models on regularity-integrability structures*. arXiv:2310.10202, (2023). *Construction of a large class of renormalised random models encompassing the case of the BPHZ renormalisation, under a spectral gap assumption on the law of the noise. No Feynman graphs involved.*

[4] N. Berglund. *An introduction to singular stochastic PDEs – Allen-Cahn equations, metastability and regularity structures*. EMS Ser. Lect. Math., EMS Press, Berlin (2022). *A very nice book centered on the study of the Allen-Cahn equation in dimensions 1, 2 and 3, for which regularity structures are introduced.*

[5] Y. Bruned, A. Chandra, I. Chevyrev and M. Hairer. *Renormalising SPDEs in regularity structures*. J. Europ. Math. Soc. (JEMS) **23**(11) (2021), 869–947. *First proof of the general form of the renormalised equation for the BPHZ renormalisation procedure.*

[6] Y. Bruned, M. Hairer and L. Zambotti. *Algebraic renormalisation of regularity structures*. Invent. Math. **215**(3) (2019), 1039–1156. *Introduces a general algebraic setting for regularity structures. The notion of subcriticality is introduced here.*

[7] A. Chandra and M. Hairer. *An analytic BPHZ theorem for regularity structures*. arXiv:1612.08138v5 (2016). *The first proof of convergence of the BPHZ renormalised models, based on cumulant assumptions on the noise; proof done via a deep analysis of a whole class of Feynman graphs.*

[8] A. Chandra and H. Weber. *Stochastic PDEs, regularity structures and interacting particle systems*. Ann. Fac. Sci. Toulouse **26**(4) (2017), 847–909. *An early introduction to the subject.*

[9] I. Corwin and H. Shen. Some recent profess in singular stochastic partial differential equations. Bull. Amer. Math. Soc. **57**(3) (2020), 409–454. *Provides a number of motivating situations where singular SPDEs arise as scaling limit of microscopic systems.*

[10] P. Friz and M. Hairer. *A course on rough paths. With an introduction to regularity structures*. Universitext, Springer, Cham 2020. *A reference for learning the fundamentals of rough paths theory and regularity structures.*

[11] M. Hairer. *A theory of regularity structures*. Invent. Math. **198**(2) (2014), 269–504. *The groundbreaking work where regularity structures were introduced.*

[12] M. Hairer. *Introduction to regularity structures*. Braz. J. Probab. Stat. **29**(2) (2015), 175–210. *An early introduction to the theory of regularity structures.*

[13] M. Hairer and E. Pardoux. *A Wong-Zakai theorem for stochastic PDEs*. J. Math. Soc. Japan **67**(4) (2015), 1551–1604. *As the title says...*

[14] M. Hairer and R. Steele. *The BPHZ theorem for regularity structures via the spectral gap theorem*. Arch Rational Mech Anal. **248**(9) (2024). *A much more accessible (compared to [7]), and inductive, proof of convergence of the BPHZ renormalised models based on a spectral gap assumption. No Feynman graphs involved.*

# Lecture 2 – Germs and PDEs

One of the fundamental ideas of regularity structures consists in trading the objects "*functions and distributions*" $\Lambda$ for a family of $\mathbb{R}^d$-indexed local approximate descriptions $(\Lambda_x)_{x\in\mathbb{R}^d}$ of these objects. This comes under the form of the notion of *germ*. We need some norms to introduce them via an inequality that quantifies the approximation.

There are several equivalent ways of introducing the Hölder spaces of negative regularity. We choose the semigroup approach because it is convenient here. Write $(P_t)_{t>0}$ for the heat semigroup and note that the heat kernel $p_t$ satisfies for any multiindex $m = (m_1, \ldots, m_d) \in \mathbb{N}^d$ and $a \geq 0$ the integrability estimate

$$\int |\partial_x^m p_t(x,y)| \, |x-y|^a \, dy \lesssim t^{\frac{a-|m|}{2}} \tag{6}$$

where $|m| \equiv m_1 + \cdots + m_d$. For $r \leq 0$ we define $\mathcal{C}^r$ as the completion of $L^\infty$ with respect to the norm[13]

$$\|\Lambda\|_r \equiv \sup_{0 < t \leq 1} t^{-\frac{r}{2}} \|P_t(\Lambda)\|_\infty.$$

We have for all $\epsilon > 0$ and $r \leq 0$ the classical estimate

$$\|(P_t - \mathrm{Id})(\Lambda)\|_{r-\epsilon} \lesssim t^{\epsilon/2} \|\Lambda\|_r.$$

We denote by $\langle \Lambda, \varphi \rangle$ the result of applying a distribution $\Lambda$ on a test function $\varphi$.

For $m = (m_1, \ldots, m_d) \in \mathbb{N}^d$ and $z = (z_1, \ldots, z_d) \in \mathbb{R}^d$ we set $z^m \equiv \prod_{i=1}^d z_i^{m_i}$.

## 4 Germs and their reconstruction

**Definition 2.** *A **germ** is a bounded $\mathbb{R}^d$-indexed collection $(\Lambda_x)_{x\in\mathbb{R}^d}$ of elements of some fixed $\mathcal{C}^r$ space with $r \leq 0$. It is said to be **$\gamma$-coherent** ($\gamma \in \mathbb{R}$) if there is some exponents $\beta \leq \gamma$ such that*

$$|\langle \Lambda_y - \Lambda_x, p_t(x,\cdot) \rangle| \lesssim t^{\beta/2} \Big(|y-x| + t^{\frac{1}{2}}\Big)^{\gamma-\beta} \tag{7}$$

*for all $x, y \in \mathbb{R}^d$ and $0 < t \leq 1$.*

**Example.** Let $f \in C^r$ for a non-integer positive $r$. Write $F_x^\ell(z) = \sum_{|m| \leq r-|\ell|} \partial^{m+\ell} f(x) \frac{(z-x)^m}{m!}$ for the Taylor polynomial of $\partial^\ell f$ at $x$, and set $\partial^\ell f(z) = F_x^\ell(z) + R_x^\ell(z)$, with $|R_x^\ell(z)| \lesssim |z-x|^{r-|\ell|}$. One has from Newton's binomial identity $F_y(z) = F_x(z) + \sum_{|m| < r} R_x^m(y) \frac{|z-y|^m}{m!}$, so we see from (6) that the family $(F_x)_{x\in\mathbb{R}^d}$ satisfies the identity (7) with $\beta = 0$ and $\gamma = r$.

The following reconstruction theorem is due to Hairer [5, Theorem 3.10] in the setting of modelled distributions over a regularity structure. The regularity structures free form stated here is due to Caravenna & Zambotti [4, Theorem 5.1]. The proof given here is taken from [2].

**Theorem 3.** *One can associate to any $\gamma$-coherent germ $(\Lambda_x)_{x\in\mathbb{R}^d}$ with $\gamma > 0$ a unique distribution $\Lambda \in \mathcal{C}^r$ such that one has uniformly in $x \in \mathbb{R}^d$ and $0 < t \leq 1$*

$$|\langle \Lambda - \Lambda_x, p_t(x,\cdot) \rangle| \lesssim t^{\gamma/2}. \tag{8}$$

---

[13]. The space $L^\infty$ is continuously embedded in $\mathcal{C}^0$ and strictly smaller. The equivalence of this heat semigroup-based norm with the classical isotropic Littlewood-Paley-based norm on $\mathcal{B}_{\infty\infty}^r$ is presented for instance in Theorem 2.34 of Bahouri, Chemin & Danchin's textbook [1].

The distribution $\Lambda$ is called the *reconstruction of the germ* $(\Lambda_x)_{x\in\mathbb{R}^d}$; it is also denoted $\mathscr{R}(\Lambda_x)_{x\in\mathbb{R}^d}$.

**Proof.** *Existence* – Define the function

$$\mathbb{I}_s^t(x) \equiv \int p_{t-s}(x,y)\langle \Lambda_y, p_s(y,\cdot)\rangle\, dy.$$

We will construct $\Lambda$ as $\lim_{t\downarrow 0}\lim_{s\downarrow 0} \mathbb{I}_s^t$, with the limits taken in that order. We note that

$$|\mathbb{I}_t^t(x)| = |\langle \Lambda_x, p_t(x,\cdot)\rangle| \lesssim t^{r/2} \tag{9}$$

uniformly in $x\in\mathbb{R}^d$ and $0<t\leq 1$ because the family $(\Lambda_x)_{x\in\mathbb{R}^d}$ is bounded in $\mathcal{C}^r$. Pick $\epsilon>0$.

(i) *Convergence of* $\mathbb{I}_s^t(x)$ *as* $s\downarrow 0$ – For $0<s_2<s_1<t\leq 1$ we have from (6)

$$\begin{aligned}
|\mathbb{I}_{s_2}^t(x) - \mathbb{I}_{s_1}^t(x)| &= \left|\int p_{t-s_1}(x,z)p_{s_1-s_2}(z,y)\langle \Lambda_y - \Lambda_z, p_{s_2}(y,\cdot)\rangle\, dzdy\right| \\
&\leq \int p_{t-s_1}(x,z)p_{s_1-s_2}(z,y)\, s_2^{\beta/2}\bigl(|y-x|+s_2^{\frac{1}{2}}\bigr)^{\gamma-\beta}\, dzdy \\
&\lesssim \bigl((s_1-s_2)^{1/2}+s_2^{1/2}\bigr)^{\gamma-\beta} s_2^{\beta/2} \lesssim s_1^{\gamma/2}
\end{aligned} \tag{10}$$

uniformly in $x\in\mathbb{R}^d$. Hence the ($x$ and locally in $t$) uniform limit $\mathbb{I}_0^t(x)\equiv\lim_{s\downarrow 0}\mathbb{I}_s^t(x)\in\mathbb{R}$ exists *since the exponent $\gamma$ is positive*, and $|\mathbb{I}_0^t(x)|\lesssim t^{r/2}$, uniformly in $x\in\mathbb{R}^d$ and $0<t\leq 1$, from (9) and (10).

(ii) *Convergence of* $\mathbb{I}_0^t$ *in* $\mathcal{C}^{r-\epsilon}$ *as* $t\downarrow 0$. As the identity

$$P_{t'}(\mathbb{I}_0^t) = \mathbb{I}_0^{t+t'}$$

follows from the semigroup property it follows from the definition of the space $\mathcal{C}^r$ and the bounds (9) and (10) that the family $(\mathbb{I}_0^t)_{0<t\leq 1}$ is bounded in $\mathcal{C}^r$, since $r\leq 0$ and $\gamma>0$. We further have for all $0<t_2<t_1\leq 1$

$$\|\mathbb{I}_0^{t_1} - \mathbb{I}_0^{t_2}\|_{r-\epsilon} = \|(P_{t_1-t_2} - \mathrm{Id})(\mathbb{I}_0^{t_2})\|_{r-\epsilon} \lesssim (t_1-t_2)^{\epsilon/2}\|\mathbb{I}_0^{t_2}\|_r \lesssim (t_1-t_2)^{\epsilon/2}.$$

The family $(\mathbb{I}_0^t)_{0<t\leq 1}$ is thus a Cauchy family in $\mathcal{C}^{r-\epsilon}$, so it converges in that space as $t$ goes to 0. We denote by $\Lambda$ its limit. Since

$$\langle \Lambda, p_a(x,\cdot)\rangle = \lim_{t\downarrow 0}\langle \mathbb{I}_0^t, p_a(x,\cdot)\rangle = \lim_{t\downarrow 0}\mathbb{I}_0^{t+a}(x) = \mathbb{I}_0^a(x) \tag{11}$$

behaves like $a^{r/2}$, uniformly in $x\in\mathbb{R}^d$, we actually have $\Lambda\in\mathcal{C}^r$. We see from (10) and (11) with $s'=0$ and $s=t$ that we have the $x$-uniform bound $|\langle\Lambda - \Lambda_x, p_t(x,\cdot)\rangle|\lesssim t^{\gamma/2}$ from the statement.

*Uniqueness* – Two possible reconstructions $\Lambda^1$ and $\Lambda^2$ necessarily satisfy the $(x,t)$-uniform estimate

$$|\langle \Lambda^1 - \Lambda^2, p_t(x,\cdot)\rangle| \lesssim t^{\gamma/2}.$$

As for any Schwartz function $\varphi$ the functions $P_t(\varphi)$ converge to $\varphi$ in the smooth topology, one has from the symmetry of the kernels $p_t$ and the positivity of the exponent $\gamma$ that

$$\langle \Lambda^1 - \Lambda^2, \varphi\rangle = \lim_{t\downarrow 0}\int\langle \Lambda^1 - \Lambda^2, p_t(x,\cdot)\rangle\varphi(x)\, dx = \lim_{t\downarrow 0} O(t^{\gamma/2}) = 0. \qquad\square$$

The proof makes it clear that we can trade the condition (2) that defines a $\gamma$-coherent germ for the condition

$$|\langle \Lambda_y - \Lambda_x, p_t(x,\cdot)\rangle| \lesssim \sum t^{\beta_i/2}\Big(|y-x| + t^{\frac{1}{2}}\Big)^{\gamma-\beta_i}$$

for a finite sum where each $\beta_i \leq \gamma$.

*Exercice 1* – Let $(\Lambda_x)_{x\in\mathbb{R}^d}$ be a $\gamma$-coherent germ ($\gamma > 0$) such that $(x,y) \mapsto \Lambda_x(y)$ is continuous. Use the construction of the proof of Theorem 3 to prove that the reconstruction $\Lambda$ of the germ is the continuous function $\Lambda(x) = \Lambda_x(x)$.

**Proposition (Young product)** – *Let $f$ be $r_1$-Hölder continuous function on $\mathbb{R}^d$ for some non-integer $r_1 > 0$, and let $g \in \mathcal{C}^{r_2}$ for some $r_2 < 0$. Denote by $F_x(z) = \sum_{|\ell|<r_1} \frac{(\partial^\ell f)(x)}{\ell!}(z-x)^\ell$ the Taylor polynomial of $f$ of order $r_1$ based at the point $x$. One defines an $(r_1+r_2)$-coherent germ setting*

$$\Lambda_x = F_x g.$$

*If $r_1 + r_2 > 0$, its unique reconstruction defines the product $fg$.*

**Proof** – Note that $y$ and $x$ play a symmetric role in the definition of a $\gamma$-coherent germ. Since $F_y(z) = F_x(z) + \sum_{|m|<r_1} R_x^m(y)\frac{|z-y|^m}{m!}$, with $|R_x^m(y)| \lesssim |y-x|^{r_1-|m|}$, we estimate here

$$|\langle \Lambda_x - \Lambda_y, p_t(y,\cdot)\rangle| \lesssim \sum_{|m|<r_1} |y-x|^{r_1-|m|} |\langle (\cdot-y)^m g, p_t(y,\cdot)\rangle|.$$

We write $\langle (\cdot-y)^m g, p_t(y,\cdot)\rangle = \langle g, (\cdot-y)^m p_t(y,\cdot)\rangle$. An elementary induction tells us that

$$(\cdot-y)^m p_t(y,\cdot) = \sum c_{ab} t^b \partial_y^a p_t(y,\cdot)$$

for some constants $c_{ab}$ and a sum over $2b - |a| = |m|$ and $\max(b,|a|) \leq |m|$. We then write

$$\langle g, (\cdot-y)^m p_t(y,\cdot)\rangle = \sum c_{ab} t^b \langle g, \partial_y^a p_t(y,\cdot)\rangle$$

and we note that we have

$$|\langle g, \partial_y^a p_t(y,\cdot)\rangle| = \left|\int \partial_y^a p_{t/2}(y,z) \langle g, p_{t/2}(z,\cdot)\rangle dz\right| \lesssim_g t^{r_2/2} t^{-|a|/2}$$

because $g \in \mathcal{C}^{r_2}$ and from the estimate (1). This gives the bounds

$$|\langle g, (\cdot-y)^m p_t(y,\cdot)\rangle| \lesssim t^{(r_2+|m|)/2}$$

and

$$|\langle \Lambda_x - \Lambda_y, p_t(y,\cdot)\rangle| \lesssim t^{(r_2+|m|)/2} \sum_{|m|<r_1} |y-x|^{r_1-|m|}.$$

Writing $r_1 - |m| = r_1 + r_2 - (r_2 + |m|)$, we see on this estimate that the germ is $(r_1+r_2)$-coherent. □

**Remarks** – **1.** There is a natural notion of distance in the space of germs that turns the reconstruction map into a continuous map; see e.g. Theorem 12.7 in [4]. Note the obvious fact that the reconstruction map $\mathscr{R}$ is linear. It is also continuous if one uses on the space of germs the previous norm.

*2.* A slightly different setting is needed to deal with anisotropic spaces like the parabolic spaces $\mathbb{R} \times \mathbb{T}^d$ or $\mathbb{R} \times \mathbb{R}^d$. We use in that case a different semigroup $(P_t)_{t>0}$ with the same scaling properties as the parabolic space. See Hoshino's work [6] for a very nice reference.

*Exercice 2 –* Can you prove that this definition of $fg$ coincides with the definition given in Appendix 1 in terms of the Littlewood-Paley decompositions of $f$ and $g$?

## 5 Germs and PDEs

The reconstruction theorem opens the door to trading the task of looking for a function $u$ solution of some equation for the task of looking for a $\gamma$-coherent germ ($\gamma > 0$) uniquely characterized as the solution of some equation set in the space of germs. We test that idea on a well-posed equation.

### 5.1 PDEs as equations on germs

If for instance $u$ satisfies a parabolic equation on $(0, T] \times \mathbb{T}^d$ of the form

$$(\partial_t - D^2) u = F(u) \tag{12}$$

with some nice initial condition and some smooth function $F$, we would like to see if we can characterize $u$ by some equation on a germ $(U_z)_{z \in [0,T] \times \mathbb{T}^d}$ whose reconstruction is $u$. Let a deterministic function $\zeta \in C([0,T] \times \mathbb{T}^d)$ be given; it plays for the moment the role of a placeholder for an irregular function, later a random distribution. In the model situation for us

$$F(v) = f(v) \zeta. \tag{13}$$

(We would treat similarly the case where $F(v) = f(v)\zeta + g(v, \partial v)$.) Since $u$ is meant to be well approximated near each $z$ by $U_z$ one has $U_z(z) = u(z)$, because of Exercice 1, and it seems reasonable to look for $U_z$ as a solution of the equation

$$(\partial_t - D^2) U_z = F(U_z). \tag{14}$$

To make sense of that equation one needs first to define the nonlinear image of a germ by a regular enough function. This is done easily by trading the function $F$ in $F(U_z)$ by some Taylor expansion of $F$ around $U_z(z) = u(z)$. On the other hand it is fairly non-elementary to define an operator $\mathsf{K}$ in the space of ($\gamma$-coherent) germs that is intertwined with the inverse heat operator via the reconstruction operator, so that we have

$$\mathscr{R}\mathsf{K}(V_z)_{z \in (0,T] \times \mathbb{T}^d} = (\partial_t - D^2)^{-1}(\mathscr{R}((V_z)_{z \in (0,T] \times \mathbb{T}^d})) = (\partial_t - D^2)^{-1}(v)$$

for any $\gamma$-coherent germ $(V_z)_{z \in (0,T] \times \mathbb{T}^d}$ with reconstruction $v$.[14]

To make some further progress we look at some special germs of the form

$$U_z(\cdot) = \sum_{\tau \in \boldsymbol{B}} u_\tau(z)\, \Pi_z(\tau)(\cdot) \tag{15}$$

for some finite label set $\boldsymbol{B}$ to be determined. We expect here the coefficient functions $u_\tau$ to depend only on $u$, and the reference functions $\Pi_z(\tau)$ to be continuous and depend only on $\zeta$.

---
14. Yet, this can be done, as for instance in Broux, Caravenna & Zambotti's work [3].

We assume without loss of generality that there is among the symbols $\boldsymbol{B}$ a symbol whose associated reference function is the constant function equal to 1. Denote by $\boldsymbol{1}$ this symbol. We are also free to include in $\boldsymbol{B}$ some abstract monomials $X^k$. Assume as well that $\Pi_z(\tau)(z)=0$ for all the other reference functions, so $U_z(y)$ is close to $U_z(z) = u(z) = u_{\boldsymbol{1}}(z)$ for $y$ close to $z$.

Since the only thing that matters about $U_z$ is its values near $z$, it makes sense to replace $F(U_z)$ in (14) by its Taylor expansion around the constant $u_{\boldsymbol{1}}(z)$, to some fixed degree $\alpha$ to be determined. Plugging the ansatz (15) in this modified version of (14) yields the equation

$$\sum_{\tau \in \boldsymbol{B}} u_\tau(z)(\partial_t - D^2)\Pi_z(\tau) = \sum_{|\ell|\leq \alpha} \frac{\partial^\ell F(u_{\boldsymbol{1}}(z))}{\ell!} \sum_{\tau_1,\ldots,\tau_\ell \in \boldsymbol{B}\setminus\{\boldsymbol{1}\}} u_{\tau_1}(z)\cdots u_{\tau_\ell}(z)\,\Pi_z(\tau_1)\cdots\Pi_z(\tau_\ell).$$

In the model situation (13) this equation reads

$$\sum_{\tau \in \boldsymbol{B}} u_\tau(z)(\partial_t - D^2)\Pi_z(\tau) = \sum_{|\ell|\leq \alpha} \frac{\partial^\ell f(u_{\boldsymbol{1}}(z))}{\ell!} \sum_{\tau_1,\ldots,\tau_\ell \in \boldsymbol{B}\setminus\{\boldsymbol{1}\}} u_{\tau_1}(z)\cdots u_{\tau_\ell}(z)\,\Pi_z(\tau_1)\cdots\Pi_z(\tau_\ell)\zeta. \quad (16)$$

We made some real progress this time. We wish we could identify the terms on both sides of (16).

*1)* Let associate a new symbol $\circ$ to the continuous function $\zeta$. We can define two collections of symbols $\boldsymbol{B}$ and $\boldsymbol{B}^-$, together with a map $\mathcal{I}$ from $\boldsymbol{B}^-$ to $\boldsymbol{B}$, as the smallest collections such that $\{\circ,(X^k)_{k\in\mathbb{N}^d}\} \subset \boldsymbol{B}^-$ and if $\{\tau_1,\ldots,\tau_\ell\} \subset \boldsymbol{B}$ then the unordered tuple $(\tau_1,\ldots,\tau_\ell,\circ)$ belongs to $\boldsymbol{B}^-$ and $\mathcal{I}(\tau_1,\ldots,\tau_\ell,\circ) \in \boldsymbol{B}$, and $(X^k)_{k\in\mathbb{N}^d} \subset \boldsymbol{B}$ and $X^0 = \boldsymbol{1}$. The first element of $\boldsymbol{B}$ constructed here is $\mathcal{I}(\circ)$.

This rule does not necessarily defines a finite set $\boldsymbol{B}$ so we further assume that the symbols $\circ$ and the monomials $X^k$ come with a notion of degree $|\circ|, |X^k| \in \mathbb{R}$ and that the degree $|\sigma|$ of an element $\sigma = \mathcal{I}(\tau_1,\ldots,\tau_\ell,\circ)$ of $\boldsymbol{B}$ is defined inductively as $|\sigma| = |\tau_1| + \cdots + |\tau_\ell| + |\circ| + 2$, and $|X^k| = |k|$. The coefficient 2 here is consistent with Schauder estimate for the heat inverse operator. We *assume* that only finitely many elements of $\boldsymbol{B}$ have a degree smaller than any fixed constant.

*2)* This inductive definition of $\boldsymbol{B}$ allows to identify the corresponding terms in (16). This leads to

*(2a)* an inductive system of equations for the reference functions $\Pi_z(\tau)$ with $\tau \notin \{\boldsymbol{1},(X^k)_{k\geq 0}\} \cap \boldsymbol{B}\}$ and $|\tau| < 2$, where

$$(\partial_t - D^2)\Pi_z(\mathcal{I}(\circ)) = \zeta$$

and for $\tau = \mathcal{I}(\tau_1,\ldots,\tau_\ell,\circ)$

$$(\partial_t - D^2)\Pi_z(\tau) = \Pi_z(\tau_1)\cdots\Pi_z(\tau_\ell)\zeta \quad (17)$$

with each $\Pi_z(\tau)$ null to some order at point $z$,[15]

*(2b)* a system of equations for the coefficient functions $u_\tau$ that starts with $u_{\boldsymbol{1}} = u$ and $u_{\mathcal{I}(\circ)} = f(u_{\boldsymbol{1}})$.

If the system *(2b)* turns out to have a unique solution, the germ (10) will be uniquely defined from the two properties *(2a)* and *(2b)*. We will still have to see that this germ is $\gamma$-coherent for some positive $\gamma$... which is not clear at all presently. This is a complicated way of solving Equation (12) in the model situation (13), and we will certainly not do that when $\zeta$ is continuous![16] However, if successful, this procedure does something very interesting: it decouples in the analysis the task *(2a)* of defining the reference functions $\Pi_z(\tau)$ from the task *(2b)* that identifies the coefficient functions $u_\tau$ of the germ. It is that type of decoupling that is at work in item **6** of the comments following Theorem 1 in Lecture 1.

---

15. By induction the $\Pi_z(\tau)$ depend only on $\zeta$. We always set $\Pi_z(X)(y) \equiv y - z$. We can remove from $\Pi_z(\tau)$ any polynomial in the kernel of the operator $\partial_t - D^2$ without destroying the relation (17).

16. This procedure is fully run in the more complicated setting of some well-posed quasilinear equations in Linares & Otto's expository work [7].

## 5.2 A robust formulation involves an additional drift

The system (17) may no longer make sense if $\zeta$ is not a continuous function but rather a distribution, as is typically the case with stochastic PDEs. Bony's product rule may not be satisfied in that case and the product in the right hand side of (17) may be ill-defined. The above reformulation of the equations (12-13) in the form of the partly decoupled systems *(2a)* and *(2b)* opens an opportunity to deal with this fundamental product problem if one is ready to change the equation in an appropriate way.

When $F(u) = f(u)\zeta + c(u)$ has an additional drift $c(u)$ one needs to add a term $c(u)$ to the right hand side of (16). For now on we think of $\zeta$ as a regularized version $\xi^\epsilon$ of a random distribution $\xi$ and think of the drift $c(\cdot)$ as depending on the law of $\xi$ and $\epsilon$. We will record that $\epsilon$-dependence by some upper indices $\epsilon$ on the objects. Step *1)* above is left unchanged, but we now the new version

$$\sum_\tau u_\tau^\epsilon(z)(\partial_t - D^2)\Pi_z^\epsilon(\tau) = \sum_{\ell,\tau_i} \frac{\partial^\ell f(u_{\mathbf{1}}^\epsilon(z))}{\ell!} \left\{ \prod_{i=1}^\ell u_{\tau_i}^\epsilon(z) \, \Pi_z^\epsilon(\tau_i) \right\} \xi^\epsilon + c^\epsilon(u^\epsilon) \qquad (18)$$

of (16), that is

$$\sum_\tau u_\tau^\epsilon(z)(\partial_t - D^2)\Pi_z^\epsilon(\tau) = \sum_{\ell,\tau_i} \frac{\partial^\ell f(u_{\mathbf{1}}^\epsilon(z))}{\ell!} \left\{ \prod_{i=1}^\ell u_{\tau_i}^\epsilon(z) \, \Pi_z^\epsilon(\tau_i) \right\} \xi^\epsilon + \sum_{k,\tau_j} \frac{\partial^k c^\epsilon(u_{\mathbf{1}}^\epsilon(z))}{k!} \prod_{j=1}^k u_{\tau_j}^\epsilon(z) \, \Pi_z^\epsilon(\tau_j)$$

We analyse the situation inductively.

- The first term constructed in $\boldsymbol{B}$ is $\mu_1 = \mathcal{I}(\circ)$, with associated reference function[17]

$$\Pi_z^\epsilon(\mathcal{I}(\circ)) = (\partial_t - D^2)^{-1}(\zeta) - (\partial_t - D^2)^{-1}(\zeta)(z)$$

and $u_{\mu_1}^\epsilon(z) = f(u_{\mathbf{1}}^\epsilon(z)) = f(u^\epsilon(z))$.

- The second term constructed in $\boldsymbol{B}$ is $\mu_2 = \mathcal{I}(\mathcal{I}(\circ), \circ)$, corresponding to $\ell = 1$ and $\tau_1 = \mathcal{I}(\circ)$ in (17). The product $\Pi_z^\epsilon(\tau_1)\xi^\epsilon$ depends on $\epsilon$.

  – If this product is converging in some appropriate function space as $\epsilon$ goes to 0 we define $\Pi_z(\tau_2)$ as a solution of $(\partial_t - \Delta)\Pi_z^\epsilon(\mu_2) = \Pi_z^\epsilon(\mu_1)\xi^\epsilon$, which further cancels at point $z$ at order $|\mu_2|$.[18]

  – Otherwise we use the term $c^\epsilon(u^\epsilon)$ as a *reservoir* to counter-balance the divergence of $\Pi_z^\epsilon(\mu_1)\xi^\epsilon$ as $\epsilon$ goes to 0. If, for instance, $\Pi_z^\epsilon(\mu_1)\xi^\epsilon - a_2^\epsilon$ is converging in some function space, *in some probabilistic sense*(!), for some constant (or deterministic function) $a_2^\epsilon$ then we choose the function $c^\epsilon(\cdot)$ such that

$$c^\epsilon(u^\epsilon(z)) = -a_2^\epsilon (f'f)(u^\epsilon(z)) + c_3^\epsilon(u^\epsilon(z)),$$

for some other functions $c_3^\epsilon(\cdot)$,[19] left as a free parameter at this stage of the inductive procedure. We see that the converging term $\Pi_z^\epsilon(\mu_1)\zeta^\epsilon - a_2^\epsilon$ appears in the right hand side of the system (18), in the terms

$$f'(u^\epsilon(z))u_{\mu_1}^\epsilon(z)\, \Pi_z^\epsilon(\mu_1)\xi^\epsilon + c^\epsilon(u^\epsilon(z)) = (f'f)(u^\epsilon(z))\{\Pi_z^\epsilon(\mu_1)\xi^\epsilon - a_2^\epsilon\} + c_3^\epsilon(u^\epsilon(z)).$$

We set in that case $u_{\mu_2}^\epsilon(z) = (f'f)(u^\epsilon(z))$ and define $\Pi_z^\epsilon(\mu_2)$ by the identity $(\partial_t - D^2)\Pi_z^\epsilon(\mu_2) = \Pi_z^\epsilon(\mu_1)\xi^\epsilon - a_2^\epsilon$, with an appropriate cancellation property at $z$.

---

17. Recall we need to have $\Pi_z(\tau)(z) = 0$ for all the $\Pi_z(\tau)$ for $U_z$ to qualify as a germ, from Exercice 1.
18. We can remove for that purpose an affine function to any candidate $\Pi_z^\epsilon(\mu_2)$ if needed.
19. The function $c_3^\epsilon(\cdot)$ has a lower index 3 as it will be used in the third step of the iteration. We are presently at the second step of the induction.

• The next iteration step involves $\mu_3 = \mathcal{I}(\mu_1, \mu_1, \circ)$ and reads

$$u^\epsilon_{\mu_3}(z)(\partial_t - D^2)\Pi^\epsilon_z(\mu_3) + (\cdots) = \frac{f^{(2)}(u^\epsilon(z))}{2} (u^\epsilon_{\mu_1}(z))^2 \Pi^\epsilon_z(\mu_1)^2 \xi^\epsilon + (\cdots). \tag{19}$$

We have inside $\Pi^\epsilon_z(\mu_1)^2 \xi^\epsilon$ two diverging products, each of the form $\Pi^\epsilon_z(\mu_1)\xi^\epsilon$. But we know how to possibly cure them from the previous iteration step. So assume this time that the term

$$\Pi^\epsilon_z(\mu_1)^2 \xi^\epsilon - 2a^\epsilon_2 \Pi^\epsilon_z(\mu_1)$$

is converging in some function space. We thus need some term

$$-2a^\epsilon_2 \frac{f^{(2)}(u^\epsilon(z))}{2} u^\epsilon_{\mu_1}(z) = -a^\epsilon_2 f^{(2)}(u^\epsilon(z)) f(u(z)) \Pi^\epsilon_z(\mu_1)$$

from the contributions of $c^\epsilon$. Since $c^\epsilon = -a^\epsilon_2 (f'f) + c^\epsilon_3$, the term $(c^\epsilon)' = -a^\epsilon_2 f^{(2)} f - a^\epsilon_2 (f')^2 + (c^\epsilon_3)'$ provides it for free, with some other terms that we do not use at that stage of the iteration.

• You can imagine how one can iterate this type of procedure in a systematic way. As we proceed along the iteration, we build progressively a partial expansion of $u$, and the initial product that appears in the system (17) needs to be updated at each iteration step. We need to understand which counterterms make the updated product convergent in some function space, *in some probabilistic sense*, and how these counterterms can be obtained from the Taylor expansion of a function $c^\epsilon_k(u)$, given that we already have a partial expansion of $u$ at that stage of the iteration.

We expect that, after exhausting all the elements of $\boldsymbol{B}$ after a finite number $n$ of steps, the system corresponding to (17) has only some right hand sides that are converging *in some probabilistic sense* as $\epsilon > 0$ goes to 0. We note the important fact that the non-polynomial part of the system *(2b)* of equations for the coefficient functions $u_\tau$ is unaffected by the drift $c^\epsilon$.

This analysis is consistent with Theorem 1 in Lecture 1, which says that there is a choice of drift $c^\epsilon$ such that the renormalised equation

$$(\partial_t - D^2)u^\epsilon = f(u^\epsilon)\xi^\epsilon + c^\epsilon(u^\epsilon)$$

has a reformulation in the space of germs that is meaningful in the limit where $\epsilon > 0$ goes to 0. A solution of the formal $\epsilon = 0$ equation will then be defined as the reconstruction of the limit germ. The fundamental objects of regularity structures introduced in the next lectures allow to realize that program.

## References

[1]  H. Bahouri, J.-Y. Chemin and R. Danchin. *Fourier Analysis and Nonlinear Partial Differential Equations*. Grundlehren der mathematischen Wissenschaften, vol. 343, Springer (2012).

[2]  I. Bailleul and M. Hoshino. *A tourist's guide to regularity structures and singular stochastic PDEs*. EMS Surv. Math. Sci. (2025).

[3]  L. Broux, F. Caravenna and L. Zambotti. *Hairer's multilevel Schauder estimates without regularity structures*. Trans. Am. Math. Soc. **377** (2024), 6981-7035. *The title says it all.*

[4]  F. Caravenna and L. Zambotti. *Hairer's reconstruction theorem without regularity structures*. EMS Surv. Math. Sci. 7(2) (2021), 207–251. *Caravenna & Zambotti extracted from Hairer's work the notion of germ and the mechanics of the reconstruction theorem independently of the setting of regularity structures.*

[5]  M. Hairer. *A theory of regularity structures*. Invent. Math. **198**(2) (2014), 269–504.

[6]  M. Hoshino. *A semigroup approach to the reconstruction theorem and the multilevel Schauder estimate*. Ann. Henri Lebesgue **8** (2025), 151–180. *A short self-contained proof of the reconstruction and multilevel Schauder estimates in regularity structures in the very general setting of regularity-integrability structures.*

[7]  P. Linares and F. Otto. *A tree-tree approach to regularity structures: the regular case for quasi-linear equations*. arXiv:2207.10627 (2022). *This is an expository article on the multi-index approach to quasilinear parabolic equations in a setting where the noise is regular and no renormalisation is needed.*

# Lecture 3 – The mechanics of local expansions

We saw in Lecture 2 that germs are involved in a reformulation of stochastic PDEs that may be robust to some vanishing regularization of some irregular "noises" in the equations. The introduction in Lecture 2 of some particular germs of the form $\sum_{\tau \in \boldsymbol{B}} u_\tau(z) \, \Pi_z(\tau)(\cdot)$ was a real input in our analysis. We examine this special structure from a naive point of view in this lecture. This elementary picture turns out to be deep and leads directly to the notion of regularity structure: the algebraic backbone of the computations that we do when we want to change the base point of a germ and look at the corresponding increment of that germ.

## 6 Algebra as the language of local expansions

### 6.1 The Taylor expansion device

The usual notion of local description of a function near an arbitrary point involves some Taylor expansions and amounts to comparing a function to a polynomial centered at $x$

$$v(\cdot) \simeq \sum_n v_n(x) \, (\cdot - x)^n \tag{20}$$

near $x$. The sum over $n$ is finite, the approximation quantified, and we end up describing in this way the class of Hölder functions with real positive regularity exponents. From the binomial identity

$$(\cdot - x)^n = \sum_{\ell \leq n} \binom{n}{\ell} (y - x)^{n-\ell} (\cdot - y)^\ell \tag{21}$$

one gets a local description of $f$ near any other point $y$

$$v(\cdot) \simeq \sum_\ell \left( \sum_{n; \ell \leq n} v_n(x) \binom{n}{\ell} (y-x)^{n-\ell} \right) (\cdot - y)^\ell. \tag{22}$$

Writing (20) with $y$ in place of $x$, the expansion (22) brings the important insight that one has

$$v_\ell(y) \simeq \sum_{n \geq \ell} v_n(x) \binom{n}{\ell} (y - x)^{n-\ell} \tag{23}$$

for $y$ near $x$. We note that the coefficients $v_n(x)$ of $v_\ell(y)$ in (23) are associated with some multi-indices $n \geq \ell$. The relation (23) between the coefficients in the expansion at different points is a *consistency condition*.[20] We talk about the polynomials $((\cdot - x)^n)_{x \in \mathbb{R}^d, n \in \mathbb{N}^d}$ as an *expansion device*.

A more general local expansion device involves a collection of functions or distributions $(\Pi_x \tau)(\cdot)$, with labels $\tau$ in some finite set $\boldsymbol{B}$, indexed by the points of the state space. Mimicking the classical Taylor situation, we will consider below a class of functions/distributions that are locally described near each point $x$ of the state space by an expression of the form

$$v(\cdot) \simeq \sum_{\tau \in \boldsymbol{B}} v_\tau(x) \, \Pi_x(\tau)(\cdot) \tag{24}$$

for some coefficients $v_\tau(x)$ and a finite sum.[21] One has for instance $\boldsymbol{B} = \mathbb{N}^d$ and $\Pi_x(k)(\cdot) = (\cdot - x)^k$ in the polynomial setting of this section. Let me give you before a simple example where the state space is $[0,1]$ and the expansion device is naturally indexed by another set than $\mathbb{N}^d$.

---

20. Conversely, this consistency condition for the coefficients $(v_n)_{n \in \mathbb{N}^d}$ is required to ensure the existence of a function $v$ satisfying (1) – an elementary, but not so trivial, statement.

21. The above expression implicitly assumes that the coefficients $v_\tau(x)$ are some functions of $x$.

## 6.2 An expansion device for controlled ordinary differential equations

Let $h \in C^\infty([0,1], \mathbb{R}^\ell)$ be an $\mathbb{R}^\ell$-valued path and $V_1, \ldots, V_\ell$ be some smooth vector fields on $\mathbb{R}^d$ with bounded derivatives. Let $m$ be the unique solution path of the differential equation

$$\dot{m}_t = \sum_{i=1}^\ell V_i(m_t) \dot{h}_t$$

with a given initial condition. Identifying a vector field with a first order differential operator[22] one has

$$m_t \simeq m_s + \sum_{i=1}^\ell V_i(m_s)(h_t^i - h_s^i) \;+\; \sum_{1 \le j,k \le \ell} (V_j V_k)(m_s) \int_s^t \int_s^{r_1} \dot{h}_{r_2}^j \dot{h}_{r_1}^k dr_2 dr_1$$
$$+ \sum_{1 \le a,b,c \le \ell} (V_a V_b V_c)(m_s) \int_s^t \int_s^{r_1} \int_s^{r_2} \dot{h}_{r_3}^a \dot{h}_{r_2}^b \dot{h}_{r_1}^c dr_3 dr_2 dr_1$$

You can guess what would be the next term if we wanted to expand further. Even though everything is smooth here, so we could wirte $m_t \simeq m_s + m_s'(t-s) + (\cdots)$, it seems meaningful to use in an expansion device some reference objects that are related to the problem at hand, here the above dynamics for the path $m$. The above expansion for the path $m$ around time $s$ suggests that we choose here the set of *ordered* tuples

$$\boldsymbol{B} = \{\mathbf{1}\} \sqcup [\![1,\ell]\!] \sqcup [\![1,\ell]\!]^2 \sqcup [\![1,\ell]\!]^3$$

and set $\Pi_s(\mathbf{1})(t) = 1$ for all $s,t$ and

$$\Pi_s(i)(t) = h_t^i - h_s^i$$

and

$$\Pi_s((j,k))(t) = \int_s^t \int_s^{r_1} \dot{h}_{r_2}^j \dot{h}_{r_1}^k dr_2 dr_1$$

and

$$\Pi_s((a,b,c))(t) = \int_s^t \int_s^{r_1} \int_s^{r_2} \dot{h}_{r_3}^a \dot{h}_{r_2}^b \dot{h}_{r_1}^c dr_3 dr_2 dr_1.$$

We note the particular fact that this unique expansion device can be used for all the ordinary differential equations controlled by $h$, for all smooth vector fields $V_1, \ldots, V_\ell$ in some arbitrary finite dimensional (or Banach) ambiant space. This 'universal' feature of this expansion device confirms that this choice of expansion device ought to be relevant.[23]

## 6.3 The reference objects for general expansion devices

We come back to the setting of a general expansion device as in (24) and try to identify from the model situation of Taylor polynomials what features we should ask from such a device. What follows is simple. Starting from an expansion of an object based at some point $x$, we will ask that if we change the base point in an expansion for another base point $y$, in some possible way, then the new expansion coincides with the expansion of the object based at $y$. We will see that this simple coherence condition on our expansions make appear some algebraic structures on the set of symbols $\{\tau\}$ that encode the expansion procedure in a systematic way.

We proceed linearly along the analysis of the model setting of Section 6.1.

---

[22]. One has for isntance $V = V(\mathrm{Id})$, $V_2 V_1 = (DV_1)(V_2)$ and $V_3 V_2 V_1 = V_3(V_2(V_1)) = D^2 V_1(V_2, V_3) + DV_1(DV_2(V_3))$.

[23]. We would use a situation-dependent expansion device if we used the classical Taylor monomials $(t-s)^n$ to give some local expansion of a solution path.

***What counterpart to the binomial expansion (21) could we use in a more general setting?*** As in (21), it seems meaningful to impose that $\Pi_x(\tau)$ is a(n $(x, y)$-dependent) linear combination of some $\Pi_y(\sigma)$ for any other point $y$. We express the situation by introducing the real vector space $T$ spanned by $\boldsymbol{B}$ the identity and by asking that one has

$$\Pi_x(\tau)(\cdot) = \Pi_y(\Gamma_{yx}(\tau))(\cdot) \tag{25}$$

with a linear map $\Gamma_{yx}: T \to T$ that describes the above mentioned linear combination. Since the roles of $x$ and $y$ are exchangeable, the linear maps $\Gamma_{yx}$ need to be invertible with $\Gamma_{yx}^{-1} = \Gamma_{xy}$.

***A peculiarity of the Taylor monomials.*** In the classical setting of Section 6.1 one uses in the local expansions (20) and (22) the same family $((\cdot - x)^n)_{x \in \mathbb{R}^d, n \in \mathbb{N}^d}$ of reference (polynomial) functions to give a local description of both $v$ and its local coefficient functions $v_\ell$. With a more general local description device (24) there is no reason to use the same reference objects for $v$ and its local coefficient functions $v_\tau$. This is particularly clear if the $\Pi_x(\tau)$ are meant to be used to describe possibly some distributions that are not some functions, so some of the reference objects $\Pi_x(\tau)$ should be some distributions. On the other hand, it makes sense to only use some functions as some reference objects to describe the local expansion properties of the functions $v_\tau$.

We introduce for this reason another set of symbols $\boldsymbol{B}^+$ and associate to each $\mu \in \boldsymbol{B}^+$ an $(x \in \mathbb{R}^d)$-indexed reference function[24] $(y \in \mathbb{R}^d) \mapsto g_{yx}(\mu)$ that plays the role of the monomial functions $y \mapsto (y - x)^{n-\ell}$ in (22), with $\mu$ in the role of $n - \ell$ here. More precisely, we assume that each coefficient function $v_\sigma$ can be compared near an arbitrary $x$ with a finite linear combination

$$v_\sigma(y) \simeq \sum_{\tau \in \boldsymbol{B}, \mu \in \boldsymbol{B}^+} c_\mu^{\tau\sigma} g_{yx}(\mu) v_\tau(x) \tag{26}$$

of some $v_\tau(x)$, with some constants $c_\mu^{\tau\sigma}$ depending on $\sigma, \tau \in \boldsymbol{B}, \mu \in \boldsymbol{B}^+$ – typically some symmetry constants. In the polynomial setting one has $\boldsymbol{B}^+ = \boldsymbol{B} = \mathbb{N}^d$ and $g_{yx}(k) = (y - x)^k$ and $c_k^{n, \ell} = \mathbf{1}_{k=n-\ell} \binom{n}{\ell}$; this is indeed a symmetry coefficient. The next section shows in particular that a simple consistency requirement implies that the maps $\Gamma_{yx}(\cdot)$ and the maps $g_{yx}(\cdot)$ are related.

### 6.4 Consistency relations under re-expansion

We now go along with the above game of 'expansion, re-expansion then identification' described at the begining of Section 6.3 as a way of uncover the algebraic structure of an expansion device. There are two ways of changing the base point in $v(\cdot) \simeq \sum_{\sigma \in \boldsymbol{B}} v_\sigma(y) \Pi_y(\sigma)(\cdot)$, either in $v_\sigma(y)$ or in $\Pi_y(\sigma)$. We explore the two possibilites. With the above notations one has

$$v(\cdot) \simeq \sum_{\sigma \in \boldsymbol{B}} v_\sigma(y) \Pi_y(\sigma)(\cdot) \simeq \sum_{\tau, \sigma \in \boldsymbol{B}, \mu \in \boldsymbol{B}^+} c_\mu^{\tau\sigma} g_{yx}(\mu) v_\tau(x) \Pi_y(\sigma)(\cdot).$$

Comparing this expansion to the original expansion (24) at the point $x$, we have for (25) the explicit representation $\Gamma_{yx}(\tau) = \sum_{\sigma \in \boldsymbol{B}, \mu \in \boldsymbol{B}^+} c_\mu^{\tau\sigma} g_{yx}(\mu) \sigma$. To lighten the notations, it will be convenient to set

$$\tau/\sigma \equiv \sum_{\mu \in \boldsymbol{B}^+} c_\mu^{\tau\sigma} \mu$$

it is an element of the real vector space $T^+$ spanned by $\boldsymbol{B}^+$. In those terms one has

$$\Gamma_{yx}(\tau) = \sum_{\sigma \in \boldsymbol{B}} g_{yx}(\tau/\sigma) \sigma \tag{27}$$

---

[24]. A consistent notation with the notation $\Pi_x(\tau)(\cdot)$ would be $g_x(\mu)(\cdot)$. There are some good reasons to choose the notation $g_{yx}$. See the identity before (48) in Lecture 4.

and

$$v_\sigma(y) \simeq \sum_{\tau \in \boldsymbol{B}} g_{yx}(\tau/\sigma) v_\tau(x) \qquad (28)$$

We can also derive a transitivity relation similar to (25) for $\Gamma_{zy}$ and $\Gamma_{yx}$. Indeed we can develop the coefficient $v_\tau(x)$ in (25) and, re-indexing the labels, one gets

$$v_\eta(z) \simeq \sum_{\sigma \in \boldsymbol{B}} g_{zy}(\sigma/\eta) v_\sigma(y) \simeq \sum_{\sigma, \tau \in \boldsymbol{B}} g_{zy}(\sigma/\eta) g_{yx}(\tau/\sigma) v_\tau(x).$$

Since we also have

$$v_\eta(z) \simeq \sum_{\tau \in \boldsymbol{B}} g_{zx}(\tau/\eta) v_\tau(x) \qquad (29)$$

it seems natural to impose the identity

$$\sum_{\sigma \in \boldsymbol{B}} g_{zy}(\sigma/\eta) g_{yx}(\tau/\sigma) = g_{zx}(\tau/\eta). \qquad (30)$$

The relations (25) and (30) together will be our counterpart of the binomial identity (21).[25]

## 6.5 Algebra encodes the consistency relations

There is a way to describe the skeleton of the relations (25), (27) and (30) with no mention of any explicit reference functions or distributions. We introduce the splitting maps

$$\Delta : T \to T \otimes T^+ \quad \text{and} \quad \Delta^+ : T^+ \to T^+ \otimes T^+$$

by the identities[26]

$$\Delta \tau = \sum_\sigma \sigma \otimes (\tau/\sigma) \quad \text{and} \quad \Delta^+(\tau/\eta) = \sum_\sigma (\sigma/\eta) \otimes (\tau/\sigma).$$

Using this notation, the transition map (27) reads

$$\Gamma_{yx} = (\text{Id} \otimes g_{yx}) \Delta \qquad (31)$$

and (25) and (30) take the form

$$\Pi_x = (\Pi_y \otimes g_{yx}) \Delta \quad \text{and} \quad g_{zx} = (g_{zy} \otimes g_{yx}) \Delta^+$$

An important property of the splitting maps $\Delta$ and $\Delta^+$ is given by the following relations

$$(\Delta \otimes \text{Id}) \otimes \Delta = (\text{Id} \otimes \Delta^+) \Delta$$

and

$$(\text{Id} \otimes \Delta^+) \Delta^+ = (\Delta^+ \otimes \text{Id}) \Delta^+.$$

These identities express the coherence of the re-expansion procedure.

*Exercice 1* – Prove that the first equality above implies that $\Delta^+(\tau/\eta) = \sum_\sigma (\sigma/\eta) \otimes (\tau/\sigma)$.

---

25. We have two conditions rather than just one as identity (2) actually plays two different roles in Section 6.1.

26. Since there is not a unique way to express an element of $T^+$ in the form $\tau/\sigma$ this 'definition' of $\Delta^+(\tau/\eta)$ may seem problematic. However this formula holds true for the concrete regularity structures of Definition 1 defined below, not as a definiting property but rather as a consequence of the definition of $\Delta^+$ given below.

The transitivity relation $\Gamma_{zx} = \Gamma_{zy}\Gamma_{yx}$ is for instance encoded by the first relation

$$\begin{aligned}\Gamma_{zy}\Gamma_{yx} = (\operatorname{Id} \otimes g_{zy} \otimes g_{yx})(\Delta \otimes \operatorname{Id})\Delta &= (\operatorname{Id} \otimes g_{zy} \otimes g_{yx})(\operatorname{Id} \otimes \Delta^+)\Delta \\ &= (\operatorname{Id} \otimes g_{zx})\Delta = \Gamma_{zx}\end{aligned} \quad (32)$$

We want the family of reference functions $(g_{yx}(\mu))_{\mu \in \boldsymbol{B}^+}$ to be sufficiently rich to describe locally an *algebra of functions*. The cheapest way to ensure this property is to assume that the linear span $T^+$ of $\boldsymbol{B}^+$ has an *algebra* structure and that the maps $g_{yx}$ on $T^+$ are multiplicative.[27]

- We noted after (23) that, in the model situation, all the indices $n$ that appear in the expansion of some coefficient function $v_\ell$ are larger than or equal to $\ell$, and the exponent $n - \ell$ of $(x-y)$ has size $|n| - |\ell| \geq 0$. It seems legitimate to ask similarly that all the basis elements $\tau \in \boldsymbol{B}$ have a notion of size $|\tau|$ and one has the relation $|\tau/\sigma| = |\tau| - |\sigma| \geq 0$. This property will re-appear in Definition 1 under the form of some *gradings* on some vector spaces.

In the end, we find it natural to encode the consistency relations of some local expansions under a re-expansion procedure in the form of a graded vector space $T$, a graded algebra $T^+$, and some splitting maps $\Delta$ and $\Delta^+$ that satisfy some algebraic relations modelled on the consistency relations.

## 7 Regularity structures

We record the essential features of this discussion in the following definition.[28]

**Definition 4.** *A **concrete regularity structure** $\mathscr{T} = (T, T^+)$ is a pair of vector spaces*

$$T = \bigoplus_{a \in A} T_a, \qquad T^+ = \bigoplus_{\alpha \in A^+} T^+_\alpha,$$

*where each vector space $T_a$ and $T^+_\alpha$ is finite dimensional, and such that the following holds.*

- *The index set $A$ for $T$ is a locally finite subset of $\mathbb{R}$ bounded below[29], and $A^+ \subset [0, \infty)$ contains $0$.*

- *There is a splitting map[30] $\Delta^+ : T^+ \to T^+ \otimes T^+$ that turns the vector space $T^+$ into a Hopf algebra[31] with unit $\mathbf{1}_+$, counit $\mathbf{1}'_+$.*

- *There is a splitting map $\Delta : T \to T \otimes T^+$ such that*

$$(\Delta \otimes \operatorname{Id})\Delta = (\operatorname{Id} \otimes \Delta^+)\Delta \quad \text{and} \quad (\operatorname{Id} \otimes \mathbf{1}'_+)\Delta = \operatorname{Id}.^{32} \quad (33)$$

*Moreover for any $a \in A$ one has*

$$\Delta T_a \subset \bigoplus_{\alpha \in A^+} T_{a-\alpha} \otimes T^+_\alpha. \quad (34)$$

---

27. We do not impose an algebra structure on $T$ as some of the elements of $\boldsymbol{B}$ are meant to be associated with some true distributions, so the are not meant to be multiplied with one another.
28. This is a special form of Hairer's more general notion of regularity structure.
29. This condition on $A$ was called *subcriticality* in Lecture 1.
30. The splitting map $\Delta^+$ is usually called a *coproduct*.
31. A Hopf algebra is an algebra such that $\boldsymbol{B}^+_0 = \{\mathbf{1}_+\}$, the index set $A_+$ is stable by addition, there is a splitting map $\Delta^+ : T^+ \to T^+ \otimes T^+$ such that $(\Delta^+ \otimes \operatorname{Id})\Delta^+ = (\operatorname{Id} \otimes \Delta^+)\Delta^+$ and $\Delta^+(\boldsymbol{B}^+_\alpha) \subset \bigoplus_{\alpha_1 + \alpha_2 = \alpha} \boldsymbol{B}^+_{\alpha_1} \otimes \boldsymbol{B}^+_{\alpha_2}$, and there is a counit $\mathbf{1}'_+$ such that $(\mathbf{1}'_+ \otimes \operatorname{Id})\Delta = (\operatorname{Id} \otimes \mathbf{1}'_+)\Delta = \operatorname{Id}$. For such a structure there is automatically a unique algebra morphism $S : T^+ \to T^+$ such that $S(\mathbf{1}_+) = \mathbf{1}_+$ and, writing $\Delta^+ \tau = \sum \tau_1 \otimes \tau_2$, one has $\sum \tau_1 S(\tau_2) = 0$ unless $\tau = \mathbf{1}_+$. This map $S$ is called the *antipode*.
32. This means that $\Delta(\tau) = \tau \otimes \mathbf{1}_+ + (\cdots)$ with the $(\cdots)$ some linear combination of elements $\sigma \otimes \mu$ with $\mu \in T^+ \setminus \{\mathbf{1}_+\}$.

The details of that definition will not play a role right now, so we live them under the rug presently. We sometimes use the notation $\Delta(\tau) = \sum_{\sigma \leq \tau} \sigma \otimes (\tau/\sigma)$. Denote by $\boldsymbol{B}_a$ and $\boldsymbol{B}_\alpha^+$ some linear bases of $T_a$ and $T_\alpha^+$ respectively. The preceding decomposition of $\Delta(\tau)$ is unique if one asks that $\sigma$ runs in the basis $\boldsymbol{B}$. An element $\tau$ of $\boldsymbol{B}_a$ is said to have *degree* (or homogeneity) $a$, also denoted $|\tau|$.

**Two examples (a)** *The elementary example of the polynomial regularity structure.* One has $T = T^+$ and $\Delta = \Delta^+$ in the polynomial regularity structure, where $A = \mathbb{N}$ and a linear basis of $T_n$ is given by the $(X^\ell)_{\ell \in \mathbb{N}^d}$ with $|\ell| = n$, for any $n \in A = \mathbb{N}$. The splitting map $\Delta$ is defined by the relation

$$\Delta(X^n) = \sum_{\ell \leq n} \binom{n}{\ell} X^\ell \otimes X^{n-\ell}$$

where $\ell = (\ell_1, \ldots, \ell_d) \leq n = (n_1, \ldots, n_d)$ means that $\ell_i \leq n_i$ for all $1 \leq i \leq d$.

**(b)** *Controlled ordinary differential equations.* The regularity structure associated with the expansion device of controlled differential equations is described in detail in Appendix 3, together with another Hopf algebra also used for the study of such equations. As in the polynomial example, one has $T = T^+$ in both cases. The basis elements of both algebras are some trees whose vertices are decorated by an integer[33] in the interval $[\![1, \ell]\!]$, and whose splitting maps chop these trees into pieces in different ways (the sum over $\sigma \leq \tau$), with $\sigma$ in the role of the remaining bit of the trimmed tree and the trimmed branches in the role of $\tau/\sigma$.

Note that, at the level of generality of Definition 1, we do not require anything on the nature of the elements of $T$ and we do not ask that the elements of $T$ and $T^+$ are of the same nature. In the setting of singular stochastic PDE the elements of some bases $\boldsymbol{B}$ and $\boldsymbol{B}^+$ of $T$ and $T^+$ are *constructed inductively*. They have as a consequence a- representation in terms of (decorated) *trees*.

*Exercice – 2.* Denote by $G^+$ the set of characters of the algebra $T^+$. The convolution of two characters $g_1, g_2$ is defined as $(g_1 \star g_2)(\tau) \equiv (g_1 \otimes g_2)(\Delta^+(\tau))$, with $c_1 \otimes c_2 = c_1 c_2$ when $c_1$ and $c_2$ are numbers. The counit $\mathbf{1}'_+$ is a unit in $G^+$.

*(2.a)* Recall from footnote 30 the definition of the antipode $S : T^+ \to T^+$ – a multiplicative map. Prove that any character $g \in G^+$ has an inverse for the convolution given by $g \circ S$.

*(2.b)* We associate to any character $g \in G^+$ the linear map $\hat{g} = (\mathrm{Id} \otimes g)\Delta$ from $T$ into itself. Prove that one has $\widehat{g_1 \star g_2} = \hat{g}_1 \hat{g}_2$ for any $g_1, g_2 \in G^+$, as a consequence of the identity (33). So the map $(g \in G^+) \mapsto (\hat{g} \in \mathrm{End}(T))$ is a linear representation.

*The algebraic structure of Definition 1 encodes the consistency of some local expansions device under re-expansion.* The actual definition of the reference objects $\Pi_x(\tau)$ and $g_{yx}(\sigma)$ that define an expansion device played no role in this definition. We now introduce some size conditions on the reference objects $\Pi_x(\tau)$ and $g_{yx}(\sigma)$ that are the analogues of the tautological size condition $|(y-x)^n| \leq |y-x|^{|n|}$ for the polynomial reference functions.

## 8 Models and germs

We denote by $G^+$ the set of real-valued linear and multiplicative maps on the algebra $T^+$.

**Definition 5.** *A **model** over the regularity structure $\mathscr{T}$ is a pair $\mathsf{M} = (\Pi, g)$ of families of maps*

$$(\Pi_x : T \to \mathcal{S}'(\mathbb{R}^d))_{x \in \mathbb{R}^d}, \quad (g_{yx} \in G^+)_{x, y \in \mathbb{R}^d}$$

*with $\Pi_x$ linear, that satisfy*

– *some algebraic relations: One has*

$$(g_{zy} \otimes g_{yx})\Delta^+ = g_{zx}$$

---

33. For the example of Section 6.2, we identify for instance $(a, b, c)$ to the linear tree (a trunk) with three nodes decorated respectively by $a, b, c$, with $a$ the decoration of the root of thetree.

*and, setting*

$$\Gamma_{yx} = (\mathrm{Id} \otimes g_{yx})\Delta$$

*as in* (31)*, one has the compatibility relation*

$$\Pi_x(\tau) = \Pi_y(\Gamma_{yx}(\tau))$$

*for all* $x, y, z \in \mathbb{R}^d$ *and all* $\tau \in T$*;*

– *some size constraints: For each constant* $c \in \mathbb{R}$ *one has*

$$\sup_{\mu \in \boldsymbol{B}^+, |\mu| \leq c} \sup_{x,y \in \mathbb{R}^d} \frac{|g_{yx}(\mu)|}{|y-x|^{|\mu|}} < \infty \tag{35}$$

*and*

$$\sup_{\tau \in \boldsymbol{B}, |\tau| \leq c} \sup_{0 < t \leq 1} \sup_{x \in \mathbb{R}^d} t^{-\frac{|\tau|}{2}} |\langle \Pi_x(\tau), p_t(x, \cdot) \rangle| < \infty. \tag{36}$$

The condition $(g_{zy} \otimes g_{yx})\Delta^+ = g_{zx}$ ensures that the computation (32) holds true, so we indeed have the transitivity relation $\Gamma_{zx} = \Gamma_{zy}\Gamma_{yx}$ for all $x, y, z \in \mathbb{R}^d$. The $\Gamma$ maps encode the analogue of Newton's binomial identity (21).

*Exercices* – *3.* We work in the polynomial regularity structure of Example **(a)**, where $T = T^+$ and $\{\tau\} = \mathbb{N}^d$. Check that setting $\Pi_x(n)(y) = g_{yx}(n) = (y-x)^n$ for $n \in \mathbb{N}^d$ defines a model.

*4.* Show that the condition $|\Pi_x(\tau)(y)| \lesssim |y-x|^{|\tau|}$ implies the condition (36).

*Models*[34] *are introduced as the fundamental brick of some local expansion device. The following definition allows to use them to build some germs. Set* $T_{<\gamma} = \bigoplus_{a \in A, a < \gamma} T_a$ *for any* $\gamma \in \mathbb{R}$ *bigger than* $\min A$*, and denote by* $h_\tau$ *the component of* $h \in T$ *for any* $\tau$ *in the linear basis* $\boldsymbol{B}$ *of* $T$.

**Definition 6.** *Fix a regularity structure* $\mathscr{T}$ *and a model* M *on it. Fix a regularity exponent* $\gamma \in \mathbb{R}$. *One defines the space* $\mathcal{D}^\gamma(\mathscr{T}, g)$ *as the space of functions* $\boldsymbol{v} \colon \mathbb{R}^d \to T_{<\gamma}$ *such that*

$$\|\boldsymbol{v}\|_{\mathcal{D}^\gamma} \equiv \max_{\tau \in \boldsymbol{B}, |\tau| < \gamma} \sup_{x,y \in \mathbb{R}^d} \left( \|\boldsymbol{v}(x)\|_\tau + \frac{|\{\boldsymbol{v}(y) - \Gamma_{yx}(\boldsymbol{v}(x))\}_\tau|}{|y-x|^{\gamma - |\tau|}} \right) < \infty.$$

An element of a $\mathcal{D}^\gamma(\mathscr{T}, g)$ space is called a ***modelled distribution***.[35] Do not get confused: These objects are some *functions*! This choice of name comes from the fact that they can be used to construct some distributions (and some functions as well), as a consequence of Proposition 7 below. The archetype of modelled distribution is given by the lift

$$\boldsymbol{v}(x) = \sum_{|n| < r} \frac{v^{(n)}(x)}{n!} X^n \tag{37}$$

in the polynomial regularity structure of a $\gamma$-Hölder function $v$ with $\gamma > 0$ non-integer. One has $A = \mathbb{N}$ in that case and $\|\boldsymbol{v}\|_{\mathcal{D}^\gamma} < \infty$ means here that the $v^{(n)}(x)$ are bounded by a constant and $v^{(n)}(y) - \left(v^{(n)}(x) + \sum_{|\ell| < r - |n|} \frac{1}{\ell!} v^{(n+\ell)}(x) (y-x)^\ell\right) = O(|y-x|^{\gamma - |n|})$ for all $n \in \mathbb{N}^d$ with $|n| < \gamma$.

**Proposition 7.** *Let* $\boldsymbol{v} \in \mathcal{D}^\gamma(\mathscr{T}, g)$ *for some* $\gamma \in \mathbb{R}$*. One defines a* $\gamma$-*coherent germ setting*

$$\Lambda_x = \Pi_x(\boldsymbol{v}(x)) \qquad (x \in \mathbb{R}^d). \tag{38}$$

---

34. See Equation (2.21) in reference [3] from Lecture 1 for a natural notion of distance on the space of models over $\mathscr{T}$.
35. There is a natural notion of distance on the space of modelled distributions, given after item *(2.ii)* in Lecture 5.

**Proof.** Denote by $\boldsymbol{B}_{<\gamma}$ a linear basis of $T_{<\gamma}$ and by $(h)_\tau$ the component on $\tau \in \boldsymbol{B}_{<\gamma}$ of an element $h \in T_{<\gamma}$. One has

$$\begin{aligned}
|\langle \Lambda_y - \Lambda_x, p_t(x, \cdot) \rangle| &= |\langle \Pi_y(\boldsymbol{v}(y)) - \Pi_x(\boldsymbol{v}(x)), p_t(x, \cdot) \rangle| \\
&= |\langle \Pi_x(\Gamma_{xy}(\boldsymbol{v}(y)) - \boldsymbol{v}(x)), p_t(x, \cdot) \rangle| \\
&= \sum_{\tau \in \boldsymbol{B}_{<\gamma}} |(\Gamma_{xy}(\boldsymbol{v}(y)) - \boldsymbol{v}(x))_\tau| \, |\langle \Pi_x(\tau), p_t(x, \cdot) \rangle| \\
&\lesssim \sum_{\tau \in \boldsymbol{B}_{<\gamma}} |y - x|^{\gamma - |\tau|} t^{\frac{|\tau|}{2}}
\end{aligned}$$

which implies indeed that $(\Lambda_x)_{x \in \mathbb{R}^d}$ is a $\gamma$-coherent germ. $\square$

If $\gamma > 0$ the $\gamma$-coherent germ (38) has a unique reconstruction. You see here that if all the $\Pi_x(\tau)$ belong to some $\mathcal{C}^r$ space with $r < 0$ then the reconstruction theorem only ensures that the reconstruction is (a priori) only a distribution.

From Exercice 1 of Lecture 2, the reconstruction of the germ associated with the particular modelled distribution $\boldsymbol{v}$ in (37) is the $\gamma$-Hölder function $v$ itself. The germ $U_z(\cdot)$ that we introduced in Section 5.1 of Lecture 2 is meant to be the germ associated with a modelled distribution $\boldsymbol{u}(z) = \sum_{\tau \in \boldsymbol{B}} u_\tau(z) \tau$. In this example, the model itself solves a recursive system of PDEs.

*Exercice – 5.* Let $F: \mathbb{R} \to \mathbb{R}$ be a smooth function and $\boldsymbol{v} \in \mathcal{D}^\gamma(\mathcal{T}, g)$ with $\gamma > 0$ be such that $v_\tau = 0$ unless $|\tau| \geq 0$. Assume here that the product of any two symbols $\tau_1, \tau_2$ in $T$ with $|\tau_1| > 0$ and $|\tau_2| > 0$ is indeed an element of $T$. Set $F(\boldsymbol{v})(x) \equiv Q_\gamma\left(\sum_{n \geq 0} \frac{F^{(n)}(v_\mathbf{1})}{n!}(\boldsymbol{v}(x) - v_\mathbf{1}(x)\mathbf{1})^n\right)$, where $Q_\gamma(\tau) = \tau \mathbf{1}_{|\tau| < \gamma}$. Prove that $F(\boldsymbol{v}) \in \in \mathcal{D}^\gamma(\mathcal{T}, g)$.

Let us come back to the singular stochastic PDE

$$(\partial_t - D^2) u = f(u) \xi + g(u, \partial u) \tag{39}$$

of Lecture 1, with some given initial condition. The language of regularity structures, models, germs and their reconstruction will allow us to consider, for some regularity structure associated with the equation[36], for any model $M$ on that regularity structure with $\Pi_x(\circ) = \xi$ for all $x$, and for any function $v$ modelled after $M$, to make sense of the ill-defined quantity $f(v)\xi + g(v, \partial v)$ as the reconstruction of some modelled distribution. A solution to the equation (20) will then essentially be for us a modelled distribution $\boldsymbol{v} = v\mathbf{1} + (\cdots)$ such that $(\partial_t - D^2)v$ coincides with the reconstruction of $f(\boldsymbol{v}) \circ + g(\boldsymbol{v}, \partial \boldsymbol{v})$ and $v$ has some appropriate initial condition. This is a real step forward, but this picture has the conceptual drawback that the dynamics of $v$ is not defined from $v$ only since it involves $\boldsymbol{v}$. We refer to that point as a *lack of interpretability* of the solution – we would call $v$ interpretable if its dynamics involved $v$ only.

We also keep in mind that we would like a robust solution theory where the solution is unique under some appropriate conditions and it depends continuously on all the parameters in the equation. The main parameter in the meaningless equation (39) is the realization $\xi(\omega)$ of the random noise, seen as an element of some function space. Another option would be to consider the random variable $\xi(\cdot)$, with values in that function space, as a parameter. The above mentioned model-dependent formulation of the equation takes the model $M$ as a parameter, which thus needs to be random since $\Pi_x(\circ) = \xi$. Asking continuity of the solution with respect to this parameter offers the possibility that we can construct some random models $M^\epsilon$ that converge to some limit model $\overline{M}$ and such that the solutions $v^\epsilon$ associated with $M^\epsilon$ are interpretable, while the solution $\overline{v}$ associated with $\overline{M}$ will be the limit of the $v^\epsilon$, by continuity. In that setting, while $\overline{v}$ may not be interpretable it appears as a limit of interpretable $v^\epsilon$. Building such models $\overline{M}$ is the object of *renormalization*. We will describe in Lecture 4 a general class of models whose associated solutions $v$ are interpretable and out of which we could construct a limit random model $\overline{M}$.

---

36. It contains a noise symbol $\circ$, as in Section 5.1 of Lecture 2.

# Lecture 4 – Renormalised continuous models

Recall the discussion at the end of Lecture 3. Take a continuous or smooth function $\zeta$, a placeholder for a regularized version of a realization of the random noise. We introduce in the present lecture a class of continuous models that have a particular recursive structure. This point is fundamental in proving that they lead to some interpretable solutions of the equation $(\partial_t - D^2)u = f(u)\zeta + g(u, \partial u)$. Proving that fact is an important step in giving a complete proof of Theorem 1 in Lecture 1, but is beyond the content of the present set of lectures. At an elementary level, the construction of Section 12 will show that the notion of model introduced in the previous is section is not restricted to the trivial case of the Taylor monomials over the polynomial regularity structure, and that there is actually a large class of models indexed by some family of linear maps.

Let $\mathsf{M} = (\Pi, g)$ be a model on a regularity structure $\mathscr{T}$.

**1. What are the $\Pi_x$ maps useful for?** To build germs! We saw in Proposition 7 of Lecture 3 that one can associate to any modelled distribution $\boldsymbol{v} \in \mathcal{D}^\gamma(\mathscr{T}, g)$ a $\gamma$-coherent germ

$$\Lambda_x = \Pi_x(\boldsymbol{v}(x)) = \sum_{\tau \in \boldsymbol{B}} u_\tau(x) \Pi_x(\tau).$$

In case $\gamma$ is positive the reconstruction theorem gives us a unique distribution $\Lambda$ which is well approximated by $\Lambda_x$ near each point

$$\Lambda(\cdot) \simeq \sum_{\tau \in \boldsymbol{B}} u_\tau(x) \Pi_x(\tau)(\cdot).$$

By analogy with the usual Taylor expansions $f(\cdot) \simeq \sum_k f^{(k)}(x)(\cdot - x)^k$, it makes sense to think of the $\Pi_x(\tau)(\cdot)$ as some kind of monomials and to think of the coefficient functions $u_\tau$ as some *generalized derivatives*. The definition in Lecture 3 of the space $\mathcal{D}^\gamma(\mathscr{T}, g)$ of modelled distributions of regularity $\gamma$ fits nicely with that interpretation of the functions $u_\tau$. The $(x \neq y)$-uniform bound on $\frac{\|\boldsymbol{v}(y) - \Gamma_{yx}(\boldsymbol{v}(x))\|_a}{|y-x|^{\gamma-a}}$ means indeed that each $u_\tau$ has an expansion to order $\gamma - |\tau|$ given by the $\tau$-component of $\Gamma_{yx}(\boldsymbol{v}(x))$. This is exactly what happens with $\beta$-Hölder functions $f$ on $\mathbb{R}^d$, for $\beta > 0$. Not only does $f$ have an expansion to order $\beta$, but all its derivatives $f^{(k)}$ have some expansion to order $\beta - |k|$.

All the $\Pi_x(\tau)$ in a model are some distributions, say living in a $\mathcal{C}^r$ space for some $r < 0$. So they satisfy for all $x$ the $(z,t)$-uniform bound $|\langle \Pi_x(\tau), p_t(z, \cdot) \rangle| \lesssim t^{r/2}$. The defining property of a model requires that we have an $(x,t)$-uniform bound $|\langle \Pi_x(\tau), p_t(x, \cdot) \rangle| \lesssim t^{|\tau|/2}$, strictly better than the previous a priori bound when $|\tau| > r$.[37] These improved bound dovetails nicely with the idea that the $\Pi_x(\tau)$ behave as some kind of remainders.

**2. Remainders of what?!** In the model situation of the polynomial regularity structure, the monomials

$$\Pi_x(X^n)(y) = (y-x)^n = y^n - \sum_{0 \le \ell \le n-1} \frac{p_n^{(\ell)}(x)}{\ell!}(y-x)^\ell$$

are the Taylor remainders of order $n$ the functions $p_n : z \to z^n$ at point $x$. We will reinforce in this lecture the picture of $\Pi_x(\tau)$ as a remainder by introducing a map $\boldsymbol{\Pi}$ for which $\Pi_x(\tau)$ is the 'generalized Taylor remainder' of $\boldsymbol{\Pi}(\tau)$ at the point $x$.

**3. The $\Pi_x$ maps cannot reasonably be multiplicative.** We saw in Section 5.2 of Lecture 2, in defining $\Pi_x(\circ \mathcal{I}(\circ))$ as $\zeta K(\zeta) - c_\zeta$ for some constant $c_\zeta$, that we should not ask the maps $\Pi_x$ to be multiplicative. In an inductive construction, the multiplicativity assumption brings back the definition of the $\Pi_x$ map to the definition of $\Pi_x(\mathcal{I}(\sigma))$ in terms of some previously constructed terms. The notion of admissible model introduced below will keep that recursive feature without the multiplicativity assumption.

---

[37]. Note that we have the same base point $x$ both in $\Pi_x(\tau)$ and in $p_t(x, \cdot)$ in $\langle \Pi_x(\tau), p_t(x, \cdot) \rangle$. In all practical situations we have $|\tau| > r$ for all $\tau \in T$.

To work in a bounded domain without boundary and avoid in this way some technical problems, we work in this lecture in the $d$-dimensional torus $\mathbb{T}^d$. As a short hand notation we write $K$ for the convolution operator $(1-D^2)^{-1}$, where $D^2$ stands for the Laplace-Beltrami operator in $\mathbb{T}^d$. In applications to parabolic stochastic PDEs we would typically work on $(0,T]\times\mathbb{T}^d$ or $(0,T]\times\mathbb{R}^d$ and consider the parabolic operator $(\partial_t - D^2)^{-1}$ or $(\partial_t + 1 - D^2)^{-1}$ rather than the elliptic oprator $(1-D^2)^{-1}$. As we have mainly been working so far in a Euclidean setting we consider the elliptic operator $(1-D^2)^{-1}$. These analytic points have no impact on the content of the present lecture.

## 9 Tree structure of the symbols

The regularity structures that we use to study some singular stochastic PDEs have a special structure. Their basis symbols are built recursively from a noise symbol $\circ$, the abstract monomials $(X^k)_{k\in\mathbb{N}^d}$, a formal product operation and some linear operators $(\mathcal{I}_n)_{n\in\mathbb{N}^d}$. These symbols can then be represented by some decorated trees whose vertices may have a noise, a monomial decoration or both, the edges, representing the operators $\mathcal{I}_n$, having an $\mathbb{N}^d$-decoration, and the joining of two edges representing the formal product of the associated quantities. The meaning of the $n$-decoration will be clear from the identity following (3) below. We write $\mathcal{I}$ for $\mathcal{I}_0$. The tree of $T$ with only one node, without decoration, is denoted by $\mathbf{1}$.

The *degree* of a symbol, or a decorated tree, is defined recursively as the sum of the degrees of its vertex and edge decorations, where the degree of $\mathcal{I}_n$ is $2 - |n|$. Any element of $T^+$ is a linear combination of some elements of the form $X^\ell \prod \mathcal{I}_{n_i}(\tau_i)$, for a finite product of trees $\mathcal{I}_{n_i}(\tau_i)$ of positive degree.[38] Some trees thus have a double identity, as some elements of both $T$ and $T^+$. The tree of $T^+$ with only one node, without decoration, is denoted by $\mathbf{1}_+$.

The splitting map $\Delta$ is part of the definition of a concrete regularity structure; it is a priori only constrained by the conditions that appear in this definition. As we have introduced some more structure on the elements of $T$ and $T^+$, with the operators $\mathcal{I}_n$, we ask a little more from the splitting maps $\Delta$ and $\Delta^+$. We further require that $\Delta: T \to T \otimes T^+$ satisfies the identity

$$\Delta(\mathcal{I}_n(\tau)) = (\mathcal{I}_n \otimes \mathrm{Id})\Delta(\tau) + \sum_{|\ell| < |\mathcal{I}_n(\tau)|} \frac{X^\ell}{\ell!} \otimes (\mathcal{I}_{n+\ell}(\tau)) \tag{40}$$

for $\ell$ ranging in $\mathbb{N}^d$, for all $\tau \in T$.[39] One extends $\Delta$ by multiplicativity. We will see in the proof of Lemma 10 below the role played by the sum over $\ell$ in this definition and the reason why it has that form precisely.

*Exercice – 1.* Show that $\Gamma_{yx}(X\circ) = X\circ + (y-x)\circ$ and $\Gamma_{yx}(X\circ\mathcal{I}(\circ)) = X\circ\mathcal{I}(\circ) + (y-x)\circ\mathcal{I}(\circ) + g_{yx}(\mathcal{I}(\circ))X\circ + g_{yx}(X\mathcal{I}(\circ))\circ$.

## 10 Admissible models and interpretation maps

Do you remember the system of PDEs that we found in Section 5.1 of Lecture 2, when we implemented naively this idea of reformulating our singular SPDE as an equation set in a space of germs? The right hand sides of this system typically involve some products that do not make sense if $\zeta$ has low regularity. This is not satisfying as we are looking for a formulation of the equation that has the potential to be robust to replacing $\zeta$ by $\xi^\epsilon$ and sending $\epsilon > 0$ to 0. The notion of admissible that we now introduce is an efficient way of encoding part of this PDE system in a form that avoids the previous product problem. The price to pay for that is to require from the maps $\Pi_x$ that they are not multiplicative.

---

[38]. We can assume without loss of generality that the only symbols of $T$ and $T^+$ of null degree are $\mathbf{1}$ and $\mathbf{1}_+$ respectively.

[39]. One needs to justify that one can indeed construct some splitting maps on $T$ and $T^+$ that satisfy (1). This is doable and not complicated.

**1. An equivalent of the PDE system from Lecture 2.** We would like to work here with some *continuous models* for which
$$\Pi_x(\circ) = \zeta$$
and
$$\begin{aligned}\Pi_x(\mathcal{I}_n(\tau))(y) &= (\partial^n K)(\Pi_x(\tau))(y) + \text{polynomial}\\ &= (\partial^n K)(\Pi_x(\tau))(y) - \sum_{|\ell|<|\mathcal{I}_n\tau|} \frac{(\partial^{n+\ell}K)(\Pi_x(\tau))(x)}{\ell!}(y-x)^\ell\end{aligned} \qquad (41)$$

is the Taylor remainder at $x$ of the function $(\partial^n K)(\Pi_x(\tau))$, at order $|\mathcal{I}_n(\tau)| = |\tau| + 2 - |n|$, for all $x, y$ in the state space. A model satisfying this identity is said to be *admissible*.[40] It turns out that this particular property of a model allows to lift a singular stochastic PDE into an equation on a space of modelled distributions (see point **2** in Lecture 5); it plays a key conceptual and technical role in the analysis of a singular stochastic PDEs.

For $n=0$ the function $\Pi_x(\mathcal{I}(\tau))(\cdot)$ differs from the function $K(\Pi_x(\tau))$ by a polynomial of degree strictly smaller than $|\tau| + 2$. For $-2 < |\tau| < 0$ this polynomial is thus in the kernel of the Laplace-Beltrami operator $D^2$, so, for an admissible model, $(1 - D^2)\Pi_x(\mathcal{I}(\tau))$ is equal to $\Pi_x(\tau)$ up to a polynomial function. In a parabolic setting we would typically have that $(\partial_t - D^2)\Pi_x(\mathcal{I}(\tau)) = \Pi_x(\tau)$. If the $\Pi_x$ were multiplicative and $\tau$ a product of other symbols, we would end up with the system of PDEs from Lecture 2. We thus see the admissibility property as an equivalent way of formulating this system without involving any product problem. Yet, by working with some non-multiplicative models we loose the fact that it suffices to define $\Pi_x(\tau)$ for $\tau \in \{1, \circ, X\}$ to define it completely.

Before looking at how one can define some non-multiplicative models we go back to this idea that the $\Pi_x(\tau)$ should somehow be the 'remainders' of some other distribution $\boldsymbol{\Pi}(\tau)$ and see that this picture dovetails nicely with the admissibility notion just introduced.

**2. Did you say remainder of something?!** We will eventually work with some models for which the maps $\Pi_x$ have the particular form
$$\Pi_x(\tau) = (\boldsymbol{\Pi} \otimes h_x)\Delta(\tau), \qquad (42)$$

for some *fixed* map $\boldsymbol{\Pi} : T \to \mathcal{S}'(\mathbb{T}^d)$ such that $\boldsymbol{\Pi}(X^\ell)(y) = y^\ell$ and
$$\boldsymbol{\Pi}(\mathcal{I}_n(\tau)) = (\partial^n K)(\boldsymbol{\Pi}(\tau))$$

for all $\tau \in T, n \in \mathbb{N}^d$, and some $x$-dependent maps $h_x : T^+ \to \mathbb{R}$. We call the map $\boldsymbol{\Pi}$ the *interpretation map* of the model. The form (41) of $\Pi_x(\tau)$ fits the general philosophy that $\Pi_x(\tau)$ should be some kind of 'Taylor remainder' at $x$ of some quantity $\boldsymbol{\Pi}(\tau)$.[41] The splitting map $\Delta$ in (42) tells us which pieces of $\tau$ are used to build that local expansion and the function $h_x$ describes what coefficients in that expansion are added to $\boldsymbol{\Pi}(\tau)$ to give $\Pi_x(\tau)$.

**Lemma 8.** *The identities (41) and (42) together imply (40).*

**Proof.** One has indeed
$$\begin{aligned}\Pi_x(\mathcal{I}_n(\tau)) &= (\boldsymbol{\Pi} \otimes h_x)\Delta(\mathcal{I}_n(\tau))\\ &= \sum_{\sigma \leq \tau}(\boldsymbol{\Pi} \otimes h_x)(\mathcal{I}_n(\sigma) \otimes (\tau/\sigma)) + \sum_{|\ell|<|\mathcal{I}_n\tau|}\frac{(\cdot)^\ell}{\ell!}h_x(\mathcal{I}_{n+\ell}(\tau))\\ &= \sum_{\sigma \leq \tau}(\partial^n K)(\boldsymbol{\Pi}(\sigma))\,h_x(\tau/\sigma) + \sum_{|\ell|<|\mathcal{I}_n\tau|}\frac{(\cdot)^\ell}{\ell!}h_x(\mathcal{I}_{n+\ell}(\tau))\\ &= (\partial^n K)(\Pi_x(\tau)) + P_{x,\tau}(\cdot)\end{aligned}$$

---

40. We see from (40) that the $n$ in $\mathcal{I}_n$ represents the differentiation operator $\partial^n$ for admissible models.
41. Yet the form (41) does not define $\boldsymbol{\Pi}(\tau)$ right away.

for some polynomial $P_{x,\tau}(\cdot)$.[42] Since we work with a smooth model, Exercice 2 from Lecture 3 tells us that the smooth function $\Pi_x(\mathcal{I}_n(\tau))$ is null at order $|\mathcal{I}_n(\tau)|$ at point $x$. We thus need to have

$$P_{x,\tau}(\cdot) = -\sum_{|k|<|\mathcal{I}_n(\tau)|} \frac{(\partial^{n+k}K)(\Pi_x(\tau))(x)}{k!}(\cdot-x)^k$$

which is indeed formula (40). □

It is true that at that stage the particular polynomial bit in the condition (1) on $\Delta(\mathcal{I}_n(\tau))$ did not play any particular role and could have potentially been replaced by some other polynomial bit provided we adjust the definition of $h_x$. Its particular form will play a crucial role in the proof of Lemma 3 below where we construct some maps that satisfy the relation (42).

We emphasized at the end of Lecture 3 the need to work with some models that give an interpretable notion of solution of the equation

$$(\partial_t - D^2)u = \mathscr{R}(\{\Pi_x(f(\boldsymbol{u}) \circ + g(\boldsymbol{u},\partial\boldsymbol{u}))\}_x).$$

This property is not shared by all the models, even the admissible models. It turns out to be possible to construct a whole family of admissible models that give some interpretable solutions of the preceding equation. These special admissible models have two important features: They are parametrized by some particular linear maps $R:T\to T$, and they are built recursively. Both features play a key role in proving the interpretable character of the corresponding solutions. Proving this fact would require another set of lectures. We introduce here the maps $R$ and the recursive construction of their associated models.

## 11 Preparation maps

Denote by $n_\circ(\tau)$ the number of $\circ$ symbols that appear in the definition of a given decorated tree $\tau$. The following definition was first introduced by Bruned in [1].

**Definition 9.** *A **preparation map** is a linear map $R:T\to T$ such that $R(X^k\tau)=X^k\tau$ for all $k\in\mathbb{N}^d, \tau\in T$, and $R\mathcal{I}_n = \mathcal{I}_n$ for all $n\in\mathbb{N}^d$, and*

$$(R\otimes\mathrm{Id})\Delta = \Delta R, \tag{43}$$

*and for each $\tau\in T$ there exists finitely many $\tau_i\in T$ and constants $\lambda_i\in\mathbb{R}$ such that*

$$R(\tau) = \tau + \sum \lambda_i \tau_i, \quad \text{with} \quad |\tau_i|\geq|\tau| \quad \text{and} \quad n_\circ(\tau_i) < n_\circ(\tau).$$

Recall from the definition of a model that $\Gamma_{yx} = (\mathrm{Id}\otimes g_{yx})\Delta$. Identity (43) encodes the fact that $\Gamma_g = (\mathrm{Id}\otimes g)\Delta$ and $R$ commute for all maps $g$ on $T^+$, as we have

$$\Gamma_g(R\tau) = (\mathrm{Id}\otimes g)\Delta(R(\tau)) = (\mathrm{Id}\otimes g)(R\otimes\mathrm{Id})\Delta\tau = (R\otimes g)\Delta\tau = R(\Gamma_g(\tau)). \tag{44}$$

A preparation map is in particular a perturbation of the identity by some elements with degrees not smaller than the degree of $\tau$ ($|\tau_i|\geq|\tau|$) and that are defined with strictly less noises ($n_\circ(\tau_i)<n_\circ(\tau)$).

*Exercices – 2.* Prove that a preparation map is invertible, so they form a group.

---

[42]. We used (3) in the second equality and (2) in the third equality.

The preparation maps that are used for singular SPDEs behave as follows. The definition of $R(\tau)$ only requires to look at the structure of the tree $\tau$ at its root. If there is only one incoming edge then the three is of type $\mathcal{I}_n(\sigma)$ and $R$ does nothing. If there are multiple incoming edges at the root then $R$ checks whether any of the subtrees of $\tau$ that reach the root have negative size. If there are no such trees then $R$ does nothing. Otherwise one defines for each subtree $\sigma_i$ of $\tau$ of negative size that reaches the root the tree $\tau_i$ as $\tau$ pruned from $\sigma_i$, and the constant $\lambda_i$ as a constant associated to $\sigma_i$ in some way, and we define $R(\tau)$ from the formula above.

## 12 Admissible model associated with a preparation map

One can associate to any preparation map $R$ a continuous admissible model, built recursively. This construction makes clear the idea that a preparation map 'cures' only the divergencies of a tree-dependent analytic quantity that are due to what happens at the root of that tree, while the full renormalization proceeds by induction. The ad hoc manipulations that we did in Section 5.2 of Lecture 2 were implementing the first steps of the type of construction that we do here.

The continuous admissible model associated with a preparation map $R$ is defined from a side family[43] $(\Pi_x^\times(\tau)(\cdot))_{x\in\mathbb{T}^d, \tau\in T}$ of continuous functions constructed by induction on $|\tau|$, that satisfies

$$\Pi_x^\times(\mathbf{1})(y)=1, \quad \Pi_x^\times(\circ)(y)=\zeta(y), \quad \Pi_x^\times(X_i)(y)=y_i-x_i,$$

that is *multiplicative*:

$$\Pi_x^\times(\sigma\tau)(\cdot) = \Pi_x^\times(\sigma)(\cdot)\,\Pi_x^\times(\tau)(\cdot) \qquad (\forall \sigma, \tau \text{ decorated trees})$$

and that satisfies the condition

$$\Pi_x^\times(\mathcal{I}_n(\tau))(y) = (\partial^n K)(\Pi_x^\times(R(\tau)))(y) - \sum_{|\ell|\leq|\mathcal{I}_n(\tau)|} (\partial^{\ell+n} K)(\Pi_x^\times(R(\tau)))(x)\frac{(y-x)^\ell}{\ell!} \tag{45}$$

for all $n\in\mathbb{N}^d$ and $x,y\in\mathbb{T}^d$. We set

$$\Pi_x(\tau)(\cdot) \equiv \Pi_x^\times(R(\tau))(\cdot). \tag{46}$$

The inductive character of this definition of $\Pi_x$ is in (45) and (46). It takes roughly the following form. First, $R$ produces some counterterms, if needed, by turning $\tau$ into $R(\tau)$. Then, for each tree in this possible sum we do the same thing. We go up along each edge that is incoming into the root. Each vertex at the other side of these edges can be seen as the root of the subtree of the ambiant tree that is above this vertex. Then we repeat the process: apply $R$, *etc*.

For each $x\in\mathbb{T}^d$ we define a character $g_x \in G^+$ from its convolution inverse $g_x^{-1} = g_x \circ S$, a map from $T^+$ to $\mathbb{R}\backslash\{0\}$ given by the formula

$$g_x^{-1}(\mathcal{I}_n(\tau)) \equiv -\sum_{|\ell|<|\mathcal{I}_n\tau|} \frac{(-x)^\ell}{\ell!}(\partial^{n+\ell}K)(\Pi_x(\tau))(x) \tag{47}$$

and $g_x^{-1}(X^\ell) = (-x)^\ell$ for all $\ell\in\mathbb{N}^d$, so $g_x(X^\ell) = x^\ell$ since $S(X)=-X$. We complete the definition of our future continuous model by setting

$$g_{yx} \equiv g_y \star g_x^{-1} = (g_y \otimes g_x^{-1})\Delta^+$$

and

$$\Gamma_{yx} \equiv \widehat{g_{yx}} = (\mathrm{Id}\otimes g_{yx})\Delta. \tag{48}$$

---

[43]. We should more properly denote it by $(\Pi_x^{R\times}(\tau)(\cdot))_{x\in\mathbb{T}^d,\tau\in T}$, with an exponent $R$, to emphasize its dependence on $R$. Similarly we should write below $\Pi_x^R(\tau) = \Pi_x^{R\times}(R(\tau))$. We prefer here to use some lighter notations.

The next lemma says that one can construct a map $\mathbf{\Pi}$ such that the structural relation $\Pi_x(\tau) = (\mathbf{\Pi} \otimes g_x^{-1})\Delta(\tau)$ holds. This particular form will make the proof of the transitivity relation $\Pi_x \Gamma_{xy} = \Pi_y$ elementary as we will then have

$$\begin{aligned}\Pi_x \Gamma_{xy} = (\mathbf{\Pi} \otimes g_x^{-1} \otimes g_{xy})(\Delta \otimes \mathrm{Id})\Delta &= (\mathbf{\Pi} \otimes g_x^{-1} \otimes g_{xy})(\mathrm{Id} \otimes \Delta^+)\Delta \\ &= (\mathbf{\Pi} \otimes g_x^{-1} \star g_{xy})\Delta \\ &= (\mathbf{\Pi} \otimes g_y^{-1})\Delta = \Pi_y.\end{aligned}$$

The definition of $\mathbf{\Pi}$ mimicks the definition of the $\Pi_x$. We define a *multiplicative* map $\mathbf{\Pi}^\times$ on the set of all decorated trees such that $\mathbf{\Pi}^\times(\circ)(y) = \xi^\epsilon(y)$ and $\mathbf{\Pi}^\times(X^k)(y) = y^k$ and $\mathbf{\Pi}^\times(\mathcal{I}_n(\tau))(y) = (\partial^n K)(\mathbf{\Pi}^\times(R(\tau)))(y)$, for all $n, \tau, y$. Set[44]

$$\mathbf{\Pi}(\tau) \equiv \mathbf{\Pi}^\times(R(\tau)).$$

**Lemma 10.** *One has $\Pi_x(\tau) = (\mathbf{\Pi} \otimes g_x^{-1})\Delta\tau$ for any decorated tree $\tau$.*

**Proof.** The statement holds for $\tau \in \{1, \circ, X_i\}$. We now proceed by induction on the total number $\sharp(\tau)$ of edge and vertex decorations of $\tau$ to prove that one has both $\Pi_x(\tau) = (\mathbf{\Pi} \otimes g_x^{-1})\Delta(\tau)$ and $\Pi_x^\times(\tau) = (\mathbf{\Pi}^\times \otimes g_x^{-1})\Delta(\tau)$ for all $\tau \in T$.

For $\tau = \mathcal{I}_n(\sigma)$, the fact that $R\mathcal{I}_n = \mathcal{I}_n$ and the induction assumption applied to $\sigma$ justify that

$$\begin{aligned}\Pi_x(\mathcal{I}_n(\tau)) &= \Pi_x^\times(\mathcal{I}_n(\sigma)) \\ &= (\partial^n K)(\Pi_x(\sigma)) - \sum_{|\ell| < |\mathcal{I}_n \tau|} (\partial^{\ell+n} K)(\Pi_x(\sigma))(x) \frac{(\cdot - x)^\ell}{\ell!} \\ &= (\partial^n K)((\mathbf{\Pi} \otimes g_x^{-1})\Delta\sigma) - \sum_{|\ell| < |\mathcal{I}_n \tau|} (\partial^{\ell+n} K)((\mathbf{\Pi} \otimes g_x^{-1})\Delta\sigma)(x) \frac{(\cdot - x)^\ell}{\ell!} \\ &= \mathbf{\Pi}((\mathcal{I}_n \otimes g_x^{-1})\Delta\sigma) - \sum_{|\ell| < |\mathcal{I}_n \tau|} \mathbf{\Pi}((\mathcal{I}_{\ell+n} \otimes g_x^{-1})\Delta\sigma)(x) \frac{(\cdot - x)^\ell}{\ell!} \\ &= \mathbf{\Pi}((\mathrm{Id} \otimes g_x^{-1})\Delta(\mathcal{I}_n \sigma))\end{aligned}$$

with the last equality coming from the structural assumption (42) relating $\Delta$ and $\mathcal{I}_n$.

It remains to check the statement for $\tau = \prod \tau_i$ a finite product of decorated trees, assuming that all decorated trees $\sigma$ with $\sharp(\sigma) < \sharp(\tau)$ satisfy the statement – which is in particular the case of each $\tau_i$. We use the fact that $R$ commutes with the $\Gamma$ maps, identity (44), to write in the second equality

$$\begin{aligned}(\mathbf{\Pi} \otimes g_x^{-1})\Delta(\tau) &= ((\mathbf{\Pi}^\times R) \otimes g_x^{-1})\Delta(\tau) = (\mathbf{\Pi}^\times \otimes g_x^{-1})\Delta(R(\tau)) \\ &= (\mathbf{\Pi}^\times \otimes g_x^{-1})\Delta(R(\tau) - \tau) + (\mathbf{\Pi}^\times \otimes g_x^{-1})\Delta(\tau) \\ &= \Pi_x^\times(R(\tau) - \tau) + \Pi_x^\times(\tau) \\ &= \Pi_x^\times(R(\tau)) \\ &= \Pi_x(\tau)\end{aligned}$$

applying the induction assumption on $R(\tau) - \tau$ and the $\tau_i$ in the fourth equality after noting that

$$(\mathbf{\Pi}^\times \otimes g_x^{-1})\Delta(\tau) = (\mathbf{\Pi}^\times \otimes g_x^{-1})\Delta\Big(\prod \tau_i\Big) = \prod (\mathbf{\Pi}^\times \otimes g_x^{-1})\Delta(\tau_i) = \prod \Pi_x^\times(\tau_i) = \Pi_x^\times(\tau),$$

because the maps $\Delta, \mathbf{\Pi}^\times$ and $g_x^{-1}$ are multiplicative. $\square$

---

[44]. While the map $\mathbf{\Pi}^\times$ is multiplicative, the map $\mathbf{\Pi}$ is not multiplicative unless the preparation map $R$ is multiplicative. In most situations of interest $R$ is not multiplicative.

Since one has $\Gamma_{zy}\Gamma_{yx} = \Gamma_{zx}$ for all $x, y, z \in \mathbb{T}^d$, one has from Lemma 10 the transitivity relation $\Pi_x \Gamma_{xy} = \Pi_y$ for all $x, y \in \mathbb{T}^d$.

It remains to prove some size estimates on $\langle \Pi_x \tau, p_t(x, \cdot) \rangle$ and $g_{yx}$ to justify that (46) and (48) define jointly a continuous model. The following technical fact will be useful in the proof of that fact.

**Lemma 11.** *One has for all $x \in \mathbb{T}^d$ and all decorated trees $\tau \in T$*

$$g_{yx}(\mathcal{I}_n \tau) = \sum_{1 < \theta \leq \tau, |n|-2 < |\theta|} g_{yx}(\tau/\theta)(\partial^n K)(\Pi_y(\theta))(y) - \sum_{|p| < |\mathcal{I}_n \tau|} \frac{(y-x)^p}{p!} (\partial^{n+p} K)(\Pi_x \tau)(x). \quad (49)$$

For the memory, note that (49) has the form

$$g_{yx}(\mathcal{I}_n \tau) = (\partial^n K)(\Pi_y(\theta))(y) - \sum_{|p| < |\mathcal{I}_n(\tau)|} (\partial^{n+p} K)(\Pi_x \tau)(x) \frac{(y-x)^p}{p!} + \sum_{1 < \theta < \tau, |n|-2 < |\theta|} (\ldots)$$

**Proof.** We start with the formula

$$g_y(\mathcal{I}_n \sigma) = \sum_{\theta \leq \tau; 0 < |\mathcal{I}_n \theta|} g_y(\sigma/\theta)(\partial^n K)(\Pi_y(\theta))(y) \quad (50)$$

of recursive type proved below. We then have formula (49) from (50) and the formula (47) for $g_x^{-1}(\mathcal{I}_n(\tau))$ by writing

$$\begin{aligned}
g_{yx}(\mathcal{I}_n(\tau)) &= \sum_{\sigma \leq \tau} g_y(\mathcal{I}_n(\sigma)) g_x^{-1}(\tau/\sigma) + \sum_\ell \frac{y^\ell}{\ell!} g_x^{-1}(\mathcal{I}_{n+\ell}(\tau)) \\
&= \sum_{\theta \leq \sigma \leq \tau, 0 < |\mathcal{I}_n(\theta)|} g_y(\sigma/\theta)(\partial^n K)(\Pi_y(\theta))(y) \, g_x^{-1}(\tau/\sigma) \\
&\quad - \sum_{|\ell+m| < |\mathcal{I}_n(\tau)|} \frac{y^\ell}{\ell!} \frac{(-x)^m}{m!} (\partial^{n+\ell+m} K)(\Pi_x(\tau))(x) \\
&= \sum_{\theta \leq \tau, 0 < |\mathcal{I}_n(\theta)|} g_{yx}(\tau/\theta)(\partial^n K)(\Pi_y(\theta))(y) - \sum_{|p| < |\mathcal{I}_n(\tau)|} \frac{(y-x)^p}{p!} (\partial^{n+p} K)(\Pi_x(\tau))(x).
\end{aligned}$$

– It remains to prove (50). Write momentarily $g(y, \mathcal{I}_n(\sigma))$ for the right hand side of formula (50). We check that $g(y, \mathcal{I}_n(\sigma))$ satisfies the convolution identity

$$(g(y, \cdot) \star g_y^{-1}) \Delta^+(\mathcal{I}_n(\sigma)) = 0$$

that characterizes uniquely the character $g_y$. One has indeed

$$\begin{aligned}
(g(y, \cdot) \star g_y^{-1}) \Delta^+(\mathcal{I}_n(\sigma)) &= \sum_{\theta \leq \mu \leq \sigma; |\mathcal{I}_n(\theta)| > 0} g_y(\mu/\theta)(\partial^n K)(\Pi_y(\theta))(y) \, g_y^{-1}(\sigma/\mu) \\
&\quad - \sum_{|\ell+q| < |\mathcal{I}_n(\sigma)|} \frac{x^\ell}{\ell!} \frac{(-x)^q}{q!} (\partial^{n+\ell+q} K)(\Pi_y(\sigma))(y) \\
&= \sum_{\theta \leq \sigma; |\mathcal{I}_n(\theta)| > 0} (\partial^n K)(\Pi_y(\theta))(y) \mathbf{1}_{\theta=\sigma} - (\partial^n K)(\Pi_y(\sigma))(y) = 0
\end{aligned}$$

We used the identity (47) giving $g_y^{-1}(\mathcal{I}_{n+\ell}(\sigma))$ in the first equality and the relation

$$\Delta^+(\tau/\eta) = \sum_{\eta \leq \sigma \leq \tau} (\sigma/\eta) \otimes (\tau/\sigma)$$

and the fact that $g \star g^{-1} = \mathbf{1}'_+ \mathbf{1}_+$, for any character $g$ on $T^+$, in the second equality. □

**Theorem 12.** *The maps* $(\Pi_x : T \to \mathcal{S}'(\mathbb{T}^d))_{x \in \mathbb{T}^d}$ *and* $(g_{yx} \in G^+)_{x,y \in \mathbb{T}^d}$ *in* (46) *and* (48) *define jointly a continuous admissible model.*

The model associated with a preparation map is called a ***renormalized model***.

**Proof.** The admissible feature of the the model to be comes from its definition

$$\mathbf{\Pi}(\mathcal{I}_n(\tau)) \equiv \mathbf{\Pi}^\times(R(\mathcal{I}_n(\tau))) = \mathbf{\Pi}^\times(\mathcal{I}_n(\tau)) = (\partial^n K)(\mathbf{\Pi}^\times(R(\tau))) = (\partial^n K)(\mathbf{\Pi}(\tau)).$$

We prove each analytic bound separately.

**1.** We prove by induction on $n_\circ(\tau)$ that we have

$$|(\Pi_x^\times(\tau))(y)| + |(\Pi_x(\tau))(y)| \lesssim |y - x|^{|\tau|}.$$

Since we subtract in the definition (45) of $\Pi_x^\times(\mathcal{I}_n(\tau))$ the correct Taylor expansion we have indeed

$$|\Pi_x^\times(\mathcal{I}_n(\tau))| \lesssim |y - x|^{|\mathcal{I}_n \tau|}, \tag{51}$$

so the same estimate holds for $\Pi_x(\mathcal{I}_n(\tau)) = \Pi_x^\times(R(\mathcal{I}_n(\tau))) = \Pi_x^\times(\mathcal{I}_n(\tau))$.

Take now $\tau = X^a \circ \prod_{i=1}^m \mathcal{I}_{a_i}(\sigma_i)$ and write $R(\tau) = \tau + \sum \lambda_j \tau_j$, with $|\tau_j| \geq |\tau|$ and $n_\circ(\tau_j) < n_\circ(\tau)$ for each $j$. The induction assumption applies to each $\tau_j$, and we have

$$\begin{aligned}
|\Pi_x(\tau)(y)| = |\Pi_x^\times(R(\tau))(y)| &\lesssim |\Pi_x^\times(\tau)(y)| + \sum_j |\Pi_x^\times(\tau_j)(y)| \\
&\lesssim |y - x|^{|a| + \sum_i |\mathcal{I}_{a_i}(\sigma_i)|} + \sum_j |y - x|^{|\tau_j|} \lesssim |y - x|^{|\tau|}.
\end{aligned}$$

We used in the second inequality the multiplicativity of $\Pi_x^\times$, the induction assumption, and the fact that $\xi^\epsilon = \Pi_x^\times(\circ)$ is bounded on the bounded state space $\mathbb{T}^d$.

**2.** We use Lemma 11 to prove by induction on $|\mu|$ that $|g_{yx}(\mu)| \lesssim |y - x|^{|\mu|}$ for all $\mu \in T^+$. It suffices to prove that bound for $\mu = \mathcal{I}_n(\tau) \in T^+$. In the formula of Lemma 11 the term $g_{yx}(\mathcal{I}_n(\tau))$ appears as the sum of two terms. One is the Taylor remainder to order $|\mathcal{I}_n(\tau)|$ of the (smooth) function $y \mapsto (\partial^n K)(\Pi_y(\sigma))(y)$ at point $x$. This term has the expected size. The other term is the sum of the terms $g_{yx}(\tau/\sigma)(\partial^n K)(\Pi_y(\sigma))(y)$ where $1 < \sigma < \tau, |n| - 2 < |\sigma|$. The term $(\partial^n K)(\Pi_y(\sigma))(y)$ is bounded and $|g_{yx}(\tau/\sigma)| \lesssim |y - x|^{|\tau/\sigma|}$, by induction as $|\tau/\sigma| < |\tau|$. The conclusion follows as

$$|\tau/\sigma| = |\tau| - |\sigma| < |\tau| - |n| + 2 = |\mathcal{I}_n(\tau)|$$

and we work in the bounded space $\mathbb{T}^d$. □

Everything here depends on the preparation map $R$ and the continuous function $\zeta$, since $\Pi_x^\times(\circ)(y) = \zeta(y)$. The continuous admissible model associated with the trivial preparation map $R = \text{Id}$ is called the *canonical model*. Write $\Pi_x^{\text{can}, \zeta}$ for it, and we write $\Pi_x^{R, \zeta}$ for the maps associated with a preparation map $R$ and a continuous function $\zeta$. It is elementary to see that $\Pi_x^{R, \zeta}(\tau)$ is, for any $\tau \in \mathbf{B}$, a linear combination of some elements of the form $\Pi_x^{\text{can}, \zeta}(\sigma)$. We note that $\Pi_x^{\text{can}, \zeta}$ is multiplicative while the map $\Pi_x^{R, \zeta}(\cdot)$ is generally not multiplicative.

### References

[1] Y. Bruned. *Recursive formulae in regularity structures*. Stoch. PDE: Anal. Comp., **6** (2018):525–564.

This work gives a recursive construction of the regularity structure introduced by Bruned, Hairer & Zambotti. Preparation maps are introduced in this work.

# Lecture 5 – Back to the future

We have mainly been working so far in a Euclidean setting to keep concentrated on the essential features of what we saw. Everything has a parabolic counterpart. As far as space is concerned, we work in the bounded, boundaryless, setting of $\mathbb{T}^d$ rather than $\mathbb{R}^d$ or a domain with some boundary to avoid some technical problems.

Let $\xi$ be a random spacetime distribution with almost sure regularity smaller than $-1$ in the Besov-Hölder scale $(B^\alpha_{\infty\infty})_{\alpha\in\mathbb{R}}$. Denote by $(\Omega, \mathcal{F}, \mathbb{P})$ the probability space where it is defined. Denote by $\mathcal{F}(u_0)(t) = e^{tD^2}(u_0)$ the free propagation operator of the Laplace-Beltrami operator. We consider the singular stochastic PDE

$$(\partial_t - D^2)u = f(u)\xi + g(u, \partial u), \tag{52}$$

with initial condition $u_0$ in some function space, in its integral form

$$u = \mathcal{F}(u_0) + (\partial_t - D^2)^{-1}(f(u)\xi + g(u, \partial u)) \tag{53}$$

This reformulation will make the parallel with equation (55) below clearer. In both formulations, (52) and (53), these equations are presently only formal due to the fundamental product problem involved in the term $f(u)\xi$. Denote by $\xi^\epsilon$ a regularization of $\xi$ by spacetime convolution. Recall the main result of Lecture 1.

**Theorem 13.** *Under some mild conditions on the noise $\xi$ and the initial condition $u_0$, for any mollification procedure $\xi \mapsto \xi^\epsilon$ there exists a deterministic explicit function $c^\epsilon$ and a random positive time $T = T(\omega)$ such that the solutions to the classical and locally well-posed equation*

$$(\partial_t - D^2)u^\epsilon = f(u^\epsilon)\xi^\epsilon + g(u^\epsilon, \partial u^\epsilon) + c^\epsilon(u^\epsilon, \partial u^\epsilon) \tag{54}$$

*with initial condition $u_0$ converge in probability in $C([0,T] \times \mathbb{T}^d, \mathbb{R})$ to some limit function $u$ that doe not depend on the mollification procedure.*

We assume below that $u_0$ is smooth to avoid the introduction of an ad hoc space of modelled distribution that can handle some initial conditions of low regularity.

## 13  Structure of a proof of Theorem 13

We have introduced in the previous lectures all the objects that we need to describe the structure of the proof of Theorem 13. They are not yet perfectly tailored to the study of singular stochastic PDEs[45] but this is enough to have a clear idea of what is involved below.

***1. A regularity structure for each nice enough singular SPDE.*** *One can associate to each equation* (Eq) *in a large class of singular stochastic PDEs a (rich enough) regularity structure $\mathscr{T}_{(Eq)}$ and a positive exponent $\gamma > 0$.*

Following Lecture 2 we trade the task of finding a solution of a singular stochastic PDE in some classical function space for the task of finding a $\gamma$-coherent germ solution of some equation set in the space of germs. Building on Lecture 3, if we are given a model $\mathsf{M} = (\Pi, g)$ on the regularity structure $\mathscr{T}_{(Eq)}$, we look for some particular germs of the form

$$\Lambda_z = \Pi_z(\boldsymbol{u}(z)) = \sum_{\tau \in \boldsymbol{B}} u_\tau(z) \Pi_z(\tau)$$

---

45. We have mainly been working in a Euclidean setting so far. We need some parabolic objects to deal with some parabolic equations. Some weighted function spaces are needed to deal with some low regularity initial condition, etc.

for a modelled distribution $\sum_{\tau \in \boldsymbol{B}} u_\tau(z)\tau$. The problem on germs then becomes the following problem on a space of modelled distributions.

One can associate to the right hand side of (52) first a modelled distribution $f(\boldsymbol{u}) \circ + g(\boldsymbol{u})$, then a germ $\{\Pi_x(f(\boldsymbol{u}(z)) + g(\boldsymbol{u}(z)))\}_z$, as in Proposition 7 of Lecture 3. If the modelled distribution is $\gamma$-regular with $\gamma > 0$ the germ is $\gamma$-coherent, so it has a unique reconstruction. This reconstruction defines the right hand side of (52). Loosely speaking, solving Equation (52) amounts to solving the equation

$$(\partial_t - D^2)u = \mathscr{R}(\{\Pi_z(f(\boldsymbol{u}(z)) \circ + g(\boldsymbol{u}(z), \partial \boldsymbol{u}(z)))\}_z)$$

with initial condition $u_0$. This is our naive attempt at making sense of Equation (52). However we cannot proceed directly in this way as the input in this formulation is a modelled distribution $\boldsymbol{u}$ and the output a function $u$. To formulate the equation as a fixed point problem we need to lift the output into a modelled distribution. This can be done but we proceed slightly differently by lifting first the integral formulation (53) of (52) as an equation on a space of modelled distributions that is consistent with the naive attempt. This leads us to the following strategy.

For a model $\mathsf{M} = (\Pi, g)$ on $\mathscr{T}_{(\mathrm{Eq})}$ and $\gamma > 0$ denote by $\mathcal{D}^\gamma_{(0,t_0]}(\mathscr{T}_{(\mathrm{Eq})}, g)$ the corresponding set of modelled distributions defined on $(0, t_0] \times \mathbb{T}^d$. Recall that we denote by $\mathscr{R}$ the reconstruction map on germs and we set for $\boldsymbol{v} \in \mathcal{D}^\gamma_{(0,t_0]}(\mathscr{T}_{(\mathrm{Eq})}, g)$

$$\mathscr{R}^\mathsf{M}(\boldsymbol{v}) \equiv \mathscr{R}(\{\Pi_z(\boldsymbol{v}(z))\}_z).$$

This is well-defined from Proposition 7 in Lecture 3.

**2. A well-posed formulation of (2) in a space of modelled distributions.**

 *(2.i) Lift of $K$ to modelled distributions.* One can associate to any admissible model $\mathsf{M} = (\Pi, g)$ on $\mathscr{T}_{(\mathrm{Eq})}$ a continuous operator

$$\mathcal{K}^\mathsf{M} \colon \mathcal{D}^\gamma_{(0,t_0]}(\mathscr{T}_{(\mathrm{Eq})}, g) \to \mathcal{D}^{\gamma+2}_{(0,t_0]}(\mathscr{T}_{(\mathrm{Eq})}, g)$$

*that is intertwined to the operator $(\partial_t - D^2)^{-1}$ via the the reconstruction operator $\mathscr{R}^\mathsf{M}$*

$$\mathscr{R}^\mathsf{M} \mathcal{K}^\mathsf{M} = (\partial_t - D^2)^{-1} \mathscr{R}^\mathsf{M}. \tag{55}$$

We emphasize that the construction of the operator $\mathcal{K}^\mathsf{M}$ requires that the model $\mathsf{M}$ be admissible. This is one of the reasons why we introduced this concept in Lecture 4. The free propagation operator $\mathcal{F}$ can be lifted[46] into an operator $\boldsymbol{\mathcal{F}}$ that sends $u_0$ into $\mathcal{D}^\gamma_{(0,t_0]}(\mathscr{T}_{(\mathrm{Eq})}, g)$ and such that

$$\mathscr{R}^\mathsf{M} \boldsymbol{\mathcal{F}} = \mathcal{F}.$$

The fixed point equation

$$\boldsymbol{u} = \boldsymbol{\mathcal{F}}(u_0) + \mathcal{K}^\mathsf{M}(f(\boldsymbol{u}) \circ + g(\boldsymbol{u}, \partial \boldsymbol{u})) \tag{56}$$

then makes sense as an equation in a space of modelled distributions. If $\boldsymbol{u}^\mathsf{M}$ is a solution of this equation, the function $u^\mathsf{M} = \mathscr{R}^\mathsf{M}(\boldsymbol{u}^\mathsf{M})$ is a solution of the equati

$$u^\mathsf{M} = \mathcal{F}(u_0) + \mathscr{R}^\mathsf{M}\mathcal{K}^\mathsf{M}(f(\boldsymbol{u}^\mathsf{M}) \circ + g(\boldsymbol{u}^\mathsf{M}, \partial \boldsymbol{u}^\mathsf{M})) = \mathcal{F}(u_0) + (\partial_t - D^2)^{-1}\mathscr{R}^\mathsf{M}(f(\boldsymbol{u}^\mathsf{M}) \circ + g(\boldsymbol{u}^\mathsf{M}, \partial \boldsymbol{u}^\mathsf{M})),$$

from the intertwining relation (54). This is consistent with the above naive attempt.

---

[46]. We simply set $\boldsymbol{\mathcal{F}}(u_0)(z) = \sum_{|k|<\gamma} u_0^{(k)}(z) X^k$, with values in the polynomial symbols.

**(2.ii) Local well-posedness of (5).** *There exists an M-dependent positive time $t_0(M) > 0$ such that Equation (4) has a unique solution $\boldsymbol{u}^M$ in $\mathcal{D}^\gamma_{(0,t_0(M)]}(\mathcal{T}_{(Eq)}, g)$.*

Pick any two models $M_1, M_2$ on $\mathcal{T}_{(Eq)}$. Although $\boldsymbol{u}_1 \in \mathcal{D}^\gamma(M_1, g)$ and $\boldsymbol{u}_2 \in \mathcal{D}^\gamma(M_2, g)$ do not live in the same space one defines a notion of distance between $\boldsymbol{u}_1$ and $\boldsymbol{u}_2$ by setting

$$\|\boldsymbol{u}_1 : \boldsymbol{u}_2\| \equiv \max_{\tau \in \boldsymbol{B}, |\tau| < \gamma} \sup_{x,y \in \mathbb{R}^d} \left( \|(\boldsymbol{u}_1 - \boldsymbol{u}_2)(x)\|_\tau + \frac{|\{\boldsymbol{u}_1(y) - \Gamma_{yx}(\boldsymbol{u}_1(x))\}_\tau - \{\boldsymbol{u}_2(y) - \Gamma_{yx}(\boldsymbol{u}_2(x))\}_\tau|}{|y-x|^{\gamma-|\tau|}} \right).$$

With this definition, one can choose $t_0$ in *(2.ii)* uniform in $M$, for $M$ in any fixed ball of the space of models, in such a way that *the map $M \to \boldsymbol{u}^M$ is continuous* on that ball.

If further $M = M^{\mathrm{can}}(\zeta)$ is the canonical model associated with a continuous function $\zeta$ then

$$u^\zeta \equiv \mathcal{R}^{M^{\mathrm{can}}(\zeta)}(\boldsymbol{u}^{M^{\mathrm{can}}(\zeta)})$$

is a classical solution of the locally in time well-posed equation

$$(\partial_t - D^2)u^\zeta = f(u^\zeta)\zeta + g(u^\zeta, \partial u^\zeta)$$

with initial condition $u_0$.

Let now $R$ be a preparation map. Denote by

$$M(R, \zeta) = (\Pi_z^{R,\zeta}, g_{yx}^{R,\zeta})_{x,y,z}$$

the model associated with $R$ and $\zeta$, and denote by $\mathcal{R}^{M(R,\zeta)}$ the reconstruction map associated with $M(R, \zeta)$. This is a map from $\mathcal{D}^\gamma(\mathcal{T}_{(Eq)}, g^{R,\zeta})$ with values in $C((0, t_0(M(R,\zeta))] \times \mathbb{T}^d)$.

**3. Interpretability of the renormalization procedure.** One can associate to $R$, $f$ and $g$ a function[47] $c_R$ such that

$$u^{R,\zeta} \equiv \mathcal{R}^{M(R,\zeta)}(\boldsymbol{u}^{M(R,\zeta)})$$

is a classical solution of the locally in time well-posed equation

$$(\partial_t - D^2)u^{R,\zeta} = f(u^{R,\zeta})\zeta + g(u^{R,\zeta}, \partial u^\zeta) + c_R(u^{R,\zeta}, \partial u^{R,\zeta})$$

with initial condition $u_0$.

So, given a random noise $\xi$ and a regularization procedure $\epsilon \to \xi^\epsilon$, can we choose some deterministic preparation maps $R^\epsilon$ so that the admissible models $M(R^\epsilon, \xi^\epsilon)$ converge? Denote by $\boldsymbol{\Pi}^{R^\epsilon, \xi^\epsilon}$ the interpretation map associated with the admissible model $M(R^\epsilon, \xi^\epsilon)$ as in Lecture 4. Recall we denote by $(\Omega, \mathcal{F}, \mathbb{P})$ the probability space on which the noise is defined.

**4. Probabilistic convergence of the BPHZ renormalized models.** *Under some mild assumptions on the law of the random noise $\xi$ there is a unique preparation maps $R^\epsilon$ ($0 < \epsilon \leq 1$) such that*[48]

$$\mathbb{E}[\boldsymbol{\Pi}^{R^\epsilon, \xi^\epsilon}(\tau)(0)] = 0 \qquad (\forall \tau \in \boldsymbol{B}, |\tau| < 0). \tag{57}$$

*The model $M^\epsilon = M(R^\epsilon, \xi^\epsilon)$ converge in any $L^p(\mathbb{P})$ space ($1 \leq p < \infty$) to a limit model.*

---

[47]. We consider the functions $f$ and $g$ as fixed, so we do not record the dependence of $c_R$ on $f$ and $g$ in the notation.

[48]. If the law of the noise $\xi$ is translation invariant one has $\mathbb{E}[\boldsymbol{\Pi}^{R^\epsilon, \xi^\epsilon}(\tau)(0)] = \mathbb{E}[\boldsymbol{\Pi}^{R^\epsilon, \xi^\epsilon}(\tau)(z)]$ for all $z$, so the evaluation point 0 has nothing special in (57).

The unique admissible model built from the preparation map $R^\epsilon$ specified by (6) is called the **BPHZ model**. It is explicit. It follows from the continuity of the maps $M \to \mathscr{R}^M$ and $M \to \boldsymbol{u}^M$ that

$$u^\epsilon \equiv \mathscr{R}^{M^\epsilon}(\boldsymbol{u}^{M^\epsilon})$$

is a continuous function of $M^\epsilon$, with values in some traditional function space, so it converges in probability to a limit function since $M^\epsilon$ is converging in $L^p(\mathbb{P})$. From point *3*, the function $u^\epsilon$ is also the solution of the equation

$$(\partial_t - D^2)u^\epsilon = f(u^\epsilon)\zeta + g(u^\epsilon, \partial u^\epsilon) + c_{R^\epsilon}(u^\epsilon, \partial u^\epsilon)$$

with initial condition $u_0$. The function $c_{R^\epsilon}$ was denoted $c^\epsilon$ in Theorem 52.

***On the benefits of continuity.*** The continuity of the solution map $M \mapsto u$ allow to transfer *automatically* a number of probabilistic properties of the random BPHZ model $M$ to the random function $u$. Hairer & Schonbauer proved for instance in [11] a support theorem[49] for $M$. It implies automatically a support theorem for the law of $u$. Hairer & Weber [13] proved a large deviation theorem for $M$. It implies automatically a large deviation theorem for $u$. I proved with Hoshino & Takano [5] a transportation cost inequality for $M$. This implies...

***Give to Caesar what belongs to Caesar.*** The fundamental notions of regularity structures and local well-posedness result of point *2* were already in Hairer's groundbreaking work [10]. Bruned, Hairer & Zambotti gave a systematic construction of the regularity structures needed for the study of singular SPDEs in [7]. The first proof of point *3* on the renormalized equation was given by Bruned, Chandra, Chevyrev & Hairer in [6]. It is also Bruned, Hairer & Zambotti who proved the existence of a unique admissible model satisfying condition (6). The first proof of convergence of the this model, point *4*, was given by Chandra & Hairer in their unpublished work [9].

Note the amuzing fact that the initials of Bruned, Hairer & Zambotti almost match the initials BPHZ used to name the model of point *4*. There is no relation! The BPHZ acronym is a reference to Bogoliubov, Parasiuk, and then Hepp and Zimmerman, who introduced some systematic renormalization procedure in quantum field theory in the mid-fifties.[50] Yet the problems that Bogoliubov & Parasiuk solved have the same origin as ours: some ill-defined products of distributions!

## 14 The road ahead

There is of course a number of technical difficulties to run that program that I do not mention here. It is important to underline a conceptual point that is involved in the items *3* and *4*.

For any model $M$, a solution $\boldsymbol{u}^M$ of Equation (4) the function $u^M = \mathscr{R}^M(\boldsymbol{u}^M)$ is a solution of the equation

$$u^M = \mathcal{F}(u_0) + \mathscr{R}^M \mathcal{K}^M(f(\boldsymbol{u}^M) \circ + g(\boldsymbol{u}^M)) = \mathcal{F}(u_0) + (\partial_t - D^2)^{-1} \mathscr{R}^M(f(\boldsymbol{u}^M) \circ + g(\boldsymbol{u}^M)).$$

While one has $\mathscr{R}^M(g(\boldsymbol{u}^M)) = g(u^M)$ for some elementary reasons, it is not clear, and generally wrong, that $\mathscr{R}^M(f(\boldsymbol{u}^M) \circ)$ is given by a function of $u^M$, as claimed in point *3*. With $M = (\Pi, g)$ and $\boldsymbol{u}^M(\cdot) = \sum_\tau u_\tau^M(\cdot)\,\tau$, the distribution $\mathscr{R}^M(f(\boldsymbol{u}^M) \circ)$ is well approximated near an arbitrary point $z$ by an expansion of the form

$$f(u^M(z))\zeta + \sum_{\tau \neq \mathbf{1}} f'(u^M(z))u_\tau^M(z)\,\Pi_z(\tau \circ) + \cdots$$

---

49. That is an exact description of the support of the law of $\boldsymbol{X}$.
50. This was Bogoliubov & Parasiuk's achievement in 1955, which were polished mathematically by Hepp (1966) and then Zimmerman (1969) a decade later.

If we had $\Pi_z(\tau_1 \cdots \tau_k \circ) = \Pi_z(\tau_1) \cdots \Pi_z(\tau_k)\zeta$ for all $\tau_1, \ldots, \tau_k$ then the expansion would be equal to $f(u^M)\zeta$, up to some negligeable quantity, and $\mathscr{R}^M(f(\boldsymbol{u}^M) \circ)$ would be equal to $f(u^M)\zeta$. However the maps $\Pi_z$ are *not multiplicative* except for the canonical model from Lecture 4, so the preceding argument and conclusion break down. This explains why $\mathscr{R}^M(f(\boldsymbol{u}^M) \circ)$ is not a priori a function of $u^M$ for an arbitrary model M.

The model $M(R, \zeta)$ that we build from a preparation map $R$ and $\zeta$ thus needs to be very special if one can indeed prove that they are interpretable in the sense that point *3* holds true. To decipher their interpretable character requires the understanding of one of two different additional algebraic structures that are implicitly involved in the construction of the set of decorated trees. The first is a multi-pre-Lie algebra structure uncovered by Bruned, Chandra, Chevyrev & Hairer in [5]. The other was uncovered by Bruned & Manchon in [8] and put to work in Bailleul & Bruned's work [1] to prove point *3*.

Understanding the mechanics behind point *3* allows to understand the reason why the choice of counterterm in the renormalized equation is not unique.[51] The fact that the counterterms can be parametrized by a finite dimensional Lie group is related to the fact that there is a group structure on the set of preparation maps.

– The convergence result of point *4* is fairly non-trivial. It was first proved in an unpublished difficult work of Chandra & Hairer [8] under some moment-type assumptions on the law of the noise. After some deep original work of Linares, Otto, Tempelmayr & Tsatsoulis [14] on an alternative to Hairer's regularity structures developped by Otto and a number of co-authors, Hairer & Steele [12] gave an alternative proof of point *4* based on some crucial insight of [14], assuming a spectral gap assumption on the law of the noise. A different alternative proof was given by Bailleul & Hoshino in [3] under the same assumption. We refer the reader to the review [4] of Bailleul & Hoshino on this subject.

We considered in these lectures the model equation (1). One can also deal with some more general equations of the form

$$(\partial_t - D^2)u = f(u)\xi + g(u) + h(u)(\nabla u, \nabla u) + k(u)\nabla u$$

for some bilinear function $h(u)(\cdot, \cdot)$ and some linear function $k(u)(\cdot)$. One can also deal with systems of singular stochastic PDEs involving possibly some different differential operators, with multi-dimensional noises whose components may have some different regularity properties, etc.

# References

[1] I. Bailleul and Y. Bruned. *Locality for singular stochastic PDEs*. arXiv:2109.00399. To appear in Ann. Probab. (2025+).

[2] I. Bailleul and M. Hoshino. *A tourist's guide to regularity structures and singular stochastic PDEs*. EMS Surv. Math. Sci. (2025).

[3] I. Bailleul and M. Hoshino. *Random models on regularity-integrability structures*. arXiv:2310.10202, (2023).

[4] I. Bailleul and M. Hoshino. *Renormalization of random models: a review*. arXiv:2409.15984, to appear in *Stoch. PDE: Anal. Comp.* (2025). Describes and compares the works [7], [9] and [3].

[5] I. Bailleul, M. Hoshino and R. Takano. *Transportation cost inequalities for singular SPDEs*. arXiv:2511.19216 (2025). It is proved here that the BPHZ model satisfies a certain transportation cost inequality if the noise satisfies also a similar inequality. A number of consequences follow.

[6] Y. Bruned, A. Chandra, I. Chevyrev and M. Hairer. *Renormalising SPDEs in regularity structures*. J. Europ. Math. Soc. **23**(11) (2021), 869–947.

---

51. You may consult the end of Section 7 of [2] for some comments on this point.

[7] Y. Bruned, M. Hairer and L. Zambotti. *Algebraic renormalisation of regularity structures*. Invent. Math. **215**(3) (2019), 1039–1156.

[8] Y. Bruned and D. Manchon. *Algebraic deformations for (S)PDEs*. J. Math. Soc. Japan **75**(2) (2023):485–526. This work describes the algebraic structures of recentering and renormalization from a unified point of view involving the algebraic concept of deformation of a pre-Lie product.

[9] A. Chandra and M. Hairer. *An analytic BPHZ theorem for regularity structures*. arXiv:1612.08138v5 (2016).

[10] M. Hairer. *A theory of regularity structures*. Invent. Math. **198**(2) (2014), 269–504.

[11] M. Hairer and P. Schonbauer. *The support of singular stochastic partial differential equations*. Forum of Mathematics, Pi (2022), Vol. 10:e1 1–127. The support of the BPHZ model is described.

[12] M. Hairer and R. Steele. *The BPHZ theorem for regularity structures via the spectral gap theorem*. Arch Rational Mech Anal. **248**(9) (2024).

[13] M. Hairer and H. Weber. *Large deviations for white-noise driven, nonlinear stochastic PDEs in two and three dimensions*. Annales de la Faculté des sciences de Toulouse : Mathématiques, Série 6, Tome 24 (2015) no. 1, pp. 55-92. A large deviation result for the BPHZ model is obtained.

[14] P. Linares, F. Otto, M. Tempelmayr and P. Tsatsoulis. *A diagram-free approach to the stochastic estimates in regularity structures*. Invent. Math. **237** (2024):1469–1565. A different algebraic and analytic setting was developed by Otto and his co-authors as an alternative to Hairer's regularity structures in a number of works. This work proves the equivalent of point *4* in their setting assuming that the law of the noise satisfies a spectral inequality.

# Appendix 1 – Basics on function spaces

## 15  Littlewood, Paley and Bony

We denote throughout by $B$ a ball of the dual space $(\mathbb{R}^d)'$ with center 0. Write $\mathscr{F}$ for the Fourier transform operator. The following result is called *Bernstein inequality*.

**Proposition 14.** *If $f \in L^p$ has its Fourier transform with support in $cB$ for some $c > 0$, and $p \leq r \leq \infty$, then*

$$\|\partial^\ell f\|_{L^r} \lesssim c^{|\ell| + d\left(\frac{1}{p} - \frac{1}{r}\right)} \|f\|_{L^p}. \tag{58}$$

**Proof.** Recall that one of Young's inequality[52] asserts the continuity of the convolution operator from $L^p \times L^q$ into $L^r$ for $\frac{1}{p} + \frac{1}{q} = 1 + \frac{1}{r}$. Pick a function $\varphi$ with compact support equal to 1 on $B$ and set $\varphi_c(x) = \varphi(x/c)$. One has from Young's inequality

$$\|\partial^\ell f\|_{L^r} = \|\partial^\ell \mathscr{F}^{-1}(\varphi_c \mathscr{F} f)\|_{L^r} = \|\partial^\ell \mathscr{F}^{-1}(\varphi_c) \star f\|_{L^r} \lesssim \|\partial^\ell \mathscr{F}^{-1}(\varphi_c)\|_{L^q} \|f\|_{L^p}$$

and the claim follows by an elementary scaling argument. □

We choose a dyadic decomposition of unit

$$1 = \rho_{-1} + \sum_{j \geq 0} \rho_j \tag{59}$$

where $\rho_{-1}$ has support in $B$ and $\rho_j(\lambda) = \rho_0(2^{-j}\lambda)$ for all $\lambda \in (\mathbb{R}^d)'$, $\mathrm{supp}(\rho_0)$ included in an annulus, and define $\Delta_j$ as the convolution operator with the inverse Fourier function of $\rho_j$, which is of the form $K_j(x) = 2^{jd} K_0(2^j x)$. The operators $\Delta_j$ send continuously $L^p$ into itself, with a norm uniformly bounded in $j$.

**Definition 15.** *For $r \in \mathbb{R}$, $p, q \in [1, \infty]$ we set*

$$\mathcal{B}^r_{pq} = \left\{ u \in \mathcal{S}'(\mathbb{R}^d) \, ; \, \|u\|_{\mathcal{B}^r_{pq}} \equiv \left( \sum_{j \geq -1} (2^{rj} \|\Delta_j(u)\|_{L^p})^q \right)^{\frac{1}{q}} < \infty \right\}$$

One has for instance $\delta_0 \in \mathcal{B}^{-d}_{\infty\infty}$ but also $\delta_0 \in \mathcal{B}^{-\epsilon - d\left(1 - \frac{1}{p}\right)}_{p\infty}$ for any $\epsilon > 0$ and $p \in [1, \infty]$.

The space $\mathcal{B}^r_{pq}$ does not depend on which partition of the identity is chosen in (59) while its norm does. This can be seen easily for $p = q = \infty$ as a consequence of two qualitative facts.

*(i)* Let $A$ be an annulus of $\mathbb{R}^d$. If $(u_n)_{n \geq -1}$ is a sequence of (necessarily) smooth functions with $\mathrm{supp}(\mathscr{F}(u_n)) \subset 2^n A$ and $\|u_n\|_\infty \lesssim 2^{-rn}$, for some $r \in \mathbb{R}$, then $u = \sum_{n \geq -1} u_n \in \mathcal{B}^r_{\infty\infty}$ with $\|u\|_{\mathcal{B}^r_{\infty\infty}} \lesssim \sup_{n \geq -1} \{2^{rn} \|u_n\|_\infty\}$.

*(ii)* Let $B'$ be any other ball of $\mathbb{R}^d$ with center 0. If $(u_n)_{n \geq -1}$ is a sequence of (necessarily) smooth functions with $\mathrm{supp}(\mathscr{F}(u_n)) \subset 2^n B'$ and $\|u_n\|_\infty \lesssim 2^{-rn}$, for some $r > 0$, then $u = \sum_{n \geq -1} u_n \in B^r_{\infty\infty}$ with $\|u\|_{B^r_{\infty\infty}} \lesssim \sup_{n \geq -1} \{2^{rn} \|u_n\|_\infty\}$.

---

[52]. A conceptual proof that involves the Riesz-Thorin interpolation theorem goes as follows. From Fubini theorem, the convolution operator with $f$ sends continuously $L^1$ into itself, with norm $\|f\|_{L^1}$. It also sends $L^\infty$ into itself with the same norm. So we have from interpolation $\|f \star g\|_{L^q} \leq \|f\|_{L^1} \|g\|_{L^q}$. This inequality also tells us as well that the convolution operator with $g$ is bounded from $L^1$ into $L^q$. On the other hand, from Hölder's inequality, this operator is also bounded from $L^{q'}$ into $L^\infty$, with $q'$ the conjugate exponent of $q$. Young's inequality follows by interpolation.

More generally, if one uses another function $\rho_0'$ instead of $\rho_0$ for defining the partition of unity, and with obvious notation, there is an integer $n_0$ such that $\rho_{j'}'$ and $\rho_j$ have disjoint supports for any $j, j'$ such that $|j'-j| \geq n_0$. It follows that $\Delta_{j'}'(u) = \sum_{|j'-j|<n_0} \Delta_{j'}'(\Delta_j(u))$, and the $j'$-uniform continuity of $\Delta_{j'}'$ on $L^p$ gives the result.

Besov embedding follow from Bernstein inequality (58). It says that for all $r \in \mathbb{R}$ and all $1 \leq p_1 < p_2 \leq \infty$ and $1 \leq q_1 \leq q_2 \leq \infty$ the space $\mathcal{B}_{p_1 q_1}^r$ is continuously embedded in $\mathcal{B}_{p_2 q_2}^{r-d\left(\frac{1}{p_1} - \frac{1}{p_2}\right)}$.

We set $\mathcal{C}^r \equiv \mathcal{B}_{\infty\infty}^r$.

*Exercice 1* – Space white noise in $\mathbb{T}^d$ is the defined formally from a sequence $(\gamma^i)_{i \geq 0}$ of independent normal random variables and an orthonormal basis $(e_i)_{i \geq 0}$ of $L^2(\mathbb{T}^d)$ by the formula

$$\xi = \sum_{i \geq 0} \gamma^i e_i.$$

We want to prove that $\xi$ is almost surely an element of $\mathcal{C}^{-\frac{d}{2}-\epsilon}$ for all $\epsilon > 0$. Why is it sufficient to see that

$$(\star) = \mathbb{E}\left[\|\xi\|_{\mathcal{B}_{2p2p}^{-d/2-\epsilon}}^{2p}\right] = \sum_{j \geq -1} 2^{-2pj\left(\frac{d}{2}+\epsilon\right)} \int \mathbb{E}[|\Delta_j \xi|^{2p}(x)] \, dx < \infty$$

for an arbitrary finite $p \geq 1$? Prove that $(\star) < \infty$ using the fact that, for a real-valued Gaussian random variable $X$, one has $\mathbb{E}[|X|^{2p}] \lesssim \mathbb{E}[|X|^2]^p$.

**Proposition 16.** *For $0 < r < 1$ the space $\mathcal{C}^r$ coincides with the space of $r$-Hölder continuous functions and the norms are equivalent.*

**Proof.** First, if $u$ is $r$-Hölder we have for $j \geq 0$ the identity

$$\Delta_j(u)(x) = \int K_j(x-y)(u(y)-u(x)) \, dy.$$

since $\int K_j = 0$. That $u \in C^r$ then follows from the inequality

$$|\Delta_j(u)(x)| \lesssim 2^{-jr} 2^{jd} \int K_0(2^j(x-y)) \, |2^j(y-x)|^r \|u\|_{r-\text{Holder}} \, dy.$$

If now $u \in \mathcal{C}^r$ we first see that $\|u\|_\infty \lesssim \sum_{j \geq -1} \|\Delta_j(u)\|_\infty \lesssim \|u\|_r$ since $r > 0$. Pick $x \neq y$ at distance less than 1 and an integer $j_0$ such that $1 \geq |y-x| \simeq 2^{-j_0}$. We estimate

$$(\star)_j(y,x) \equiv \Delta_j(u)(y) - \Delta_j(u)(x)$$

in two different ways, depending on whether $j$ is greater than or smaller than $j_0$. For $j \leq j_0$ we use Bernstein inequality (1) to write

$$|(\star)_j(y,x)| \leq \|\partial \Delta_j(u)\|_\infty |y-x| \lesssim 2^{-j(r-1)} \|u\|_r |y-x|$$

and since $r < 1$ and $|y-x| \simeq 2^{-j_0}$

$$\left|\sum_{j \leq j_0} (\star)_j(y,x)\right| \lesssim 2^{j_0(1-r)} \|u\|_r |y-x| \lesssim \|u\|_r |y-x|^r.$$

For $j > j_0$ one has $|(\star)_j(y,x)| \leq 2\|\Delta_j(u)\|_\infty \lesssim 2^{-jr}\|u\|_r$, so

$$\sum_{j > j_0} |(\star)_j(y,x)| \lesssim 2^{-j_0 r}\|u\|_r \simeq |y-x|^r \|u\|_r. \qquad \square$$

We warn the reader that the spaces $\mathcal{C}^r$ do not coincide with the classical $C^r$ spaces when $r \in \mathbb{N}$. The space $\mathcal{C}^0$ is in particular slightly bigger than $L^\infty$.

For any $u, v \in \mathcal{S}'(\mathbb{R}^d)$ set

$$P_u(v) \equiv \sum_{i+2 \leq j} \Delta_i(u)\, \Delta_j(v)$$

and

$$\Pi(u,v) = \sum_{|j-i| \leq 1} \Delta_i(u)\, \Delta_j(v).$$

We have Bony's formal decomposition

$$uv = P_u(v) + \Pi(u,v) + P_v(u).$$

Here is an elementary remark.

**Remark 17.** There is an annulus $A$ such that for all distributions $u,v$ and all indices with $i+2 \leq j$ one has $\operatorname{supp} \mathscr{F}(\Delta_i(u)\Delta_j(v)) \subset 2^j A$. There is also a centered ball $B'$ such that for all distributions $u,v$ and all indices with $|j-i| \leq 1$ one has $\operatorname{supp} \mathscr{F}(\Delta_i(u)\Delta_j(v)) \subset 2^j B'$.

The following continuity estimates due to Bony follow as a consequence of Remark 4.

**Proposition 18.** *For every $r_1, r_2 \in \mathbb{R}$ and $u_1 \in \mathcal{C}^{r_1}$ and $u_2 \in \mathcal{C}^{r_2}$ we have*

$$\|P_{u_1}(u_2)\|_{r_2 + \min(0, r_1)} \lesssim \|u_1\|_{\min(0, r_1)} \|u_2\|_{r_2}.$$

*We further have*

$$\|\Pi(u_1, u_2)\|_{r_1 + r_2} \lesssim \|u_1\|_{r_1} \|u_2\|_{r_2} \qquad (60)$$

*if and only if if $r_1 + r_2 > 0$.*

**Proof.** Set $\Delta_{\leq j-2} \equiv \sum_{-1 \leq i \leq j-2}$. From Remark 17, there is an annulus $A$ such that $\operatorname{supp} \mathscr{F}(\Delta_{\leq j-2}(u_1)\Delta_j(u_2)) \subset 2^j A$, and $\|\Delta_i(u_1)\Delta_j(u_2)\|_\infty$ is bounded above by $\|u_1\|_\infty 2^{-jr_2}\|u_2\|_{r_2}$ if $r_1 \geq 0$, otherwise it is bounded above by $2^{-jr_1}\|u_1\|_{r_1} 2^{-jr_2}\|u_2\|_{r_2}$. In both situations we can apply point *(i)* above. From Remark 5, for $|i-j| \leq 1$ there is a ball $B'$ such that $\operatorname{supp} \mathscr{F}(\Delta_i(u_1)\Delta_j(u_2)) \subset 2^j B'$, while $\|\Delta_i(u_1)\Delta_j(u_2)\|_\infty \leq 2^{-j(r_1+r_2)}\|u_1\|_{r_1}\|u_2\|_{r_2}$. Point *(ii)* above gives (60) if $r_1 + r_2 > 0$. To see that (3) is sharp, consider $u_1^{(n)}(x) = n^{-r_1} e^{inx}$ and $u_2^{(n)}(x) = n^{-r_2} e^{-inx}$. One has for any $r_1' < r_1$ and $r_2' < r_2$ that $u_1^{(n)}$ is converging to 0 in $\mathcal{C}^{r_1'}$ and $u_2^{(n)}$ is converging to 0 in $\mathcal{C}^{r_2'}$, while $u_1^{(n)} u_2^{(n)} = n^{-r_1 - r_2}$ diverges when $r_1 + r_2 < 0$ and is constant when $r_1 + r_2 = 0$. We do not give the proof of the non-continuity of the resonant operator $\Pi$ when $r_1 + r_2 = 0$; we refer for that case to Theorem 2.52 of Bahouri, Chemin & Danchin's textbook [1]. $\square$

**Bony's product rule** states that *the product map $(u_1, u_2) \in \mathcal{C}^{r_1} \times \mathcal{C}^{r_2} \to u_1 u_2 \in \mathcal{C}^{\min(r_1, r_2)}$ is well-posed and continuous if and only if $r_1 + r_2 > 0$.*

*Exercice 2* – Prove that the sequence of functions $x \to \sin(nx)$ converges to 0 in $\mathcal{D}'(\mathbb{R})$ but the sequence of functions $x \to \sin^2(nx)$ does not converge to 0 in $\mathcal{D}'(\mathbb{R})$.

Write $f_x^\lambda(\cdot) = \lambda^{-d} f\left(\frac{\cdot - x}{\lambda}\right)$ and $f^\lambda(\cdot) = f_0^\lambda(\cdot)$ for any function $f: \mathbb{R}^d \to \mathbb{R}$ and $\lambda > 0$. We say the function $f_x^\lambda$ is localized at point $x$ at scale $\lambda$.

The Besov-Hölder spaces $B_{\infty\infty}^r$ are characterized by the $L^\infty$ norm of $\lambda_j^{-r} K_j \star u = \lambda_j^{-r}(K_0)^{\lambda_j} \star u$ with $\lambda_j = 2^{-j}$. The precise definition of the function $K_0$ and the scaling parameter $2^{-j}$ are not so important. As a matter of fact one can give the following equivalent definition of $\mathcal{C}^r$ in which $\mathcal{E}^n$ stands for the classical set of $C^n$ functions with support in the unit ball of $\mathbb{R}^d$ and $C^n$-norm at most 1.

**Proposition 19.** *The space $\mathcal{C}^r$ coincides, for any $r \in \mathbb{R}$, with the set of distributions $\Lambda$ on $\mathbb{R}^d$ such that*

$$|\Lambda(f_x^\lambda)| \lesssim \lambda^r$$

*uniformly in $0 < \lambda \leq 1$, $x \in \mathbb{R}^d$ and $f \in \mathcal{E}^{[|r|]+1}$ such that $\int fP = 0$ for all polynomials $P$ of degree less than $[|r|]+1$. The smallest implicit constant in $\lesssim$ defines a norm on $\mathcal{C}^r$ equivalent to the norm $\|\cdot\|_r$.*

Replacing $[|r|]+1$ by any other integer greater than $[|r|]+1$ defines on $\mathcal{C}^r$ a norm equivalent to the norm $\|\cdot\|_r$. Furthermore, for $r > 0$ non-integer, the space $\mathcal{C}^r$ coincides with the space of $r$-Hölder continuous functions and the norms are equivalent.

## 16 The parabolic world

The parabolic space $\mathbb{R} \times \mathbb{R}^{d-1}$ has norm $|(t,x)| \equiv \sqrt{|t|} + |x|$ and parabolic dilations $\delta_c(t,x) \equiv (c^2 t, cx)$ for all $c \in \mathbb{R}$. For a multi-index $\ell = (\ell_0, \ell_1, \ldots, \ell_{d-1}) \in \mathbb{N}^d$ we set $|\ell|_\mathfrak{p} \equiv 2\ell_0 + \ell_1 + \cdots + \ell_{d-1}$ and $|\mathfrak{p}| \equiv 2 + (d-1) = d+1$.

Bernstein's estimate in Proposition 14 is still true if we replace $cB$ by $\delta_c(B)$ and $|\ell|$ in (1) by $|\ell|_\mathfrak{p}$. We define a parabolic partition of unit by setting $\rho_j(\lambda) = \rho_0(\delta_{2^j}(\lambda))$. This function has Fourier inverse of the form $K_j(z) = (2^j)^{|\mathfrak{p}|} K_0(\delta_{2^j}(z))$, for any $z = (t,x)$. We also write $\Delta_j$ for the convolution operator with the previous kernel $K_j$. The Besov spaces are defined as in Definition 2 from the operators $\Delta_j$. Spacetime white noise is almost surely in the parabolic Hölder space of exponent $-\frac{|\mathfrak{p}|}{2} - \epsilon$ for any $\epsilon > 0$.

The rest of Section 15 has a direct parabolic counterpart where Euclidean dilations are turned into parabolic dilations, with for instance $f_z^\lambda(z') = \lambda^{-|\mathfrak{p}|} f(\delta_\lambda(z'-z))$.

We refer the reader to the excellent Section 2 of Jörg Martin's PhD dissertation [3] for a detailed presentation of the above elementary facts in a setting with a generic scaling.

## 17 Schauder estimate

We prove in this section a version of the classical Schauder estimate, following Friz & Hairer's account in Section 14.3 of [2].

To have some unified notations we set $|\mathfrak{p}| = d$ and $|x|_\mathfrak{p} = |x|$ in the Euclidean setting, and $|\mathfrak{p}| = d+1$ in the parabolic setting, with $|z|_\mathfrak{p}$ defined as above.

A smooth function $K: \mathbb{R}^d \setminus \{0\} \to \mathbb{R}$ with compact support is said to be **2-regularizing** if there exists for any $\ell \in \mathbb{N}^d$ a constant $C_\ell > 0$ such that

$$|(\partial^\ell K)(z)| \leq C_\ell |z|_\mathfrak{p}^{2-|\mathfrak{p}|-|\ell|}.$$

Any version of the heat kernel $K(t,x) = t^{-(d-1)/2}\exp\left(-\frac{|x|^2}{4t}\right)\mathbf{1}_{t>0}$ smoothly localized in a compact set is 2-regularizing for instance.

Write $\mathcal{E}^n_\lambda$ for the set $\{f^\lambda_x ; f \in \mathcal{E}^n, x \in \mathbb{R}^d\}$. The following fact is reminiscent of the Littlewood-Paley decomposition.

**Remark 20.** *Pick $r \geq 0$. A 2-regularizing kernel $K$ can be decomposed as*

$$K = \sum_{i \geq -1} K_{(i)}, \quad \text{with} \quad 2^{2i} K_{(i)} \in \{c\mathcal{E}^n_{2^{-i}}\}$$

for some constant $c > 0$. Indeed using the partition of unit (2) set $K_{(i)}(z) \equiv \rho_{i+1}(|z|_\mathfrak{p})K(z)$. It satisfies $|(\partial^j K_{(i)})(z)| \leq (2^{-i})^{2-|\mathfrak{p}|-|j|}$ for all $|j| \leq n$.

The following elementary fact says that the convolution of two localized functions is localized at scale the sum of the scales: *There is a constant $C > 0$ such that for all $f, g \in \mathcal{E}^n$ we have $(f^\lambda \star g^\mu)^{1/(\lambda+\mu)} \in \{c\mathcal{E}^n\}$.* Prove it. This result can be refined as follows.

**Lemma 21.** *Assume $\lambda \leq \mu$. There is some constant $c > 0$ with the following property. For all $f, g \in \mathcal{E}^n$ such that $\int Pg = 0$ for all polynomials of degree smaller than some $\gamma > 0$, we have $(f^\lambda \star g^\mu)^{1/(2\mu)} \in \{c(\lambda/\mu)^\gamma \mathcal{E}^{[n-\gamma]}\}$.*

**Proof.** We need to estimate $\partial^\ell (f^\lambda \star g^\mu)$ for all $|\ell| \leq r - \gamma$. We have

$$\partial^\ell (f^\lambda \star g^\mu)(z) = \int \partial^\ell g^\mu(z-z') f^\lambda(z') \, dz' = \int \left(\partial^\ell g^\mu(z-z') - P_z^{\ell,\lambda,\gamma}(z-z')\right) f^\lambda(z') \, dz'$$

where $P_z^{\gamma,\ell}$ stands for the Taylor expansion at point $z$ of the function $\partial^\ell g^\mu$ at (Euclidean or parabolic) order $\gamma$. The corresponding Taylor remainder has size at most $\|g^\mu\|_{\gamma+|\ell|}|z'|^\gamma \lesssim \mu^{-d-\gamma-|\ell|}|z'|^\gamma$. We obtain the factor $\lambda^\gamma$ from $\int |z'|^\gamma |f^\lambda(z')| \, dz'$ by a change of variable. $\square$

Here is Schauder continuity estimate.

**Theorem 22.** *For any $r \in \mathbb{R}$ the map $u \mapsto K \star u$ sends continuously $\mathcal{C}^r$ into $\mathcal{C}^{r+2}$.*

**Proof.** We use the characterization of the Besov-Hölder spaces from Proposition 19. Set $n_r \equiv 2(|r|+2)+2$. Pick $u \in \mathcal{C}^r$ and $f \in \mathcal{E}^{n_r}_\lambda$ such that $\int Pf = 0$ for all polynomials of (Euclidean or parabolic) degree less than $n_r + 2$. Set $\overline{K}(z) = K(-z)$. We have

$$(K \star u)(f) = u(\overline{K} \star f) = \sum_{i \geq -1} u(\overline{K_{(i)}} \star f) = \sum_{i \geq -1} 2^{-2i} u(2^{2i} \overline{K_{(i)}} \star f). \tag{61}$$

It follows from Lemma 21 that $|u(2^{2i}\overline{K_{(i)}} \star f)|$ is bounded above by $\lambda^r$ if $2^{-i} \leq \lambda$, and by $(2^i \lambda) 2^{-ir}$ otherwise. We see as a consequence that $|(K \star u)(f)| \lesssim \lambda^{r+2}$ by splitting the sum (61) in two sums, depending on whether $2^{-i} \leq \lambda$ or not. $\square$

## References

[1] H. Bahouri, J.-Y. Chemin and R. Danchin. *Fourier Analysis and Nonlinear Partial Differential Equations*. Grundlehren der mathematischen Wissenschaften, vol. 343, Springer (2012). *A nice textbook on Fourier analysis and its applications to some equations from fluid mechanics*.

[2] P. Friz and M. Hairer. *A course on rough paths. With an introduction to regularity structures*. Universitext, Springer (2020).

[3] J. Martin. *Refinements of the Solution Theory for Singular SPDEs*. PhD dissertation. Available at https://edoc.hu-berlin.de/server/api/core/bitstreams/fc81fc3b-d387-4601-abab-b1883bb0d72d/content.

# Appendix 2 – Sewing, rough paths and their dynamics

The model situation for the study of singular stochastic PDEs is the study of stochastic differential equations $dz_t = f(z_t)dX_t$, driven by some random controls of low regularity, like Brownian motion. Ito's stochastic calculus can be used for semimartingale controls but there are situations where this setting is not adequate from a modelization point of view. T. Lyons's theory of rough paths [5] provides a framework for the study of such equations. The pathwise theories of singular stochastic PDEs are modelled on the insights of rough paths theory. This appendix is a light introduction to this setting. The parallel with what we have seen about singular stochastic PDEs will be clear: The reconstruction theorem takes the form of the sewing lemma, models are now rough paths, the notion of modelled distribution is now the notion of controlled path, and there is a notion of integral with respect to a rough path that bypasses the product problem involved in $f(z_t)dX_t$.

Although we will not see the details here, it is noticeable that the semimartingale drivers have a natural lift into some rough paths whose construction does not require any renormalization procedure. This allows to relate the pathwise and Itô approach in a natural way. Renormalization may be needed for some wilder random controls outside of the realm of Itô calculus. We note as well that unlike what happens with singular stochastic PDEs, where one associates to each equation a regularity structure, the algebraic structure involved in the study of controlled ordinary differential equations does not depend on which equation we look at. We will not touch here on that side of rough paths theory.

We follow closely Zorin-Kranich's nice account [6] of the basics of rough paths theory throughout this appendix. We also refer the reader to Chapter 16 of the amazing book [3] of Feyel & de la Pradelle. You can read the lecture notes [1] for an introduction to the general theory of rough differential equations.

## 18  The sewing lemma

The reconstruction theorem for germs stated in Theorem 2 of Lecture 2 has its roots in Gubinelli's sewing lemma [4]. The form given here is due to Feyel & de la Pradelle [2].

**Theorem 23 (Sewing lemma)** – *Let $(E, |\cdot|)$ be a Banach space and $(\mu_{ts})_{0 \le s \le t \le 1}$ be a family of elements of $E$ such that*

$$|\mu_{ts} - (\mu_{tu} + \mu_{us})| \le c_1 |t-s|^\theta \tag{62}$$

*for all $0 \le s \le t \le 1$, for some constants $c_1 \in (0, +\infty)$ and $\theta \in (1, +\infty)$. There exists a unique path $(\varphi_t)_{0 \le t \le 1}$ with $\varphi_0 = 0$ such that*

$$|\varphi_t - \varphi_s - \mu_{ts}| \lesssim |t-s|^\theta. \tag{63}$$

We can talk of $(\mu_{ts})_{0 \le s \le t \le 1}$ as a germ and talk of the function $\varphi$ as its reconstruction. A **function** $\mu$ which satisfies (62) is said to be **almost additive**.

**Proof** – The proof is constructive and constructs the increments $\varphi_t - \varphi_s$ of the path $\varphi$ as a limit of the $\mu_{ts}^\pi = \sum_{i=0}^{n-1} \mu_{t_i t_{i-1}}$, for a partition $\pi = \{s = t_0 < t_1 < \cdots < t_{n-1} < t_n = t\}$ of $[s,t]$, as the mesh of the partition goes to 0.

**Step 1.** We fix $0 \le s < t \le 1$ and prove by induction on $n$ that

$$|\mu_{ts} - \mu_{ts}^\pi| \le c_1 \left( \sum_{k=1}^{n-1} (2/k)^\theta \right) |t-s|^\theta. \tag{64}$$

One has $\mu_\pi = \mu_{ts}$ if $n=1$, otherwise assume the result is proved for all partitions with $n$ sub-intervals and take a partition $\pi$ with $n+1$ sub-intervals. Since

$$\sum_{j=0}^{n-1}(t_{j+2}-t_j) = \sum_{j\leq n-1,\text{even}}(t_{j+2}-t_j) + \sum_{j\leq n-1,\text{odd}}(t_{j+2}-t_j) \leq 2(t-s) \tag{65}$$

there is an index $0\leq k\leq n-1$ such that $t_{k+2}-t_k \leq 2(t-s)/n$. The partition $\pi'=\pi\setminus\{t_k\}$ has $n$ sub-intervals and

$$|\mu_{ts}^\pi - \mu_{ts}^{\pi'}| = |(\mu_{t_{k+2}t_{k+1}} + \mu_{t_{k+1}t_k}) - \mu_{t_{k+2}t_k}| \leq c_1 |t_{k+2}-t_k|^\theta \leq c_1(2/n)^\theta |t-s|^\theta$$

so (64) follows from applying the induction hypothesis to $\pi'$.

**Step 2.** For any finite partitions $\pi_1 \subset \pi_2$ of $[s,t]$ with $\pi_1 = \{s=t_0<t_1<\cdots<t_{n-1}<t_n=t\}$ one can apply in each interval $[t_i,t_{i+1}]$ the result of step 1 with partition of $[t_i,t_{i+1}]$ induced by $\pi_2$. This gives

$$|\mu_{ts}^{\pi_1} - \mu_{ts}^{\pi_2}| \leq \sum |\mu_{t_{i+1}t_i} - \mu_{t_{i+1}t_i}^{\pi_2}| \lesssim \sum |t_{i+1}-t_i|^\theta \lesssim o_{t-s}(1)(t-s),$$

which shows that the family $(\mu_{ts}^\pi)_\pi$ is a Cauchy net in the Banach space $E$. We denote by $\bar\varphi_{ts} \in E$ its limit. It follows from this strong convergence result that for any $0\leq s\leq u\leq t\leq 1$ we have $\bar\varphi_{ts} = \bar\varphi_{tu} + \bar\varphi_{us}$, because we have $\mu_{ts}^\pi = \mu_{tu}^\pi + \mu_{us}^\pi$ for all partitions of $[s,t]$ that have $u$ as one of their points. This additivity property of the 2-variable function $\bar\varphi$ implies that $\bar\varphi_{ts} = \bar\varphi_{t0} - \bar\varphi_{s0}$. The function $\varphi_r = \bar\varphi_{r0}$ satisfies the condition (63). To see the uniqueness of such a function, writing $\psi_t - \psi_s = \sum(\psi_{t_{i+1}} - \psi_{t_i})$ for $\pi=\{t_i\}$, we see that any function $\psi$ which satisfies (63) also satisfies

$$|\psi_t - \psi_s - \mu_{ts}^\pi| \lesssim \sum |t_{i+1}-t_i|^\theta \lesssim \max(t_{i+1}-t_i)^{\theta-1}. \tag{66}$$

This concludes the proof of the theorem. $\square$

*Example (Young integral)* – Take $1\leq p<2$. For $a,b:[0,1]\to\mathbb{R}$ two $(1/p)$-Hölder functions set

$$\mu_{ts} = a_s(b_t - b_s).$$

One checks that one has for all $0\leq s\leq u\leq t\leq 1$ the identity

$$\mu_{ts} - (\mu_{tu} + \mu_{us}) = (a_s - a_t)(b_u - b_t)$$

so $|\mu_{ts} - (\mu_{tu}+\mu_{us})| \leq |t-s|^{2/p}$ with $2/p > 1$. We can set $\int_0^t a\,db \equiv \varphi_t$ for the function $\varphi$ of Theorem 1. The $(\mu_{ts})$ are not almost additive anymore when $2\leq p$. One can still define $\int_0^t a\,db$ if $b$ is replaced by a richer object and we have more information on $a$ than its mere regularity. This is what we do in the next section. We make a useful remark before.

**Remark (A cheap generalization)** – The only property of the function $\omega(s,t) = |t-s|$ that we used in the above proof is the fact that it is continuous, null on the diagonal and superadditive:

$$\omega(s,u) + \omega(u,t) \leq \omega(s,t) \qquad (\forall s\leq u\leq t).$$

These three properties give the definition of a *control*. If we replace in the proof $|t-s|$ by $\omega(s,t)$ for a control $\omega$, the inequality (65) now takes the form

$$\sum_{j\leq n-1,\text{even}} \omega(t_j,t_{j+2}) + \sum_{j\leq n-1,\text{odd}} \omega(t_j,t_{j+2}) \leq 2\omega(s,t),$$

which holds from the subadditivity of the function $\omega$. The continuity of $\omega$ is used in Equation (66), which reads now

$$\sum \omega(t_{i+1} - t_i)^\theta \lesssim \max \omega(t_i, t_{i+1})^{\theta-1} \sum \omega(t_i, t_{i+1}) \leq \max \omega(t_i, t_{i+1})^{\theta-1} \omega(s,t)$$

here again from subadditivity in the last inequality. So Theorem 1 holds if we replace $|t-s|$ in (62) and (63) by $\omega(s,t)$ for a continuous control $\omega$.

A typical example of subadditive function is given by the $q$-variation $\omega_{q,g}$ of a function $g$, defined by

$$\omega_{q,g}(s,t) \equiv \|g\|_{q-\text{var},[s,t]} \equiv \|g\|_{q,(s,t)} \equiv \left( \sup_\pi \sum |g_{t_{i+1} t_i}|^q \right)^{1/q}$$

where we write $g_{t_{i+1} t_i} \equiv g_{t_{i+1}} - g_{t_i}$, for a supremum over all the finite partitions $\pi$ of $[s,t]$. This supremum is finite if for instance $g$ is $(1/q)$-Hölder, in which case $\|g\|_{q,(s,t)}$ is smaller than the $(1/q)$-Hölder norm of $g$.

## 19 Rough paths and controlled paths

Below, it is convenient to use the notation $(\mathbb{R}^d)^{\otimes 2}$ for $L(\mathbb{R}^d, \mathbb{R}^d)$ and to identify the matrix $(a^i b^j)_{1 \leq i,j \leq d}$ with the element $a \otimes b$, for $a, b \in \mathbb{R}^d$. For those of you who are not familiar with notation, this will be just a book-keeping device for us here.[53]

• Consider a controlled ordinary differential equation

$$z_t = x_0 + \int_0^t f(z_s) dh_s = x_0 + \int_0^t f(z_s) \dot{h}_s \, ds$$

on $\mathbb{R}^n$ driven by a $C^1$ control $h$ with values in $\mathbb{R}^d$ and a one form $f \colon \mathbb{R}^n \to L(\mathbb{R}^d, \mathbb{R}^n)$ of class $C_b^3$, say. It has global solutions and the solution paths satisfy[54] for any $0 \leq s < t$

$$z_t = z_s + f(x_x)(h_t - h_s) + (f'f)(x_s) \left\{ \int_s^t (h_r - h_s) \otimes h'_r \, dr \right\} + O(|t-s|^3).$$

We note that $|h_t - h_s| \lesssim |t-s|$ and $|\int_s^t (h_r - h_s) \otimes h'_r \, dr| \lesssim |t-s|^2$. The notion of rough path that we now introduce plays the role of the coefficients $(h_t - h_s, \int_s^t (h_r - h_s) \otimes h'_r \, dr)$ in the above local expansion and allows to give a meaning to the controlled differential equation when the control $h$ is not $C^1$ anymore.

• Fix $2 < p < 3$. A ***p-rough path*** over the time interval $[0,1]$ is a continuous map

$$\boldsymbol{X} \colon \{0 \leq s \leq t \leq 1\} \mapsto (X_{st}, \mathbb{X}_{st}) \in \mathbb{R}^d \times ((\mathbb{R}^d)^{\otimes 2})$$

such that one has Chen's relations for all $s \leq u \leq t$

$$X_{st} = X_{su} + X_{ut}, \qquad \mathbb{X}_{st} = \mathbb{X}_{su} + \mathbb{X}_{ut} + X_{su} \otimes X_{ut}$$

and $X$ has finite $p$-variation and $\mathbb{X}$ has finite $(p/2)$-variation $\|\mathbb{X}\|_{p/2} \equiv \sup_\pi \sum |\mathbb{X}_{t_{i+1} t_i}|^p$ for a supremum over all the finite partitions $\pi$ of $[0,1]$.[55] This holds for instance if $X$ is $(1/p)$-Hölder and $|\mathbb{X}_{st}| \lesssim |t-s|^{2/p}$.

---

53. This is more than that in a general setting, but we do not care here.
54. For $f = (f_1, \ldots, f_d)$ and each $f_i$ with values in $\mathbb{R}^n$, for $a, b \in \mathbb{R}^d$, we set $(f'f)(a \otimes b) = \sum (df_i)(f_j) a^i b^j$.
55. A $p$-rough path can be seen to be a model over a particular regularity structure; we will not use that fact here.

*Exercice* – Take $h \in C_b^1$. Check that $(h_t - h_s, \int_s^t (h_r - h_s) \otimes h_r' dr)$ satisfies Chen's relations.

We denote by $\|\mathbb{X}\|_{p/2,(s,t)}$ the $(p/2)$-variation of $\mathbb{X}$ on a time interval $[s,t]$, and we define a continuous control setting

$$\omega_{\boldsymbol{X}}(s,t) \equiv \|X\|_{p,[s,t]} + \|\mathbb{X}\|_{p/2,[s,t]}.$$

By a solution of the rough differential equation

$$dz_t = f(z_t) dX_t$$

we mean a path $z$ that satisfies everywhere the local expansion property

$$z_t = \{z_s + f(x_x)X_{st} + (f'f)(x_s)\mathbb{X}_{st}\} + O(\omega_{\boldsymbol{X}}(s,t)^{3/p}). \tag{67}$$

It turns out to be useful to relate this property to an integral formulation of the equation

$$z_t = x_0 + \int_0^t f(z_t) dX_t \tag{68}$$

to be defined, via the introduction of the following notion.

• An $(E \times L(\mathbb{R}^d, E))$-valued **path** $\boldsymbol{a} = (a_t, a_t')_{0 \leq t \leq 1}$ is said to be **controlled by** $X$ if the path $a'$ is continuous and has finite $p$-variation and

$$R_{st}^{\boldsymbol{a}} \equiv a_t - a_s - a_s' X_{st}$$

has finite $(p/2)$-variation. We define the size of $\boldsymbol{a}$ as $\|a'\|_p + \|R^{\boldsymbol{a}}\|_{p/2}$.

**Lemma 24 (Rough integral)** – *Let $(a_t)_{0 \leq t \leq 1}$ be an $L(\mathbb{R}^d, \mathbb{R}^n)$-valued path controlled by $X$. The $\mathbb{R}^n$-valued function $\mu_{ts} = a_s X_{st} + a_s' \mathbb{X}_{st}$ is almost additive.*

We denote by $t \to \int_0^t \boldsymbol{a} \, d\boldsymbol{X}$, with bold $\boldsymbol{a}$ and $\boldsymbol{X}$, the function $\varphi$ associated to $\mu$ by Theorem 23.

**Proof** – We check that one has

$$\begin{aligned}
\mu_{ts} - (\mu_{tu} + \mu_{us}) &= (a_s - a_u)X_{ut} + a_s'(\mathbb{X}_{st} - \mathbb{X}_{su}) - a_t'\mathbb{X}_{ut} \\
&= (a - a_u)X_{ut} + a_s'(\mathbb{X}_{ut} + X_{su} \otimes X_{ut}) - a_t'\mathbb{X}_{ut} \\
&= (a_s + a_s' X_{su} - a_u)X_{ut} + (a_s' - a_t')\mathbb{X}_{ut}
\end{aligned}$$

using Chen's relation for $\mathbb{X}$ in the second equality. We define a control setting

$$\omega = \omega_{\boldsymbol{X}} + \omega_{p,a'}.$$

In the last display we have $|a_s + a_s' X_{su} - a_u| \leq \omega_{\boldsymbol{X}}(s,u)^{2/p}$ and $|X_{ut}| \leq \omega_{\boldsymbol{X}}(u,t)^{1/p}$ and $|(a_s' - a_t')| \leq \omega_{p,a'}(s,t)^{1/p}$ and $|\mathbb{X}_{ut}| \leq \omega_{\boldsymbol{X}}(u,t)^{2/p}$. This shows that $|\mu_{ts} - (\mu_{tu} + \mu_{us})| \leq \omega(s,t)^{3/p}$; we are done. □

**Corollary 25 (On the integration map)** – *The path $\boldsymbol{\varphi} = (\varphi, \varphi')$ with $\varphi' = a$ is controlled by $X$ and we have*

$$\|R^{\boldsymbol{\varphi}}\|_{p/2} \lesssim \|R^{\boldsymbol{a}}\|_{p/2} \|X\|_p + \|a'\|_p \|\mathbb{X}\|_{p/2} + \|a'\|_\infty \|\mathbb{X}\|_{p/2}.$$

**Proof** – One has $R_{st}^{\boldsymbol{\varphi}} = \varphi_t - \varphi_s - a_s X_{st} = (\varphi_t - \varphi_s - a_s X_{st} - a_s' \mathbb{X}_{st}) + a_s' \mathbb{X}_{st}$, so the bound comes from the version of the proof of Lemma 2. □

One obtains the following statement from a direct application of Taylor formula – see Proposition 4.13 in Section 4 of [2] for a proof. We denote below by $\|\cdot\|_{\sup}$ the supremum norm of a bounded function.

**Lemma (On the composition of controlled paths with smooth functions)** – *Take $f \in C_b^2$ and a path $\boldsymbol{a} = (a_t, a'_t)_{0 \leq t \leq 1}$ controlled by $X$. Then*

$$f(\boldsymbol{a}) = (f(a_t), f'(a_t)a'_t)_{0 \leq t \leq 1}$$

*is controlled by $X$ and*

$$\|f'(a)a'\|_p \leq \|f\|_{C^2}(\|a'\|_p + \|a\|_p \|a'\|_{\sup})$$

*and*

$$\|R^{f(\boldsymbol{a})}\|_{p/2} \leq \|f\|_{C^2}(\|R^{\boldsymbol{a}}\|_{p/2} + \|a\|_p^2).$$

It follows from these facts that the formula

$$\mathbb{I} \colon \boldsymbol{z} = (z, z') \to \left( x_0 + \int_0^\cdot f(\boldsymbol{z}_t) d\boldsymbol{X}_t \,, f(z) \right)$$

defines a map from the space of $\mathbb{R}^n$-valued paths controlled by $X$ into itself. A solution of the differential equation (7) is defined as a fixed point of the map $\mathbb{I}$. It can be seen, and this is elementary, that the definition of a solution to the differential equation (6) as a fixed point of the map $\mathbb{I}$ is equivalent to the definition of a solution as a path satisfying the local expansion property (6) everywhere.

## 20 Rough differential equations

We prove that the map $\mathbb{I}$ is a contraction of a large enough ball of the space of paths contolled by $X$, provided we work on a small enough $\boldsymbol{X}$-dependent time interval. The existence time turns out to be independent of the initial condition so we have globall well-posedness. As usual we need some stability results to make the fixed point formulation efficient.

### 20.1 Stability results for the integration and composition maps

One can quantify the continuity of the integration map $\boldsymbol{a} \to \boldsymbol{\varphi}$ of Lemma 2 as follows. For any two functions $g^1$ and $g^2$ we write below $\Delta g$ for $g^1 - g^2$.

**Lemma (Stability of the integration map)** – *Let $\boldsymbol{a}^1, \boldsymbol{a}^2$ be controlled by $X$, with $\boldsymbol{\varphi}^1$ and $\boldsymbol{\varphi}^2$ their associated integrals. One has*

$$\|R^{\boldsymbol{\varphi}^1} - R^{\boldsymbol{\varphi}^2}\|_{p/2} \lesssim (\|\Delta a'\|_{\sup} + \|\Delta a'\|_p + \|\Delta R^{\boldsymbol{a}}\|_{p/2})\,\omega_{\boldsymbol{X}}$$

Here the variation norms are computed on some arbitrary interval $[s,t]$ and $\omega_{\boldsymbol{X}} = \omega_{\boldsymbol{X}}(s,t)$. Similarly one can quantify the continuity of the composition map $\boldsymbol{a} \to f(\boldsymbol{a})$ as follows.

**Lemma 26 (Stability of the composition map)** – *For $\boldsymbol{a}^1, \boldsymbol{a}^2$ controlled by $X$ one has*

$$\|\Delta f(a)\|_p \lesssim_f \|\Delta a\|_p + \|\Delta a\|_{\sup} \|a^2\|_p$$

*and*

$$\|\Delta f'(a)a'\|_p \lesssim_f \|\Delta a\|_p + \|\Delta a'\|_{\sup}\|a^1\|_p + \|\Delta a\|_{\sup}\|a^{2'}\|_p + \|\Delta a'\|_p \|a^{2'}\|_{\sup}$$
$$+ \|\Delta a\|_{\sup} \|a^2\|_p \|a^2\|_{\sup}$$

*and*

$$\|\Delta R^{\boldsymbol{\varphi}}\|_{p/2} \lesssim_f \|\Delta R^a\|_{p/2} + \|\Delta a\|_{\sup}(\|R^{a^2}\|_{p/2} + \|a^2\|_p^2) + (\|a^1\|_p + \|a^2\|_p)\|\Delta a\|_p \,.$$

We refer the reader to Lemma 4.16 and Lemma 4.15 in Section 4 of [2] for a proof of these two stability results.

## 20.2 Fixed point of $\mathbb{I}$

We use the stability results of Section 20.1 to prove the contractive character of the map $\mathbb{I}$, with $\omega_{\boldsymbol{X}}(0,T)$ used as a small parameter for an appropriate ($\boldsymbol{X}$-dependent) choice of horizon $T$. One needs to quantify the dependence with respect to $\boldsymbol{X}$ of all the operations involved in the fixed point problem to prove that the unique solution path to Equation (67) is a continuous function of $\boldsymbol{X}$, with the solution path seen as a path of finite $p$-variation. We do not give these details here[56] but give now some of the details of the fixed point problem.

Denote by $c_f = \|f\|_{\sup} + \|f\|_{\operatorname{Lip}}$ the $C^1$ norm of $f$. We define a ball on the space of controlled paths

$$\mathcal{B}_T(m) = \{\boldsymbol{z} = (z,z'); \|z\|_p \leq m/c_f, \|R^z\|_{p/2} \leq m^2, \|z'\|_p \leq m, \|z'\|_{\sup} \leq \|f\|_{\sup}\}$$

defined on the time interval $[0,T]$.

**Lemma** – *There is a function $\epsilon : [1, \infty) \to (0, \infty)$ such that the ball $\mathcal{B}_T(m)$ is invariant by the map $\mathbb{I}$ if $\omega_{\boldsymbol{X}}(0,T) \leq \epsilon(m)$.*

**Proof** – For $\boldsymbol{z} \in \mathcal{B}_T(m)$ one has $\|f(z)\|_p \leq \|f\|_{\operatorname{Lip}} \|z\|_p \leq m$ and $\|f(z)\|_{\sup} \leq \|f\|_{\sup}$. Write here $\mathbb{I}(\boldsymbol{z}) = (\varphi, f(z))$. We have from Lemma 26

$$\|f'(z)z'\|_p \lesssim_f m \text{ and } \|R^{f(\boldsymbol{z})}\|_{p/2} \lesssim_f m^2,$$

so we have $\|R^{\varphi}\|_{p/2} \lesssim_f \omega_{\boldsymbol{X}}(1+m^2)$ from Corollary 25, which entails that $\|\varphi\|_p \lesssim_f \omega_{\boldsymbol{X}}(1+m^2)$. $\square$

So, given a rough path $\boldsymbol{X}$ and $m \geq 1$, both fixed, we can choose $T = T(\boldsymbol{X},m)$ to have $\omega_{\boldsymbol{X}}(0, T) \leq \epsilon(m)$, as $\omega_{\boldsymbol{X}}$ is a continuous function of $T \geq 0$ with limit $0$ in $0$. We now work on that time interval $[0, T]$. The set

$$\mathcal{B}_T^0(m) = \mathcal{B}_T(m) \cap \{\boldsymbol{z} = (z, z'); z_0 = x_0, z_0' = f(x_0)\}$$

is turned into a Polish space by the metric

$$d(\boldsymbol{z}^1, \boldsymbol{z}^2) \equiv \max\left(\|\Delta z\|_p, \|\Delta z'\|_p, \|\Delta R^z\|_{p/2}\right).$$

A direct use of the stability results from Section 20.1 then gives the following fact.

**Lemma (Contraction)** – *A possible redefinition of the map $\epsilon(\cdot)$ turns the map $\mathbb{I}$ into a contraction of $(\mathcal{B}_T^0(m), d)$.*

The map $\mathbb{I}$ thus has a unique fixed point in $\mathcal{B}_T^0(m)$. It actually has a unique fixed point in the set of all controlled path $\boldsymbol{z}$ with $z_0 = x_0$ and $z_0' = f(x_0)$ from some elementary classical reasoning. It is customary to say from the fixed point of the map $\mathbb{I}$ that it is the solution of the rough differential equation

$$z_t = x_0 + \int_0^t f(\boldsymbol{z}_t) d\boldsymbol{X}_t$$

Finding a (unique) solution to this equation is the equivalent here of the problem of finding a (unique) solution of the equation

$$\boldsymbol{u} = \boldsymbol{\mathcal{F}}(u_0) + \mathcal{K}^M(f(\boldsymbol{u}) \circ + g(\boldsymbol{u}, \partial \boldsymbol{u}))$$

---

56. See Section 4.3.3 of [6] for these statements and their proofs.

on a space of modelled distribution.

Things get a different flavour when the driving rough path is random. You can consult Chapter 5 of [1] to see how one can relate the classical Itô calculus to rough paths calculus, and for some elementary applications of rough paths to the analysis of stochastic differential equations. The benefits of the rough path approach compared to the Itô setting are similar to those that we pointed out in the more general setting of singular SPDEs. The continuity of the solution map $\boldsymbol{X} \mapsto z$ allow to transfer automatically a number of probabilistic properties of $\boldsymbol{X}$ to the random path $z$. If we have for instance a support theorem[57] for $\boldsymbol{X}$ we will have a support theorem for $z$. If we have a large deviation theorem for $\boldsymbol{X}$ we will have a large deviation theorem for $z$, *etc*.

No renormalization is needed when we lift Brownian motion or some other classes of reasonable random processes into some rough paths. However on can devilishly construct some smooth random processes $h^\epsilon$ whose canonical lift does not converge in a probabilistic sense and for which one needs renormalization.

**References**

[1] I. Bailleul, *A flow-based approach to rough differential equations*, http://lmba.math.univ-brest.fr/perso/ismael.bailleul/Files/M2Course.pdf Some lecture notes on rough differential equations with an emphasis on solution flows. Contains a sewing lemma for approximate flows of independent interest.

[2] D. Feyel and A. de la Pradelle, *Curvilinear Integrals Along Enriched Paths*. Elec. J. Probab. **11**:860–892 (2006). The authors Simplified Gubinelli's approach in [4] and gave what is now canonical proof of the sewing lemma.

[3] D. Feyel and A. de la Pradelle, *Researches in Stochastic Analysis*, https://spartacus-idh.com/pdfs/026/ A fantastic book of deep interst on stochastic analysis. A jewell!

[4] M. Gubinelli, *Controlling rough paths*. J. Funct. Anal. **216**:86–140 (2004). The work where the notions of sewing lemma, path controlled by another path and rough integral were introduced first.

[5] T. Lyons, *Differential equations driven by rough signals*, Rev. Mat. Iberoamericana **14**(2) (1998). The groundbreaking work on rough paths.

[6] P. Zorin-Kranich, *Lecture notes on martingale inequalities*. arXiv:2404.17197 (2024). A super nice, up-to-date, set of lectures on martingale inequalities.

---

57. That is an exact description of the support of the law of $\boldsymbol{X}$.

# *Solutions to the exercices*

**Lecture 2**

*Exercice 1* – If $\Lambda_x(y)$ is a continuous function of $(x,y)$ then $\mathbb{I}_s^t(x) \xrightarrow[s\downarrow 0]{} p_t(x,\Lambda_x) \xrightarrow[t\downarrow 0]{} \Lambda_x(x)$.

*Exercice 2* – It suffices to check that both maps coincide on $C_b(\mathbb{R}^d)$ and that they are both continuous on $C^{r_1} \times C^{r_2}$. The continuity of the Littlewood-Paley definition is proved in Appendix 1. To see the continuity of the Young product operation we need see first that

$$(f,g) \in C^{r_1} \times C^{r_2} \mapsto (\Lambda_x)_x$$

is continuous and then use the continuity of the reconstruction map. See Theorem 12.7 in [4].

**Lecture 3**

*Exercice 1* – Recall that in the shorthand notation $\Delta(\tau) = \sum_{\sigma \leq \tau} \sigma \otimes (\tau/\sigma)$ the $\sigma$ runs over a certain $\tau$-dependent subset of the vectors of the linear basis $\boldsymbol{B}$. For $\tau \in T$ one has

$$(\Delta \otimes \mathrm{Id})\Delta(\tau) = (\mathrm{Id} \otimes \Delta^+)\Delta(\tau).$$

Writing $\Delta(\tau) = \sum_{\sigma \leq \tau} \sigma \otimes (\tau/\sigma)$, and writing $\Delta(\sigma) = \sum_{\eta \leq \sigma}$ for each of these $\eta$ one has from the left hand side of the identity

$$(\Delta \otimes \mathrm{Id})\Delta(\tau) = \sum_{\eta \leq \sigma \leq \tau} \eta \otimes (\sigma/\eta) \otimes (\tau/\sigma).$$

Change notation and write $\Delta(\tau) = \sum_{\eta \leq \tau} \eta \otimes (\tau/\eta)$ so

$$(\mathrm{Id} \otimes \Delta^+)\Delta(\tau) = \sum_{\eta \leq \tau} \eta \otimes \Delta^+(\tau/\eta).$$

Since an element in the finite dimensional tensor product $T \otimes T^+$ has a unique decomposition over the set of elements of the form $\sum_{\eta \in \boldsymbol{B}} \eta \otimes (\ldots)$ we obtain the identity

$$\Delta^+(\tau/\eta) = \sum_{\eta \leq \sigma \leq \tau} (\sigma/\eta) \otimes (\tau/\sigma)$$

here a sum over the $\sigma$.

*Exercice 2* – The results of this exercice are used in Section 12 of Lecture 4. As a shorthand notation we write $\Delta^+(\sigma) = \sum \sigma_1 \otimes \sigma_2$, for any $\sigma \in T^+$.

*(a)* The antipode $S: T^+ \to T^+$ is characterized by the property that $S(\boldsymbol{1}_+) = \boldsymbol{1}_+$ and one has for any $\sigma \in T^+ \setminus \{\boldsymbol{1}_+\}$

$$\sum \sigma_1 S(\sigma_2) = 0$$

with the above notation for $\Delta^+(\sigma)$. First, as $\Delta^+(\boldsymbol{1}_+) = \boldsymbol{1}_+ \otimes \boldsymbol{1}_+$ one has

$$(g \star (g \circ S))(\boldsymbol{1}_+) = g(\boldsymbol{1}_+) g(S(\boldsymbol{1}_+)) = 1$$

as $g(\boldsymbol{1}_+) = 1$ for all characters. Second, because $g$ is multiplicative and linear, one has for $\sigma \neq \boldsymbol{1}_+$

$$\begin{aligned}
(g \star (g \circ S))(\sigma) &= (g \otimes (g \circ S))\Delta^+(\sigma) \\
&= (g \otimes (g \circ S))\left(\sum \sigma_1 \otimes \sigma_2\right) \\
&= \sum g(\sigma_1) g(S(\sigma_2)) \\
&= \sum g(\sigma_1 S(\sigma_2)) \\
&= g\left(\sum \sigma_1 S(\sigma_2)\right) = g(0) = 0
\end{aligned}$$

*(b)* Recall the identity $(\mathrm{Id} \otimes \Delta^+)\Delta = (\Delta \otimes \mathrm{Id})\Delta$. For $\tau \in T^+$ one has

$$\begin{aligned}
(\widehat{g_1 \star g_2})(\sigma) &= (\mathrm{Id} \otimes \{g_1 \star g_2\})\Delta(\tau) \\
&= (\mathrm{Id} \otimes \{(g_1 \otimes g_2)\Delta^+\})\Delta(\tau) \\
&= (\mathrm{Id} \otimes g_1 \otimes g_2)(\mathrm{Id} \otimes \Delta^+)\Delta(\tau) \\
&= (\mathrm{Id} \otimes g_1 \otimes g_2)(\Delta \otimes \mathrm{Id})\Delta(\tau) \\
&= (\{(\mathrm{Id} \otimes g_1)\Delta\} \otimes g_2)\Delta(\tau) \\
&= ((\mathrm{Id} \otimes g_1)\Delta)(\mathrm{Id} \otimes g_2)\Delta(\tau) \\
&= \hat{g}_1(\hat{g}_2(\tau))
\end{aligned}$$

*Exercice 3* – We work on the polynomial regularity structure of example **(a)** in Section 7. We need to check the $x, y, t$ uniform bounds

$$|\langle \Pi_x(n), p_t(x, \cdot)\rangle| \lesssim t^{n/2}$$

and

$$|g_{yx}(n)| \lesssim |y-x|^{|n|}.$$

The first comes from the bound (1) in Lecture 2 about the Gaussian moments; the second is obvious.

*Exercice 4* – A direct consequence of the Gaussian moment estimates: bound (6) in Lecture 2.

*Exercice 5* – Let me start with two preliminary remarks. First, since we assume that the product of any two symbols $\tau_1, \tau_2$ of $T$ of positive sizes is an element of $T$, the multiplicativity of the map $\Delta$ implies that the maps $\Gamma_{yx}$ satisfy, with $\Delta(\tau_i) = a_i \otimes b_i$, that $\Gamma_{yx}(\tau_1 \tau_2)$ is equal to

$$(\mathrm{Id} \otimes g_{yx})\Delta(\tau_1 \tau_2) = (\mathrm{Id} \otimes g_{yx})\Delta(\tau_1)\Delta(\tau_2) = (\mathrm{Id} \otimes g_{yx})((a_1 a_2) \otimes (b_1 b_2)) = a_1 a_2 g(b_1 b_2) = \Gamma_{yx}(\tau_1)\Gamma_{yx}(\tau_2).$$

Second, please check that the fact that

$$\Delta(T_a) \subset \bigoplus_{\alpha \in A^+} T_{a-\alpha} \otimes T_\alpha^+$$

for all $a \in A$ implies that

$$\Gamma_{yx} Q_\gamma = Q_\gamma \Gamma_{yx}.$$

Both remarks hold for all $x, y$. Now, the expression

$$F(\boldsymbol{v})(x) \equiv Q_\gamma \left( \sum_{n \geq 0} \frac{F^{(n)}(v_{\mathbf{1}})}{n!} (\boldsymbol{v}(x) - v_{\mathbf{1}}(x)\mathbf{1})^n \right)$$

only keeps from the formal Taylor series of $F(\boldsymbol{v})$ around the point $\boldsymbol{v}_1(x)$ its components of size smaller than $\gamma$, so the sum is actually finite. Here is a sketch of the argument; I leave you fill in the blanks. The details are fully spelled out in Hairer's 2014 work, Theorem 4.16 therein. Set $\overline{\boldsymbol{v}}(x) \equiv \boldsymbol{v}(x) - v_{\mathbf{1}}(x)\mathbf{1}$.

The index set $A$ of the regularity intersects any bounded set into a discrete set; we define $\gamma_0 \equiv \min(a \in A; a \neq 0) > 0$.

**1.** We have

$$\Gamma_{yx}(F(\boldsymbol{v})(x)) = \sum_{n \leq \gamma/\gamma_0} \frac{F^{(n)}(v_{\mathbf{1}}(x))}{n!} \{\Gamma_{yx}(\overline{\boldsymbol{v}}(x))\}^n + \sharp_{yx}$$

with $|(\sharp_{yx})_\tau| \lesssim |y-x|^{\gamma - |\tau|}$ for all $\tau \in \boldsymbol{B}$ with $|\tau| < \gamma$. All the remainders $\sharp_{yx}^1, \sharp_{yx}^2, \sharp_{yx}^3$ have the same property.

**2.** Write

$$\Gamma_{yx}(\overline{\boldsymbol{v}}(x)) = \Gamma_{yx}(\boldsymbol{v}(x)) - v_{\mathbf{1}}(x)\mathbf{1} = \boldsymbol{v}(y) + \sharp'_{yx} - v_{\mathbf{1}}(x)\mathbf{1} = \overline{\boldsymbol{v}}(y) + \sharp^1_{yx} + (v_{\mathbf{1}}(y) - v_{\mathbf{1}}(x))\mathbf{1}.$$

**3.** These two bounds together give

$$\Gamma_{yx}(F(\boldsymbol{v})(x)) = \sum_{n \leq \gamma/\gamma_0} \frac{F^{(n)}(v_{\mathbf{1}}(x))}{n!} \{\overline{\boldsymbol{v}}(y) + (v_{\mathbf{1}}(y) - v_{\mathbf{1}}(x))\mathbf{1}\}^n + +\sharp^2_{yx}. \tag{69}$$

By definition of the $\mathcal{D}^\gamma$ spaces, and given that $\boldsymbol{v}$ has no components on the trees of negative size, the function $v_{\mathbf{1}}$ is $\gamma_0$-Hölder. It follows that

$$|(\{\overline{\boldsymbol{v}}(y) + (v_{\mathbf{1}}(y) - v_{\mathbf{1}}(x))\mathbf{1}\}^n)_\tau| \lesssim |y - x|^{n\gamma_0 - |\tau|}$$

**4.** Combining the Taylor expansion

$$F^{(n)}(v_{\mathbf{1}}(x)) = \sum_{|n+\ell| < \gamma/\gamma_0} \frac{F^{(n+\ell)}(v_{\mathbf{1}}(y))}{\ell!} (v_{\mathbf{1}}(y) - v_{\mathbf{1}}(x))^\ell + O(|y - x|^{\gamma - n\gamma_0})$$

with (1) gives

$$\Gamma_{yx}(F(\boldsymbol{v})(x)) = \sum_{n \leq \gamma/\gamma_0} \frac{F^{(n+\ell)}(v_{\mathbf{1}}(y))}{n!\ell!} (v_{\mathbf{1}}(y) - v_{\mathbf{1}}(x))^\ell \{\overline{\boldsymbol{v}}(y) + (v_{\mathbf{1}}(y) - v_{\mathbf{1}}(x))\mathbf{1}\}^n + \sharp^3_{yx}.$$

We conclude with the identity

$$\sum_{n+\ell=m} \frac{(v_{\mathbf{1}}(y) - v_{\mathbf{1}}(x))^\ell \{\overline{\boldsymbol{v}}(y) + (v_{\mathbf{1}}(y) - v_{\mathbf{1}}(x))\mathbf{1}\}^n}{n!\ell!} = \frac{\overline{\boldsymbol{v}}(y)^m}{m!}.$$

**Lecture 4**

*Exercice 1* – A direct application of the definition of $\Gamma_{yx} = (\mathrm{Id} \otimes g_{yx})\Delta$ and the identity (40) in Lecture 4. Note that $|\mathcal{I}(\circ)| = 2 + |\circ|$.

*Exercice 2* – We can see that preparation maps are some perturbations of the identity by a nilpotent map: $R = \mathrm{Id} + N$. The inverse of such a map is given by the a posteriori finite sum $\sum (-N)^a$, where $a \in \mathbb{N}$. For $R(\tau) = \tau + \sum \lambda_i \tau_i$, since the number of noises in each $\tau_i$ is strictly smaller than the number of noise in $\tau$, one sees that $N$ sends the linear space spanned by symbols with $b$ noise symbols into the direct sum of the linear spaces spanned by symbols with at most $b-1$ noise symbols. The result follows.